\documentclass[12pt]{article}
\usepackage{amsthm,amsfonts,amssymb,epsfig,graphics,amsmath,amsbsy,color,enumerate,tikz}

\usepackage{calligra}

\DeclareMathAlphabet{\mathcalligra}{T1}{calligra}{m}{n}
\DeclareFontShape{T1}{calligra}{m}{n}{<->s*[2.2]callig15}{}

\theoremstyle{plain}
\newtheorem{theorem}{Theorem}
\newtheorem{lemma}{Lemma}

\newtheorem{remark}{Remark}
\newtheorem{definition}{Definition}

\newcommand{\mas}{\operatorname{Mas}}

\newcommand{\loc}{\operatorname{loc}}

\newcommand{\AC}{\operatorname{AC}}

\newcommand{\ess}{\operatorname{ess}}
\newcommand{\diag}{\operatorname{diag}}

\numberwithin{equation}{section}
\numberwithin{lemma}{section}
\numberwithin{theorem}{section}
\numberwithin{remark}{section}
\numberwithin{claim}{section}
\numberwithin{corollary}{section}
\numberwithin{proposition}{section}
\numberwithin{definition}{section}
\numberwithin{condition}{section}
\numberwithin{figure}{section}

\title{The Maslov index and spectral counts for 
linear Hamiltonian systems on $\mathbb{R}$}

\author{Peter Howard}

\begin{document}

\maketitle

\begin{abstract} 
Working with a general class of linear Hamiltonian systems 
specified on $\mathbb{R}$, we develop a framework 
for relating the Maslov index to the number of eigenvalues
the systems have on intervals of the form $[\lambda_1, \lambda_2)$
and $(-\infty, \lambda_2)$. We verify that our framework can 
be implemented for Sturm-Liouville systems, fourth-order potential
systems, and a family of systems nonlinear in the spectral
parameter. The analysis is primarily motivated by applications to 
the analysis of spectral stability for nonlinear 
waves, and aspects of such analyses are emphasized. 
\end{abstract}

\section{Introduction}\label{introduction}

For values $\lambda$ confined to an interval $I \subset \mathbb{R}$, 
we consider linear Hamiltonian systems 
\begin{equation} \label{hammy}
J y' = \mathbb{B} (x; \lambda) y; \quad x \in \mathbb{R}, 
\quad y(x; \lambda) \in \mathbb{C}^{2n},
\end{equation}
where $J$ denotes the symplectic matrix 
\begin{equation*}
J 
=
\begin{pmatrix}
0_n & - I_n \\
I_n & 0_n
\end{pmatrix},
\end{equation*}
and throughout the analysis we will make the following 
assumptions on $\mathbb{B} (x; \lambda)$:

\medskip
{\bf (A)} For each $\lambda \in I$, $\mathbb{B} (\cdot; \lambda) \in L^1_{\loc} (\mathbb{R}, \mathbb{C}^{2n \times 2n})$,
with $\mathbb{B} (x; \lambda)$ self-adjoint for a.e. $x \in \mathbb{R}$, and additionally
the partial derivatives $\mathbb{B}_{\lambda} (x; \lambda)$ exist for a.e. $x \in \mathbb{R}$,
with $\mathbb{B}_{\lambda} (\cdot; \lambda) \in L^1_{\loc} (\mathbb{R}, \mathbb{C}^{2n \times 2n})$.
\medskip

For this analysis, we will say that $\lambda$ is an eigenvalue of 
(\ref{hammy}) provided there exists a function
\begin{equation*}
    y (\cdot; \lambda) \in (\AC_{\loc} (\mathbb{R}, \mathbb{C}^{2n})
    \cap L^2 (\mathbb{R}, \mathbb{C}^{2n})) \backslash \{0\}
\end{equation*}
that satisfies (\ref{hammy}) for a.e. $x \in \mathbb{R}$, and 
we will take the geometric multiplicity of $\lambda$ to be the 
dimension of the space of such solutions. (Here, 
$\AC_{\loc} (\cdot)$ refers to the space of functions absolutely 
continuous on compact subsets of $\mathbb{R}$.) Our primary goal 
for the analysis is to use the Maslov index to count the number
of eigenvalues that (\ref{hammy}) has on intervals of the form 
$[\lambda_1, \lambda_2)$ and $(- \infty, \lambda_2)$ (assumed, 
in each case, to be a subset of $I$).

We are primarily motivated by applications to the spectral 
stability of nonlinear waves arising in certain nonlinear 
evolutionary PDE such as Allen-Cahn systems
\begin{equation} \label{allen-cahn-equation}
    u_t + DF(u) = u_{xx},
    \quad (x,t) \in \mathbb{R} \times \mathbb{R}_+,
    \quad u(x, t) \in \mathbb{C}^n,
\end{equation}
and higher-order analogues 
\begin{equation*}
    u_t + DF(u) = - u_{xxxx}, 
    \quad (x,t) \in \mathbb{R} \times \mathbb{R}_+,
    \quad u(x, t) \in \mathbb{C}^n.
\end{equation*}
(Here, $D$ denotes the Jacobian operator.)
In the former case, if $\bar{u} (x)$ denotes a stationary 
solution, then we can linearize about $\bar{u} (x)$
(setting $u = \bar{u} + v$ and dropping terms nonlinear in 
$v$) to obtain a linear equation 
\begin{equation*}
    v_t + D^2 F(\bar{u}) v = v_{xx},
    \quad (x,t) \in \mathbb{R} \times \mathbb{R}_+,
    \quad v(x, t) \in \mathbb{C}^n.
\end{equation*}
In this setting, spectral stability is determined 
by the spectrum of the associated eigenvalue problem 
\begin{equation} \label{schrodinger-intro}
- \phi_{xx} + D^2F (\bar{u}) \phi = \lambda \phi,
\quad x \in \mathbb{R}, 
\quad \phi (x) \in \mathbb{C}^n, 
\end{equation}
which we can put in the form 
of (\ref{hammy}) by setting $y = {y_1 \choose y_2} = {\phi \choose \phi'}$.
Precisely, we find  
\begin{equation*}
    Jy' = \mathbb{B} (x; \lambda) y; 
    \quad 
    \mathbb{B} (x; \lambda) = 
    \begin{pmatrix}
    \lambda I - D^2F (\bar{u} (x)) & 0 \\
    0 & I
    \end{pmatrix}.
\end{equation*}
We would like to determine whether (\ref{schrodinger-intro})
has any negative eigenvalues, and this information clearly 
follows from a count of the number of eigenvalues
that  (\ref{schrodinger-intro}) has on $(- \infty, 0)$.

We are particularly interested in stationary solutions
that approach fixed endstates $u_{\pm}$ as $x \to \pm \infty$, 
and such cases provide us with additional structure 
that will be necessary for our general analysis. In 
order to keep the analysis as applicable as possible, 
we will make three general assumptions (in addition to 
Assumptions {\bf (A)}), which we will 
subsequently verify in a selection of important cases.
Prior to stating these assumptions, we need to 
develop some notation and terminology that will be 
used througout the discussion. We begin with the 
following definition.

\begin{definition} \label{left-right-definition}
We say that a solution $y (\cdot;\lambda) \in \AC_{\loc} (\mathbb{R},\mathbb{C}^{2n})$
of (\ref{hammy}) {\it lies left} in $\mathbb{R}$ if for any $c \in \mathbb{R}$, 
the restriction of $y(\cdot; \lambda)$ to $(-\infty,c)$ is in 
$L^2 ((-\infty, c), \mathbb{C}^{2n})$. Likewise, 
we say that a solution $y (\cdot;\lambda) \in \AC_{\loc} (\mathbb{R},\mathbb{C}^{2n})$
of (\ref{hammy}) {\it lies right} in $\mathbb{R}$ if for any $c \in \mathbb{R}$, 
the restriction of $y(\cdot; \lambda)$ to $(c,+\infty)$ is in 
$L^2 ((c,+\infty), \mathbb{C}^{2n})$. 
\end{definition}

Our primary tool for this analysis will be the 
Maslov index, and as a starting point for a discussion of this object, we 
define what we will mean by a Lagrangian subspace of $\mathbb{C}^{2n}$. 

\begin{definition} \label{lagrangian_subspace}
We say $\ell \subset \mathbb{C}^{2n}$ is a Lagrangian subspace of $\mathbb{C}^{2n}$
if $\ell$ has dimension $n$ and
\begin{equation} 
(J u, v) = 0, 
\end{equation} 
for all $u, v \in \ell$. (Here, and throughout, $(\cdot, \cdot)$ 
denotes the usual inner product on $\mathbb{C}^{2n}$.) In addition, we denote by 
$\Lambda (n)$ the collection of all Lagrangian subspaces of $\mathbb{C}^{2n}$, 
and we will refer to this as the {\it Lagrangian Grassmannian}. 
\end{definition}

Any Lagrangian subspace of $\mathbb{C}^{2n}$ can be
spanned by a choice of $n$ linearly independent vectors in 
$\mathbb{C}^{2n}$. We will generally find it convenient to collect
these $n$ vectors as the columns of a $2n \times n$ matrix $\mathbf{X}$, 
which we will refer to as a {\it frame} for $\ell$. Moreover, we will 
often coordinatize our frames as $\mathbf{X} = {X \choose Y}$, where $X$ and $Y$ are 
$n \times n$ matrices. Following \cite{F2004} (p. 274), we specify 
a metric on $\Lambda (n)$ in terms of appropriate orthogonal projections. 
Precisely, let $\mathcal{P}_i$ 
denote the orthogonal projection matrix onto $\ell_i \in \Lambda (n)$
for $i = 1,2$. I.e., if $\mathbf{X}_i$ denotes a frame for $\ell_i$,
then $\mathcal{P}_i = \mathbf{X}_i (\mathbf{X}_i^* \mathbf{X}_i)^{-1} \mathbf{X}_i^*$.
We take our metric $d$ on $\Lambda (n)$ to be defined 
by 
\begin{equation*}
d (\ell_1, \ell_2) := \|\mathcal{P}_1 - \mathcal{P}_2 \|,
\end{equation*} 
where $\| \cdot \|$ can denote any matrix norm. We will say 
that a path of Lagrangian subspaces 
$\ell: \mathcal{I} \to \Lambda (n)$ is continuous provided it is 
continuous under the metric $d$. 

Suppose $\ell_1 (\cdot), \ell_2 (\cdot)$ denote continuous paths of Lagrangian 
subspaces $\ell_i: \mathcal{I} \to \Lambda (n)$, $i = 1,2$, for some parameter interval 
$\mathcal{I}$. The Maslov index associated with these paths, which we will 
denote $\mas (\ell_1, \ell_2; \mathcal{I})$, is a count of the number of times
the paths $\ell_1 (\cdot)$ and $\ell_2 (\cdot)$ intersect, counted
with both multiplicity and direction. (In this setting, if we let 
$t_*$ denote the point of intersection (often referred to as a 
{\it conjugate point}), then multiplicity corresponds with the dimension 
of the intersection $\ell_1 (t_*) \cap \ell_2 (t_*)$; a precise definition of what we 
mean in this context by {\it direction} will be
given in Section \ref{maslov-section}.) 

We are now prepared to state the three general assumptions (in addition 
to Assumptions {\bf (A)}) that will be required for our analysis. For 
convenient reference, some notational conventions will be embedded 
in the statements of these assumptions. 

\medskip
{\bf (B1)} For each $\lambda \in I$, there exists an $n$-dimensional space
of solutions to (\ref{hammy}) that lie left in $\mathbb{R}$, and likewise
an $n$-dimensional space of solutions to (\ref{hammy}) that lie right 
in $\mathbb{R}$. We will denote by $\mathbf{X} (x; \lambda)$
a frame comprising a choice of basis for the $n$-dimensional space
of solutions to (\ref{hammy}) that lie left in $\mathbb{R}$, and 
we will denote by $\tilde{\mathbf{X}} (x; \lambda)$
a frame comprising a choice of basis for the $n$-dimensional space
of solutions to (\ref{hammy}) that lie right in $\mathbb{R}$. 
We will show that when constructed in this way, $\mathbf{X} (x; \lambda)$
and $\tilde{\mathbf{X}} (x; \lambda)$ constitute frames for 
Lagrangian subspaces of $\mathbb{C}^{2n}$, which we will 
respectively denote $\ell (x; \lambda)$ 
and $\tilde{\ell} (x; \lambda)$. We assume that  
$\ell, \tilde{\ell} \in C (\mathbb{R} \times I, \Lambda (n))$,
with additionally 
$\mathbf{X} (x; \cdot), \tilde{\mathbf{X}} (x; \cdot) \in \AC_{\loc} (I, \mathbb{C}^{2n \times n})$
for all $x \in \mathbb{R}$. 
\medskip

\medskip
{\bf (B2)} For each $\lambda \in I$, the asymptotic frames
\begin{equation*}
    \mathbf{X}_{\pm} (\lambda) := \lim_{x \to \pm \infty} \mathbf{X} (x; \lambda); 
    \quad {\rm and} \quad
    \tilde{\mathbf{X}}_+ (\lambda) := \lim_{x \to +\infty} \tilde{\mathbf{X}} (x; \lambda)
\end{equation*}
are well-defined, and are respectively frames for Lagrangian subspaces 
$\ell_{\pm} (\lambda)$ and $\tilde{\ell}_+ (\lambda)$. In addition,
\begin{equation*}
    \ell_- (\lambda) \cap \tilde{\ell}_+ (\lambda) = \{0\}
    \quad \forall \, \lambda \in I.
\end{equation*}

\medskip
{\bf (B3)} There exists a constant $c_0 > 0$ sufficiently large so that 
for any $c > c_0$, the matrix 
\begin{equation*}
    \int_{-\infty}^c \mathbf{X} (x; \lambda)^* \mathbb{B}_{\lambda} (x; \lambda) \mathbf{X} (x; \lambda) dx
\end{equation*}
is positive definite for all $\lambda \in I$. 
\medskip

Assumptions {\bf (B1)}, {\bf (B2)}, and {\bf (B3)}, along with 
Assumptions {\bf (A)}, hold in many important cases. As specific
examples, we will verify them for linear Hamiltonian 
systems associated with Sturm-Liouville Systems
\begin{equation} \label{sls}
    - (P(x) \phi')' + V(x) \phi = \lambda Q (x) \phi,
    \quad x \in \mathbb{R}, \quad \phi (x) \in \mathbb{C}^n,
\end{equation}
fourth-order potential equations, 
\begin{equation} \label{fourth-ode}
    \phi'''' + V(x) \phi = \lambda \phi,
     \quad x \in \mathbb{R}, \quad \phi (x) \in \mathbb{C}^n,
\end{equation}
and a family of systems nonlinear in the spectral 
parameter $\lambda$,
\begin{equation} \label{sls-a}
    - (P_{11} (x) \phi')' + V_{11} (x) \phi + V_{12} (x) (\lambda I - V_{22} (x))^{-1} V_{12} (x)^* \phi
    = \lambda \phi, 
     \quad x \in \mathbb{R}, \quad \phi (x) \in \mathbb{C}^n,
\end{equation}
with appropriate assumptions on the coefficient matrices in all cases. (Equation
(\ref{sls-a}) arises in the analysis of differential-algebraic Sturm-Liouville
systems; see Section \ref{sls-a-section} for details.)

We can state our main theorem as follows. 

\begin{theorem} \label{main-theorem} Let Assumptions 
{\bf (A)}, {\bf (B1)}, {\bf (B2)}, and {\bf (B3)}
hold. If $\mathcal{N} ([\lambda_1, \lambda_2))$ denotes
the number of eigenvalues that (\ref{hammy}) has 
on the interval $[\lambda_1, \lambda_2)$, counted with 
geometric multiplicity, then 
\begin{equation*}
    \mathcal{N} ([\lambda_1, \lambda_2))
    = - \mas (\ell (\cdot; \lambda_2), \tilde{\ell}_+ (\lambda_2), (-\infty, + \infty])
    + \mas (\ell (\cdot; \lambda_1), \tilde{\ell}_+ (\lambda_1), (-\infty, + \infty]).
\end{equation*}
\end{theorem}

\begin{remark} \label{brackets}
The inclusive bracket on $+ \infty$ indicates that we use Assumption 
{\bf (B2)} to compactify $\mathbb{R}$ for our Maslov index calculations. 
In particular, this means that, at least in principle, $\pm \infty$ 
can serve as conjugate points. According to Assumption {\bf (B2)}, we have 
$\ell_- (\lambda_1) \cap \tilde{\ell}_+ (\lambda_1) = \{0\}$, so $- \infty$
will never serve as a conjugate point for our analysis (hence the open 
parentheses on $-\infty$), but it may be the case that  
$\ell_+ (\lambda_1) \cap \tilde{\ell}_+ (\lambda_1) \ne \{0\}$,
in which case $+ \infty$ will serve as a conjugate point. 
The same remark holds if $\lambda_1$ is replaced by $\lambda_2$. 
\end{remark}

For specific applications such as the ones we will discuss in 
detail, we can establish additional properties that may not 
hold in the generality Theorem \ref{main-theorem}. 
Among these, we will emphasize the following: 

\begin{enumerate}

\item In some cases we can replace the target spaces 
$\tilde{\ell}_+ (\lambda_1)$ and $\tilde{\ell}_+ (\lambda_2)$ with 
target spaces for which the flow associated with the relevant 
Maslov index is monotonic (i.e., the sign associated with each 
conjugate point is the same). As an important example, we will 
show in Section \ref{sls-section} that for 
Sturm-Liouville systems, this is the case for 
the Dirichlet Lagrangian subspace $\ell_D$ 
(with frame $\mathbf{X}_D = {0 \choose I}$). 

\item As discussed in our motivating applications to 
the stability of nonlinear waves, we are often interested in counting 
the number of eigenvalues that (\ref{hammy}) has below 
some fixed value $\lambda_2$, and for this it's convenient 
to show that we can take $\lambda_1$ sufficiently negative
so that 
\begin{equation*}
    \mas (\ell (\cdot; \lambda_1), \tilde{\ell}_+ (\lambda_1), (-\infty, + \infty])
    = 0.
\end{equation*}
In this case, 
\begin{equation*}
    \mathcal{N} ((-\infty, \lambda_2))
    = - \mas (\ell (\cdot; \lambda_2), \tilde{\ell}_+ (\lambda_2), (-\infty, + \infty]).
\end{equation*}

\item In certain specialized cases, we can apply our results 
to operators that are not self-adjoint. The most important 
such case arises when (\ref{allen-cahn-equation}) is linearized 
about a traveling wave solution $\bar{u} (x - st)$, leading
to the eigenvalue problem 
\begin{equation} \label{schrodinger-intro-shifted}
- \phi_{xx} - s \phi_x + D^2 F(\bar{u}) \phi = \lambda \phi,
\quad x \in \mathbb{R}, \quad \phi (x) \in \mathbb{C}^n.
\end{equation}
This case will be discussed in Section \ref{traveling-waves-section}.

\end{enumerate}

We now state specific results obtained for (\ref{sls}), (\ref{fourth-ode}),
and (\ref{sls-a}). In all cases, we refer to later sections, where 
detailed assumptions are stated.

\begin{theorem} \label{sls-theorem}
For (\ref{sls}), let Assumptions {\bf (SL1)} and {\bf (SL2)} 
from Section \ref{sls-section} hold, and express (\ref{sls})
in the form (\ref{hammy}) (giving (\ref{sls-hammy})). Then
for $\kappa$ specified as in (\ref{kappa-defined}), {\bf (A)},
{\bf (B1)}, {\bf (B2)}, and {\bf (B3)} all hold for 
(\ref{sls-hammy}) with $I = (- \infty, \kappa)$, and so 
the result of Theorem \ref{main-theorem} holds for 
all intervals $[\lambda_1, \lambda_2]$, $\lambda_1 < \lambda_2 < \kappa$. 
In addition, if $\mathcal{N} ([\lambda_1, \lambda_2))$
denotes the number of eigenvalues, counted with multiplicity,
that (\ref{sls}) has on the interval $[\lambda_1, \lambda_2)$, 
and we express the frame $\mathbf{X} (x; \lambda)$ from 
{\bf (B1)} as $\mathbf{X} (x; \lambda) = {X (x; \lambda) \choose Y (x; \lambda)}$,
then we have 
\begin{equation*}
    \mathcal{N} ([\lambda_1, \lambda_2))
    = \sum_{x \in \mathbb{R}} \dim \ker X(x; \lambda_2)
    - \sum_{x \in \mathbb{R}} \dim \ker X(x; \lambda_1),
\end{equation*}
and 
\begin{equation*}
    \mathcal{N} ((-\infty, \lambda_2)) 
    = \sum_{x \in \mathbb{R}} \dim \ker X(x; \lambda_2). 
\end{equation*}
\end{theorem}

\begin{theorem} \label{sls-a-theorem}
For (\ref{sls-a}), let Assumptions {\bf (DA1)} and {\bf (DA2)} 
from Section \ref{sls-a-section} hold, and express (\ref{sls-a})
in the form (\ref{hammy}) (giving (\ref{sls-a-hammy})). Then
for any interval $I \subset \mathbb{R}$ satisfying 
(\ref{sls-a-i}), Assumptions {\bf (A)},
{\bf (B1)}, {\bf (B2)}, and {\bf (B3)} all hold for 
(\ref{sls-a-hammy}), and so 
the result of Theorem \ref{main-theorem} holds for 
all intervals $[\lambda_1, \lambda_2] \subset I$. 
In addition, if $\mathcal{N} ([\lambda_1, \lambda_2))$
denotes the number of eigenvalues, counted with multiplicity,
that (\ref{sls-a}) has on the interval $[\lambda_1, \lambda_2)$,
and we express the frame $\mathbf{X} (x; \lambda)$ from 
{\bf (B1)} as $\mathbf{X} (x; \lambda) = {X (x; \lambda) \choose Y (x; \lambda)}$,
then we have 
\begin{equation*}
    \mathcal{N} ([\lambda_1, \lambda_2))
    = \sum_{x \in \mathbb{R}} \dim \ker X(x; \lambda_2)
    - \sum_{x \in \mathbb{R}} \dim \ker X(x; \lambda_1).
\end{equation*}
Finally, if $\lambda_2 \in I$ lies entirely below 
$\sigma_{\ess} (\mathcal{L}_a)$ (with $\mathcal{L}_a$
as specified in Section \ref{sls-a-section}), then 
\begin{equation*}
    \mathcal{N} ((-\infty, \lambda_2)) 
    = \sum_{x \in \mathbb{R}} \dim \ker X(x; \lambda_2). 
\end{equation*}
\end{theorem}

\begin{theorem} \label{fourth-theorem}
For (\ref{fourth-ode}), let Assumptions {\bf (FP1)} and {\bf (FP2)} 
from Section \ref{fourth-section} hold, and express (\ref{fourth-ode})
in the form (\ref{hammy}) (giving (\ref{fourth-hammy})). Then
for $\kappa$ specified as in (\ref{fourth-kappa-defined}), {\bf (A)},
{\bf (B1)}, {\bf (B2)}, and {\bf (B3)} all hold for 
(\ref{fourth-hammy}) with $I = (- \infty, \kappa)$, and so 
the result of Theorem \ref{main-theorem} holds for 
all intervals $[\lambda_1, \lambda_2]$, $\lambda_1 < \lambda_2 < \kappa$. 
In addition, if $\mathcal{N} ([\lambda_1, \lambda_2))$
denotes the number of eigenvalues, counted with multiplicity,
that (\ref{fourth-ode}) has on the interval $[\lambda_1, \lambda_2)$, 
then
\begin{equation*}
     \mathcal{N} ([\lambda_1, \lambda_2))
     = \sum_{x \in \mathbb{R}} \dim \ker \Phi (x; \lambda_2)
        - \sum_{x \in \mathbb{R}} \dim \ker \Phi (x; \lambda_1),
\end{equation*}
where (for $i = 1, 2$)
\begin{equation*}
    \Phi (x; \lambda_i)
    = \begin{pmatrix}
    \phi_1 (x; \lambda_i) & \phi_2 (x; \lambda_i) & \dots & \phi_{2n} (x; \lambda_i) \\
    \phi_1' (x; \lambda_i) & \phi_2' (x; \lambda_i) & \dots & \phi_{2n}' (x; \lambda_i)
    \end{pmatrix},
\end{equation*}
with $\{\phi_j (x; \lambda_i)\}_{j=1}^{2n}$, $i = 1, 2$, comprising a collection of $2n$ 
linearly independent solutions of (\ref{fourth-ode}) that lie left in $\mathbb{R}$.  
Finally, 
\begin{equation*}
    \mathcal{N} ((-\infty, \lambda_2)) 
    = \sum_{x \in \mathbb{R}} \dim \ker \Phi (x; \lambda_2). 
\end{equation*}
\end{theorem}

In the remainder of this introduction, we provide some background
and context for our analysis and also set out a plan for the 
paper. For the former, our results serve as natural generalizations
of Sturm's Oscillation Theorem and the Morse Index Theorem, 
and so go respectively back to \cite{Sturm1836} and \cite{Morse1934}.
The earliest result readily identifiable with our methods 
is due to Raoul Bott in \cite{Bott1956}, followed 
by the work of Victor Maslov in \cite{Maslov1965a} and 
V. I. Arnol'd in \cite{Arnold1967}. Specific applications 
to the stability of nonlinear waves were carried about by 
Chris Jones in \cite{J1988a, J1988b}, by Jones and 
collaborators in \cite{BCJLMS2018, BJ1995, CJ2018, JM2012},
and subsequently by numerous others, including 
\cite{CB2015, CH2007, CH2014, Co2019}. 
The Maslov index is amenable to numerical computations, 
and several analyses have emphasized this aspect of the 
theory, including \cite{BM2013, CDB2009, CDB2011, Chardard2009}. 
These results have all addressed applications to 
equations of the form (\ref{allen-cahn-equation})
(though \cite{CH2007, CH2014, Co2019} address a skew-gradient
reaction term) and to  
nonlinear waves associated with homoclinic orbits 
(i.e., with $u_- = u_+$). In addition, the target Lagrangian 
subspace in the relevant calculations has typically been 
taken to be the Dirichlet Lagrangian subspace rather than the 
``natural" targets $\tilde{\ell}_+ (\lambda_1)$ 
and $\tilde{\ell}_+ (\lambda_2)$ (an exception is \cite{Co2019}).

In \cite{DJ2011}, Jones and Jian Deng applied the Maslov 
index in the setting of multidimensional Schr\"odinger 
equations,
instigating a resurgence of interest in the methods
(see also \cite{CJLS2016, LS2018}). Motivated by this work, 
the author, along with Yuri Latushkin and Alim Sukhtayev 
revisited implementations of 
the Maslov index in a single space dimension, adapting 
the spectral-flow formulation of \cite{Ph1996} to obtain
a specification of the Maslov index especially suitable 
to general linear Hamiltonian systems associated with 
either homoclinic or heteroclinic orbits \cite{HLS2017, HS2016}.
This approach was employed in \cite{HLS2018} to establish 
a result for heteroclinic traveling-wave solutions
arising in equations of the form (\ref{allen-cahn-equation}),
and was employed in \cite{HJK2018} in an analysis of 
general linear Hamiltonian systems of the form (\ref{hammy})
on finite domains.

The primary goal of the current analysis is to adapt the approach 
taken in \cite{HLS2018} (addressing equations of the form 
(\ref{allen-cahn-equation})) to the more general setting of
(\ref{hammy}). Auxiliary to this, we hope to clarify the 
mechanism by which the target Lagrangian subspaces
$\tilde{\ell}_+ (\lambda_1)$ and $\tilde{\ell}_+ (\lambda_2)$ can be
replaced by target Lagrangian subspaces for which 
all conjugate points for the Maslov index calculations have
the same direction (i.e., target Lagrangian subspaces for which 
the flow is monotonic). 

The paper is organized as follows. In Section \ref{maslov-section},
we review elements of the Maslov index that will be used 
in our development, and in Section \ref{main-section} we prove
Theorem \ref{main-theorem}. In the subsequent three sections, 
we apply Theorem \ref{main-theorem} to prove 
Theorems \ref{sls-theorem}, \ref{fourth-theorem}, and 
\ref{sls-a-theorem}.

\section{The Maslov Index} \label{maslov-section}

Our framework for computing the Maslov index is adapted from 
Section 2 of \cite{HS2019}, which is based on the spectral 
flow formulation of \cite{Ph1996}. Rather than repeating that
development here, we will only highlight the points most 
salient to the current analysis. For a full discussion of the Maslov 
index in the current setting, we refer the reader to 
\cite{HS2019}, and for a broader view of the Maslov 
index we refer the reader to \cite{BF1998, CLM1994, RS1993}.

Given any pair of Lagrangian subspaces $\ell_1$ and 
$\ell_2$ with respective frames $\mathbf{X}_1 = {X_1 \choose Y_1}$
and $\mathbf{X}_2 = {X_2 \choose Y_2}$, we consider the matrix
\begin{equation} \label{tildeW}
\tilde{W} := - (X_1 + iY_1)(X_1-iY_1)^{-1} (X_2 - iY_2)(X_2+iY_2)^{-1}. 
\end{equation}
In \cite{HS2019}, the authors establish: (1) the inverses 
appearing in (\ref{tildeW}) exist; (2) $\tilde{W}$ is independent
of the specific frames $\mathbf{X}_1$ and $\mathbf{X}_2$ (as long
as these are indeed frames for $\ell_1$ and $\ell_2$); (3) $\tilde{W}$ is 
unitary; and (4) the identity 
\begin{equation} \label{key}
\dim (\ell_1 \cap \ell_2) = \dim (\ker (\tilde{W} + I)).
\end{equation}
Given two continuous paths of Lagrangian subspaces 
$\ell_i: [0, 1] \to \Lambda (n)$, $i = 1, 2$, with 
respective frames $\mathbf{X}_i: [0,1] \to \mathbb{C}^{2n \times n}$,
relation (\ref{key}) allows us to compute the Maslov 
index $\mas (\ell_1, \ell_2; [0,1])$ as a spectral flow
through $-1$ for the path of matrices 
\begin{equation} 
\tilde{W} (t) := - (X_1 (t) + iY_1 (t))(X_1 (t)-iY_1 (t))^{-1} 
(X_2 (t) - iY_2 (t))(X_2 (t)+iY_2 (t))^{-1}. 
\end{equation}

If $-1 \in \sigma (\tilde{W} (t_*))$ for some $t_* \in [0, 1]$, 
then we refer to $t_*$ as a conjugate point, and we 
see from (\ref{key}) that the multiplicity 
of $-1$ as an eigenvalue of $\tilde{W} (t_*)$
corresponds with $\dim (\ell_1 (t_*) \cap \ell_2 (t_*))$.
We compute the Maslov index $\mas (\ell_1, \ell_2; [0, 1])$ 
by allowing $t$ to increase from $0$ to $1$ and incrementing 
the index whenever an eigenvalue crosses $-1$ in the 
counterclockwise direction, while decrementing the index
whenever an eigenvalue crosses $-1$ in the clockwise
direction. These increments/decrements are counted with 
multiplicity, so for example, if a pair of eigenvalues 
crosses $-1$ together in the counterclockwise direction, 
then a net amount of $+2$ is added to the index. Regarding
behavior at the endpoints, if an eigenvalue of $\tilde{W}$
rotates away from $-1$ in the clockwise direction as $t$ increases
from $0$, then the Maslov index decrements (according to 
multiplicity), while if an eigenvalue of $\tilde{W}$
rotates away from $-1$ in the counterclockwise direction as $t$ increases
from $0$, then the Maslov index does not change. Likewise, 
if an eigenvalue of $\tilde{W}$ rotates into $-1$ in the 
counterclockwise direction as $t$ increases
to $1$, then the Maslov index increments (according to 
multiplicity), while if an eigenvalue of $\tilde{W}$
rotates into $-1$ in the clockwise direction as $t$ increases
to $1$, then the Maslov index does not change. Finally, 
it's possible that an eigenvalue of $\tilde{W}$ will arrive 
at $-1$ for $t = t_*$ and stay. In these cases, the 
Maslov index only increments/decrements upon arrival or 
departure, and the increments/decrements are determined 
as for the endpoints (departures determined as with $t=0$,
arrivals determined as with $t = 1$).

One of the most important features of the Maslov index is homotopy invariance, 
for which we need to consider continuously varying families of Lagrangian 
paths. To set some notation, we let $\mathcal{I}$ be a closed 
interval in $\mathbb{R}$, and 
we denote by $\mathcal{P} (\mathcal{I})$ the collection 
of all paths $\mathcal{L} (t) = (\ell_1 (t), \ell_2 (t))$, where 
$\ell_1, \ell_2: \mathcal{I} \to \Lambda (n)$ are continuous paths in the 
Lagrangian--Grassmannian. We say that two paths 
$\mathcal{L}, \mathcal{M} \in \mathcal{P} (\mathcal{I})$ are homotopic provided 
there exists a family $\mathcal{H}_s$ so that 
$\mathcal{H}_0 = \mathcal{L}$, $\mathcal{H}_1 = \mathcal{M}$, 
and $\mathcal{H}_s (t)$ is continuous as a map from $(t,s) \in \mathcal{I} \times [0,1]$
into $\Lambda (n) \times \Lambda (n)$. 

The Maslov index has the following properties. 

\medskip
\noindent
{\bf (P1)} (Path Additivity) If $\mathcal{L} \in \mathcal{P} (\mathcal{I})$
and $a, b, c \in \mathcal{I}$, with $a < b < c$, then 
\begin{equation*}
\mas (\mathcal{L};[a, c]) = \mas (\mathcal{L};[a, b]) + \mas (\mathcal{L}; [b, c]).
\end{equation*}

\medskip
\noindent
{\bf (P2)} (Homotopy Invariance) If $\mathcal{L}, \mathcal{M} \in \mathcal{P} (\mathcal{I})$ 
are homotopic, with $\mathcal{L} (a) = \mathcal{M} (a)$ and  
$\mathcal{L} (b) = \mathcal{M} (b)$ (i.e., if $\mathcal{L}, \mathcal{M}$
are homotopic with fixed endpoints) then 
\begin{equation*}
\mas (\mathcal{L};[a, b]) = \mas (\mathcal{M};[a, b]).
\end{equation*} 

Straightforward proofs of these properties appear in \cite{HLS2017}
for Lagrangian subspaces of $\mathbb{R}^{2n}$, and proofs in the current setting of 
Lagrangian subspaces of $\mathbb{C}^{2n}$ are essentially identical. 

As noted previously, the direction we associate with a 
conjugate point is determined by the direction in which eigenvalues
of $\tilde{W}$ rotate through $-1$ (counterclockwise is positive, 
while clockwise is negative). In order to understand the nature
of this rotation in specific cases, we will use the following 
lemma from \cite{HS2019}. 

\begin{lemma} \label{monotonicity1}
Suppose $\ell_1, \ell_2: \mathcal{I} \to \Lambda (n)$ denote paths of 
Lagrangian subspaces of $\mathbb{C}^{2n}$ with absolutely 
continuous frames $\mathbf{X}_1 = {X_1 \choose Y_1}$
and $\mathbf{X}_2 = {X_2 \choose Y_2}$ (respectively),
and $t_0$ is any value in the interior of $\mathcal{I}$. If there exists 
$\delta > 0$ so that the matrices 
\begin{equation*}
- \mathbf{X}_1 (t)^* J \mathbf{X}_1' (t) = X_1 (t)^* Y_1' (t) - Y_1 (t)^* X_1'(t)
\end{equation*}
and (noting the sign change)
\begin{equation*}
\mathbf{X}_2 (t)^* J \mathbf{X}_2' (t) = - (X_2 (t)^* Y_2' (t) - Y_2 (t)^* X_2'(t))
\end{equation*}
are both a.e.-non-negative in $(t_0-\delta,t_0+\delta)$, and at least one is 
a.e.-positive definite in $(t_0-\delta,t_0+\delta)$ then the eigenvalues of 
$\tilde{W} (t)$ rotate in the counterclockwise direction as $t$ increases through $t_0$. 
Likewise, if both of these matrices are a.e.-non-positive, and at least one is 
a.e.-negative definite, 
then the eigenvalues of $\tilde{W} (t)$ rotate in the clockwise direction as 
$t$ increases through $t_0$.
\end{lemma}

\begin{remark}
In Theorem \ref{main-theorem}, the Maslov indices are computed on 
the unbounded interval $(-\infty, + \infty)$, and the notation 
$(-\infty, +\infty]$ is used to signify that the limit
$+ \infty$ can serve as a conjugate point. Precisely, 
under our limit assumptions in {\bf (B2)}, we can compactify 
$(-\infty, \infty)$ with a map such as 
\begin{equation*}
    x = \ln (\frac{1+\tau}{1-\tau});
    \quad \tau \in [-1, 1],
\end{equation*}
and subsequently compute the relevant Maslov indices 
on the bounded interval $\mathcal{I}=[-1, 1]$, employing the considerations 
discussed in this section. We recall from Remark 
\ref{brackets} that due to our Assumption {\bf (B2)},
$-\infty$ cannot serve as a conjugate point, and so is omitted
from the square-bracket notation.
\end{remark}

\section{Proof of Theorem \ref{main-theorem}} \label{main-section}

Before proving Theorem \ref{main-theorem}, we verify the assertion 
made in the statement of Assumption {\bf (B1)} that $\mathbf{X} (x; \lambda)$
and $\tilde{\mathbf{X}} (x; \lambda)$ are necessarily frames for Lagrangian 
subspaces of $\mathbb{C}^{2n}$ for all $(x,\lambda) \in \mathbb{R} \times I$. 
We will carry out the demonstration for $\mathbf{X} (x; \lambda)$; the case of 
$\tilde{\mathbf{X}} (x; \lambda)$ is essentially identical. First, 
we note that under our Assumption {\bf (A)} we have 
$\mathbf{X} (\cdot; \lambda), \tilde{\mathbf{X}} (\cdot; \lambda) \in \AC_{\loc} (\mathbb{R}, \mathbb{C}^{2n \times n})$
for all $\lambda \in I$ (see, e.g., Theorem 2.1 in \cite{We1987}). 
Next, according to Proposition 2.1 of \cite{HS2019}, it's sufficient to show that 
\begin{equation} \label{lagrangian-demo}
   \mathbf{X} (x; \lambda)^* J \mathbf{X} (x; \lambda) = 0,
   \quad \forall \, (x, \lambda) \in \mathbb{R} \times I.
\end{equation}
In order to verify this, we fix any $\lambda \in I$ and 
compute 
\begin{equation*}
\begin{aligned}
    \frac{\partial}{\partial x} &(\mathbf{X} (x; \lambda)^* J \mathbf{X} (x; \lambda))
    = \mathbf{X}' (x; \lambda)^* J \mathbf{X} (x; \lambda) 
    + \mathbf{X} (x; \lambda)^* J \mathbf{X}' (x; \lambda) \\
    &= - (J \mathbf{X}' (x; \lambda))^* \mathbf{X} (x; \lambda) 
    + \mathbf{X} (x; \lambda)^* J \mathbf{X}' (x; \lambda) \\
    &= - (\mathbb{B} (x; \lambda) \mathbf{X} (x; \lambda))^* \mathbf{X} (x; \lambda)
    +  \mathbf{X} (x; \lambda)^* \mathbb{B} (x; \lambda) \mathbf{X} (x; \lambda) \\
    &= 0, \quad {\rm a.e.} \,\, x \in \mathbb{R},
\end{aligned}
\end{equation*}
where in obtaining the final equality we've observed from Assumption {\bf (A)}
that $\mathbb{B} (x; \lambda)$ is self-adjoint for a.e. $x \in \mathbb{R}$. 
Recalling that $\mathbf{X} (x; \lambda)^* J \mathbf{X} (x; \lambda)$ is locally 
absolutely continuous in $\mathbb{R}$, we see that it is constant on $\mathbb{R}$.
But the columns of $\mathbf{X} (x; \lambda)$ lie left in $\mathbb{R}$, so 
\begin{equation*}
    \lim_{x \to - \infty} \mathbf{X} (x; \lambda)^* J \mathbf{X} (x; \lambda) = 0.
\end{equation*}
This calculation holds for all $\lambda \in I$,
allowing us to conclude (\ref{lagrangian-demo}). 

Turning now to the proof of Theorem \ref{main-theorem}, we 
begin by fixing any pair $\lambda_1, \lambda_2 \in \mathbb{R}$, 
$\lambda_1 < \lambda_2$, so that $[\lambda_1, \lambda_2] \subset I$,
and for all $(x, \lambda) \in \mathbb{R} \times [\lambda_1, \lambda_2]$, we let $\ell (x, \lambda)$
and $\tilde{\ell} (x; \lambda)$ denote the Lagrangian subspaces 
described in {\bf (B1)} and {\bf (B2)}. We will fix some $c > 0$
to be chosen sufficiently large during the analysis, and 
we will establish Theorem \ref{main-theorem} by 
considering the Maslov index for $\ell (x; \lambda)$ and 
$\tilde{\ell} (c; \lambda)$ along a path designated as the 
{\it Maslov box} in the next paragraph. As described in 
Section \ref{maslov-section}, this Maslov index is computed as
a spectral flow for the matrix 
\begin{equation} \label{tildeW-c}
\begin{aligned}
\tilde{W}_c (x; \lambda) &:= - (X (x; \lambda) + i Y (x; \lambda))
(X (x; \lambda) - i Y (x; \lambda))^{-1} \\
& \times (\tilde{X} (c; \lambda) - i \tilde{Y} (c; \lambda))
(\tilde{X} (c; \lambda) + i \tilde{Y} (c; \lambda))^{-1}.
\end{aligned}
\end{equation}

By Maslov Box, in this case we mean the following sequence of contours: 
(1) fix $x = -c$ and let $\lambda$ increase from $\lambda_1$ to $\lambda_2$ 
(the {\it bottom shelf}); 
(2) fix $\lambda = \lambda_2$ and let $x$ increase from $-c$ to $c$ 
(the {\it right shelf}); (3) fix $x = c$ and let $\lambda$
decrease from $\lambda_2$ to $\lambda_1$ (the {\it top shelf}); and (4) fix
$\lambda = \lambda_1$ and let $x$ decrease from $c$ to $-c$ (the 
{\it left shelf}). 

{\it The Bottom Shelf}. For the bottom shelf, the Maslov index detects
intersections between $\ell (-c; \lambda)$ and $\tilde{\ell} (c, \lambda)$
as $\lambda$ increases from $\lambda_1$ to $\lambda_2$. Since $[\lambda_1, \lambda_2]$
is compact, it follows from our Assumption {\bf (B2)} that we can take $c$ 
sufficiently large so that 
\begin{equation*}
    \ell(-c; \lambda) \cap \tilde{\ell}(c; \lambda) = \{0\},
    \quad \forall \, \lambda \in [\lambda_1, \lambda_2].
\end{equation*}
In this way, we see that 
\begin{equation*}
    \mas (\ell(-c; \cdot), \tilde{\ell} (c; \cdot); [\lambda_1, \lambda_2])
    = 0.
\end{equation*}

{\it The Top Shelf}. For the top shelf, the Maslov index detects
intersections between $\ell (c; \lambda)$ and $\tilde{\ell} (c; \lambda)$
as $\lambda$ decreases from $\lambda_2$ to $\lambda_1$. These Lagrangian 
subspaces will intersect if and only if $\lambda$ is an eigenvalue 
of (\ref{hammy}), and the multiplicity of this intersection will 
correspond with the geometric multiplicity of $\lambda$ as an 
eigenvalue of (\ref{hammy}). We would like to conclude that the 
Maslov index for the top shelf is precisely a count, including 
geometric multiplicity, of the number of eigenvalues that (\ref{hammy})
has on the interval $[\lambda_1, \lambda_2)$, and in order to draw 
this conclusion we need to know that conjugate points in this 
case all have the same (positive) direction. For this, we observe
from Lemma \ref{monotonicity1} that the direction of rotation 
associated with the Maslov index along the top shelf will be 
determined by the signs of the matrices 
\begin{equation} \label{matrix1}
    - \mathbf{X} (c; \lambda)^* J \partial_{\lambda} \mathbf{X} (c; \lambda)
\end{equation}
and 
\begin{equation} \label{matrix2}
    \tilde{\mathbf{X}} (c; \lambda)^* J \partial_{\lambda} \tilde{\mathbf{X}} (c; \lambda)
\end{equation}
in the following sense: if both of these matrices are non-positive 
at some $\lambda \in (\lambda_1, \lambda_2)$, and at least one of them 
is negative definite at $\lambda$, then the rotation at that value $\lambda$ for all 
eigenvalues of $\tilde{W} (c; \lambda)$ will be in the clockwise 
direction (with $\lambda$ increasing). 

For the first of these matrices, we compute 
\begin{equation*}
    \begin{aligned}
    \frac{\partial}{\partial x} &\mathbf{X} (x; \lambda)^* J \partial_{\lambda} \mathbf{X} (x; \lambda)
    = \mathbf{X}' (x; \lambda)^* J \partial_{\lambda} \mathbf{X} (x; \lambda)
    + \mathbf{X} (x; \lambda)^* J \partial_{\lambda} \mathbf{X}' (x; \lambda) \\
    &= - (J \mathbf{X}' (x; \lambda))^* \partial_{\lambda} \mathbf{X} (x; \lambda)
    + \mathbf{X} (x; \lambda)^* \partial_{\lambda} (J \mathbf{X}' (x; \lambda)) \\
    &= - (\mathbb{B} (x; \lambda) \mathbf{X} (x; \lambda))^* \partial_{\lambda} \mathbf{X} (x; \lambda)
    + \mathbf{X} (x; \lambda)^* \partial_{\lambda} (\mathbb{B} (x; \lambda) \mathbf{X} (x; \lambda)) \\
    &=  - \mathbf{X} (x; \lambda)^* \mathbb{B} (x; \lambda) \partial_{\lambda} \mathbf{X} (x; \lambda)
    + \mathbf{X} (x; \lambda)^* \mathbb{B}_{\lambda} (x; \lambda) \mathbf{X} (x; \lambda) \\
    &\quad \quad \quad 
    + \mathbf{X} (x; \lambda)^* \mathbb{B} (x; \lambda) \partial_{\lambda} \mathbf{X} (x; \lambda) \\
    &= \mathbf{X} (x; \lambda)^* \mathbb{B}_{\lambda} (x; \lambda) \mathbf{X} (x; \lambda).
    \end{aligned}
\end{equation*}
Upon integrating this relation on $(-\infty, c)$ and observing that 
\begin{equation*}
    \lim_{x \to - \infty} \mathbf{X} (x; \lambda)^* J \partial_{\lambda} \mathbf{X} (x; \lambda)
    = 0,
\end{equation*}
we obtain the relation 
\begin{equation*}
    \mathbf{X} (c; \lambda)^* J \partial_{\lambda} \mathbf{X} (c; \lambda)
    = \int_{-\infty}^c \mathbf{X} (x; \lambda)^* \mathbb{B}_{\lambda} (x; \lambda) \mathbf{X} (x; \lambda) d\xi.
\end{equation*}
According to Assumption {\bf (B3)}, this matrix is positive definite for $c$ sufficiently 
large, and we can conclude that (\ref{matrix1}) is negative 
definite. By a similar calculation, we can check that 
(\ref{matrix2}) is negative definite as well.  We can conclude 
from Lemma \ref{monotonicity1} that the eigenvalues of $\tilde{W} (c; \lambda)$
rotate monotonically in the clockwise direction as $\lambda$ increases
from $\lambda_1$ to $\lambda_2$, and consequently that 
\begin{equation*}
    \mathcal{N} ([\lambda_1, \lambda_2))
    = - \mas (\ell (c; \cdot), \tilde{\ell} (c; \cdot); [\lambda_1, \lambda_2]).
\end{equation*}
The inclusion of $\lambda_1$ on the left-hand side is due to the clockwise rotation as $\lambda$ 
increases, leading to a decrement of the Maslov index if $(c, \lambda_1)$ is conjugate, 
and the exclusion of $\lambda_2$ follows similarly. 

{\it The left and right shelves}. The left and right shelves are both left as computations
in Theorem \ref{main-theorem}, but in order to eliminate the arbitrary value $c$,
we need to show that by taking $c$ sufficiently large we can ensure that 
\begin{equation} \label{right-left-relation}
    \mas (\ell (\cdot; \lambda_1), \tilde{\ell} (c; \lambda_1); [-c, c])
    = \mas (\ell (\cdot; \lambda_1), \tilde{\ell}_+ (\lambda_1); (-\infty, +\infty]),
\end{equation}
and similarly for $\lambda_2$. In order to understand why 
(\ref{right-left-relation}) holds, it's convenient to observe
that the left-hand hand side is computed via the matrix
$\tilde{W}_c (x; \lambda_1)$ (i.e., (\ref{tildeW}) with 
$\lambda = \lambda_1$), while the right-hand side is computed
via the matrix
\begin{equation*}
\begin{aligned}
    \tilde{\mathcal{W}} (x; \lambda_1) 
    &= - (X (x; \lambda_1) + i Y (x; \lambda_1)) (X (x; \lambda_1) - i Y (x; \lambda_1))^{-1} \\
    &\times (\tilde{X}_+ (\lambda_1) - i \tilde{Y}_+ (\lambda_1)) (\tilde{X}_+ (\lambda_1) + i \tilde{Y}_+ (\lambda_1))^{-1}.
\end{aligned}
\end{equation*}
Comparing expressions for $\tilde{W}_c (x; \lambda_1)$ and $\tilde{\mathcal{W}} (x; \lambda_1)$,
we see that we can write $\tilde{W}_c (x; \lambda_1) = \tilde{\mathcal{W}} (x; \lambda_1) \tilde{V} (c; \lambda_1)$,
where 
\begin{equation*}
    \begin{aligned}
     \tilde{V} (c; \lambda_1) 
     &= (\tilde{X}_+ (\lambda_1) + i \tilde{Y}_+ (\lambda_1)) (\tilde{X}_+ (\lambda_1) - i \tilde{Y}_+ (\lambda_1))^{-1} \\   
     &\times (\tilde{X} (c; \lambda_1) - i \tilde{Y} (c; \lambda_1)) (\tilde{X} (c; \lambda_1) + i \tilde{Y} (c; \lambda_1))^{-1}.
    \end{aligned}
\end{equation*}
Here, $\tilde{V} (c; \lambda_2)$ is a continuous function of $c$, satisfying
\begin{equation*}
    \lim_{c \to + \infty} \tilde{V} (c; \lambda_2) = I.
\end{equation*}

Let $\{w_j^c (x; \lambda_1)\}_{j=1}^n$ denote the eigenvalues of 
$\tilde{W}_c (x; \lambda_1)$, and let $\{\omega_j (x; \lambda_1)\}_{j=1}^n$
denote the eigenvalues of $\tilde{\mathcal{W}} (x; \lambda_1)$. Using 
Assumption {\bf (B2)}, we see that the limits 
\begin{equation*}
    \begin{aligned}
    \tilde{W}_c^- (\lambda_1) &:= \lim_{x \to - \infty} \tilde{W}_c (x; \lambda_1);
    \quad \tilde{\mathcal{W}}^- (\lambda_1) := \lim_{x \to - \infty} \tilde{\mathcal{W}} (x; \lambda_1) \\
    \tilde{W}_c^+ (\lambda_1) &:= \lim_{x \to + \infty} \tilde{W}_c (x; \lambda_1);
    \quad \tilde{\mathcal{W}}^+ (\lambda_1) := \lim_{x \to + \infty} \tilde{\mathcal{W}} (x; \lambda_1)
    \end{aligned}
\end{equation*}
are all well-defined. It follows that given any $\epsilon > 0$, we can find 
$L, c_0 > 0$ sufficiently large so that  
for each $j \in \{1,2,\dots, n\}$ (with an appropriate choice of labeling)
\begin{equation*}
    |w_j^c (x; \lambda_1) - \omega_j (x; \lambda_1)| < \epsilon,
    \quad \forall \, |x| > L, \, c > c_0. 
\end{equation*}
Likewise, given any $\epsilon > 0$ we can use compactness of 
$[-L, L]$ to take $c_1 > c_0$ sufficiently large so that  
\begin{equation*}
    |w_j^c (x; \lambda_1) - \omega_j (x; \lambda_1)| < \epsilon,
    \quad \forall \, x \in [-L,L], \, c > c_1. 
\end{equation*}
Combining these observations, we see that given any $\epsilon > 0$
we can take $c_1$ sufficiently large so that 
\begin{equation} \label{eig-diff}
    |w_j^c (x; \lambda_1) - \omega_j (x; \lambda_1)| < \epsilon,
    \quad \forall \, x \in \mathbb{R}, \, c > c_1.
\end{equation}
We also note that according to the second part of Assumption {\bf (B2)}, we can take $c$ large 
enough so that we have both $-1 \notin \sigma (\tilde{W}_c (-c; \lambda_1))$ 
and $-1 \notin \sigma (\tilde{\mathcal{W}} (-c; \lambda_1))$. 

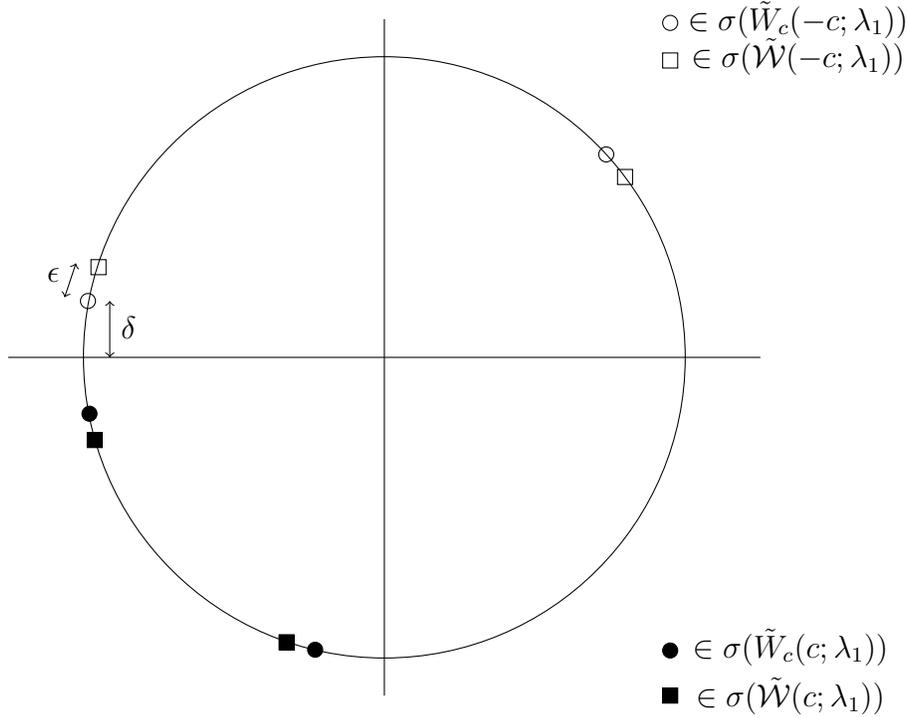
\begin{figure}[ht]
\begin{center}
\begin{tikzpicture}
\draw (-5,0) -- (5,0);	%horizontal line
\draw (0,-4.5) -- (0,4.5);	%vertical line
\draw (0,0) circle (4);	%circle
\draw (3.8,4.45) circle (.1);
\node at (5.5,4.5) {$\in \sigma (\tilde{W}_c (-c; \lambda_1))$};
\draw (3.70,3.85) rectangle (3.90,4.05);
\node at (5.5,4.0) {$\in \sigma (\tilde{\mathcal{W}} (-c; \lambda_1))$};
\filldraw [black] (3.8,-3.90) circle (.1);
\node at (5.4,-3.85) {$\in \sigma (\tilde{W}_c (c; \lambda_1))$};
\filldraw [black] (3.70,-4.6) rectangle (3.90,-4.4);
\node at (5.4,-4.5) {$\in \sigma (\tilde{\mathcal{W}} (c; \lambda_1))$};
%
%PLACE SHAPES ON THE CIRCLE
%PAIR 1
\draw (-3.94,.75) circle (.1);
\draw (-3.9,1.1) rectangle (-3.7,1.3);
%
%PAIR 2
\draw (2.95,2.7) circle (.1);
\draw (3.1,2.3) rectangle (3.3,2.5);
%
%PAIR 3
\filldraw [black] (-3.92,-.75) circle (.1);
\filldraw [black] (-3.95,-1.2) rectangle (-3.75,-1.0);
%
%PAIR 4
\filldraw [black] (-.92,-3.89) circle (.1);
\filldraw [black] (-1.4,-3.89) rectangle (-1.2,-3.69);
%
%ADD EPSILON AND DELTA
\draw[<->] (-4.25,.8) -- (-4.1,1.25);
\node at (-4.4,1.1) {$\epsilon$};
\draw[<->] (-3.65,0) -- (-3.65,.75);
\node at (-3.4,.4) {$\delta$};
\end{tikzpicture}
\end{center}
\caption{Eigenvalues of $\tilde{W}_c (\pm c; \lambda_1)$ and 
$\tilde{\mathcal{W}} (\pm c; \lambda_1)$. \label{A}}
\end{figure}

At this point, we divide the analysis into two cases: 
(1) $\lambda_1$ is not an eigenvalue of (\ref{hammy});
and (2) $\lambda_1$ is an eigenvalue of (\ref{hammy}). 
For Case (1), suppose $\lambda_1$ is not an eigenvale of (\ref{hammy}).
Then we immediately have $-1 \notin \sigma (\tilde{W}_c (c; \lambda_1))$
(for any $c \in \mathbb{R}$), and since 
$\ell_+ (\lambda_1) \cap \tilde{\ell}_+ (\lambda_1) = \{0\}$,
we can take $c$ large enough so that 
$-1 \notin \sigma (\tilde{\mathcal{W}} (c; \lambda_1))$. 
In summary, the situation is as follows: for 
$\tilde{W}_c (x; \lambda_1)$ we have both 
$-1 \notin \sigma (\tilde{W}_c (-c; \lambda_1))$ 
and $-1 \notin \sigma (\tilde{W}_c (c; \lambda_1))$,
and likewise for $\tilde{\mathcal{W}} (x; \lambda_1)$
we have both $-1 \notin \sigma (\tilde{\mathcal{W}} (-c; \lambda_1))$
and $-1 \notin \sigma (\tilde{\mathcal{W}} (c; \lambda_1))$.
It follows that there exists some $\delta > 0$ so that 
\begin{equation*}
    |w_j^c (-c; \lambda_1) + 1| > \delta,
    \quad |w_j^c (c; \lambda_1) + 1| > \delta,
    \quad \forall \, j \in \{1, 2, \dots, n\},
\end{equation*}
and 
\begin{equation*}
    |\omega_j (-c; \lambda_1) + 1| > \delta,
    \quad |\omega_j (c; \lambda_1) + 1| > \delta,
    \quad \forall \, j \in \{1, 2, \dots, n\}.
\end{equation*}
See Figure \ref{A}, sketched for the case $n=2$.

Using (\ref{eig-diff}), we can take $c$ large enough 
so that $\epsilon < \delta$. In this way, 
as $x$ increases from $-c$ to $c$, an eigenvalue
of $\tilde{W}_c (x; \lambda)$ can complete a full loop 
around $S^1$ if and only if a corresponding 
eigenvalue of $\tilde{\mathcal{W}} (x; \lambda_1)$
also completes a full loop. In addition, since the 
distance between eigenvalues of $\tilde{W}_c (x; \lambda)$
and eigenvalues of $\tilde{\mathcal{W}} (x; \lambda_1)$
is less than the initial and final distances of 
the eigenvalues of these matrices from $-1$, 
the total count of conjugate points associated with 
the Maslov index computed via $\tilde{W}_c (x; \lambda_1)$ 
must be precisely the corresponding count
computed via $\tilde{\mathcal{W}} (x; \lambda_1)$.
We can conclude that 
\begin{equation} \label{truncated-equality}
    \mas (\ell (\cdot; \lambda_1), \tilde{\ell} (c; \lambda_1); [-c, c])
    = \mas (\ell (\cdot; \lambda_1), \tilde{\ell}_+ (\lambda_1); [-c, c]),
\end{equation}
for all $c$ sufficiently large. According to {\bf (B2)}, 
we can take $c$ sufficiently large so that 
$\ell (x; \lambda_1) \cap \tilde{\ell}_+ (\lambda_1) = \{0\}$
for all $x < -c$, and since $\lambda_1$ is 
not an eigenvalue of (\ref{hammy}), we can take $c$ sufficiently 
large so that $\ell (x; \lambda_1) \cap \tilde{\ell}_+ (\lambda_1) = \{0\}$
for all $x > c$. We conclude that 
\begin{equation*}
    \mas (\ell (\cdot; \lambda_1), \tilde{\ell}_+ (\lambda_1); (-\infty, -c]) = 0,
\end{equation*}
and 
\begin{equation*}
    \mas (\ell (\cdot; \lambda_1), \tilde{\ell}_+ (\lambda_1); [c, +\infty]) = 0.
\end{equation*}
This allows us to conclude (\ref{right-left-relation}) in the case that 
$\lambda_1$ is not an eigenvalue of (\ref{hammy}). 

For Case (2), we assume $\lambda_1$ is an eigenvalue of (\ref{hammy}),
and in order to be definite we will specify its geometric multiplicity 
as $m$. We continue to have $-1 \notin \sigma (\tilde{W}_c (-c; \lambda_1))$ 
and $-1 \notin \sigma (\tilde{\mathcal{W}} (-c; \lambda_1))$ (for $c$ sufficiently
large), but now $-1 \in \sigma (\tilde{W}_c (c; \lambda_1))$ 
with multiplicity $m$ (and it's not definite whether $-1$ is in 
the spectrum of $\tilde{\mathcal{W}} (c; \lambda_1)$). The matrix 
$\tilde{W}_c (c; \lambda_1)$ will have $n-m$ eigenvalues located 
away from $-1$, and there will correspond $n-m$ eigenvalues
of $\tilde{\mathcal{W}} (c; \lambda_1)$ located away from $-1$
(with the remaining eigenvalues of $\tilde{\mathcal{W}} (c; \lambda_1)$
necessarily near $-1$). The flow associated with these $n - m$ eigenvalues 
can be analyzed precisely as in Case (1). We now consider the 
$m$ eigenvalues of $\tilde{\mathcal{W}} (x; \lambda_1)$ near 
$-1$ at $x = c$. This group of eigenvalues will track the group
of eigenvalues of $\tilde{W}_c (x; \lambda_1)$ that 
approach $-1$ as $x \to c^-$. The evolution will proceed as in 
the previous case, except that in this case the right and left
sides of (\ref{truncated-equality}) won't necessarily agree. 
For $\tilde{W}_c (x; \lambda_1)$, the evolution stops at $x = c$,
but for $\tilde{\mathcal{W}} (x; \lambda_1)$, it continues as 
$x$ tends to $+ \infty$. Moreover, as $x$ tends to $+ \infty$
the $m$ eigenvalues of $\tilde{\mathcal{W}} (x; \lambda_1)$ that 
are not bounded away from $-1$ will necessarily approach $-1$
in the asymptotic limit. At this point, it's critical to observe
that this set of $m$ eigenvalues of $\tilde{\mathcal{W}} (x; \lambda_1)$
cannot complete a loop of $S^1$ as $x$ increases from $c$ (because
they must remain near $-1$), and so the signs associated with 
their approaches to $-1$ have already been determined by the time
$x$ arrives at $c$. In particular, these signs must agree with those
of the eigenvalues of $\tilde{W}_c (x; \lambda_1)$ that approach $-1$
as $x \to c^-$. In this way, we conclude 
\begin{equation*}
    \mas (\ell (\cdot; \lambda_1), \tilde{\ell} (c; \lambda_1); [-c, c])
    = \mas (\ell (\cdot; \lambda_1), \tilde{\ell}_+ (\lambda_1); [-c, +\infty]),
\end{equation*}
and extension of the right-hand side to $(-\infty, +\infty]$ is precisely 
as before. This gives (\ref{right-left-relation}) in Case (2). The same 
considerations hold for $\lambda_2$. 

Combining these observations, we can use catenation of paths along with 
homotopy invariance to write (respectively)
\begin{equation*}
    \begin{aligned}
        0 &= \textrm{bottom shelf} + \textrm{right shelf}
        + \textrm{top shelf} + \textrm{left shelf} \\
        &= 0 +  \mas (\ell (\cdot; \lambda_2), \tilde{\ell} (c; \lambda_2); [-c, c])
        + \mathcal{N} ([\lambda_1, \lambda_2)) -  \mas (\ell (\cdot; \lambda_1), \tilde{\ell} (c; \lambda_1); [-c, c]) \\
        &=  \mas (\ell (\cdot; \lambda_2), \tilde{\ell}_+ (\lambda_2); (-\infty, +\infty])
        + \mathcal{N} ([\lambda_1, \lambda_2))
        - \mas (\ell (\cdot; \lambda_1), \tilde{\ell}_+ (\lambda_1); (-\infty, +\infty]).
    \end{aligned}
\end{equation*}
Rearranging terms, we obtain precisely the claim of Theorem \ref{main-theorem}.
\hfill $\square$

\section{Sturm-Liouville Systems} \label{sls-section}

In this section, we apply Theorem \ref{main-theorem} to Sturm-Liouville 
systems
\begin{equation} \label{sls-again}
    - (P(x) \phi')' + V(x) \phi = \lambda Q (x) \phi;
    \quad x \in \mathbb{R}, \, \phi (x; \lambda) \in \mathbb{C}^n,
\end{equation}
and we also establish the additional properties stated in 
Theorem \ref{sls-theorem}. In order to ensure that our general 
assumptions {\bf (A)}, {\bf (B1)}, {\bf (B2)}, and {\bf (B3)} hold, 
we make the following assumptions on the coefficient matrices 
$P$, $V$, and $Q$. 

\medskip
{\bf (SL1)} We take $P \in \AC_{\loc} (\mathbb{R}, \mathbb{C}^{n \times n})$
and $V, Q \in L^1_{\loc} (\mathbb{R}, \mathbb{C}^{n \times n})$, with 
$P(x)$, $V(x)$, and $Q(x)$ self-adjoint for a.e. $x \in \mathbb{R}$.
Moreover, we assume that there exist constants $\theta_P, \theta_Q > 0$ and a 
constant $C_V \ge 0$ so that for any $v \in \mathbb{C}^n$
\begin{equation*}
    (P(x)v, v) \ge \theta_P |v|^2; 
    \quad (Q(x)v, v) \ge \theta_Q |v|^2;
    \quad (V(x)v, v) \le C_v |v|^2,
\end{equation*}
for a.e. $x \in \mathbb{R}$. 

\medskip
{\bf (SL2)} Each of the matrices $P$, $V$, and $Q$ approaches 
well-defined asymptotic endstates at exponential rate as 
$x \to \pm \infty$. Precisely, there exist self-adjoint 
matrices $P_{\pm}, V_{\pm}, Q_{\pm} \in \mathbb{C}^{n \times n}$,
with $P_{\pm}, Q_{\pm}$ positive definite, along with constants
$C, M \ge 0$, $\eta > 0$, so that 
\begin{equation*}
    \begin{aligned}
        |P(x) - P_{\pm}| \le C e^{- \eta |x|};
        \quad |V(x) - V_{\pm}| &\le C e^{- \eta |x|};
        \quad |Q(x) - Q_{\pm}| \le C e^{- \eta |x|}
        \quad {\rm a.e.}\, x \gtrless \pm M, \\
        \quad |P'(x)| &\le C e^{- \eta |x|}
        \quad  {\rm a.e.}\, x \gtrless \pm M. 
    \end{aligned}
\end{equation*}
 
\begin{remark}
 The boundedness assumption on $V(x)$ is only required for the final 
 claim in Theorem \ref{sls-theorem}, in which $\lambda_1$ is taken 
 sufficiently negative so that 
 \begin{equation*}
     \mathcal{N} ((-\infty, \lambda_2)) = 
     \sum_{x \in \mathbb{R}} \dim \ker X (x; \lambda_2).
 \end{equation*}
\end{remark} 
 
We can associate with (\ref{sls-again}) the operator 
\begin{equation} \label{sls-L}
\mathcal{L} \phi := Q(x)^{-1} \Big{\{} - (P(x)\phi')' + V(x) \phi \Big{\}},
\end{equation}
for which we assign the domain
\begin{equation} \label{sls-domain}
    \mathcal{D} := \Big{\{} \phi \in L^2 (\mathbb{R}, \mathbb{C}^n): 
    \phi, \phi' \in \AC_{\loc} (\mathbb{R}, \mathbb{C}^n), \,
    \mathcal{L} \phi \in L^2 (\mathbb{R}, \mathbb{C}^n)
    \Big{\}},
\end{equation}
and we also introduce the inner product
\begin{equation*}
    \langle \phi, \psi \rangle_Q 
    := \int_{\mathbb{R}} (Q(x) \phi (x), \psi (x)) dx.
\end{equation*}
With this choice of domain and inner product, $\mathcal{L}$
is densely defined, closed, and self-adjoint, so 
$\sigma (\mathcal{L}) \subset \mathbb{R}$ (see, e.g., 
\cite{We1987}). 

As shown in \cite{He1981, KP2013}, the essential spectrum of 
$\mathcal{L}$ is entirely determined by the asymptotic systems
\begin{equation} \label{asls}
    - P_{\pm} \phi'' + V_{\pm} \phi = \lambda Q_{\pm} \phi
\end{equation}
in the following way: the essential spectrum is precisely 
the collection of values $\lambda \in \mathbb{R}$ for which there exists 
a solution to (\ref{asls}) of the form $\phi (x) = e^{i k x} r$
for some constant scalar $k \in \mathbb{R}$ and some constant
non-zero vector $r \in \mathbb{C}^n$. Upon substitution of $\phi (x) = e^{i k x} r$
into (\ref{asls}), we obtain the relation 
\begin{equation*}
    (k^2 P_{\pm} + V_{\pm})r = \lambda Q_{\pm} r.
\end{equation*}
If we compute an inner product of this equation with $r$,
we find 
\begin{equation*}
    k^2 (P_{\pm}r,r) + (V_{\pm}r,r) = \lambda (Q_{\pm} r, r).
\end{equation*}
Since $P_{\pm}, Q_{\pm}$ are positive definite, we see that 
\begin{equation*}
    \lambda \ge \frac{(V_{\pm} r, r)}{(Q_{\pm} r, r)},
    \quad \forall \, k \in \mathbb{R}.
\end{equation*}
We'll set 
\begin{equation} \label{kappa-defined}
    \kappa := 
    \min \Big{\{}\inf_{r \in \mathbb{C}^n \backslash \{0\}} \frac{(V_{-} r, r)}{(Q_{-} r, r)}, 
    \inf_{r \in \mathbb{C}^n \backslash \{0\}} \frac{(V_{+} r, r)}{(Q_{+} r, r)} \Big{\}}.
\end{equation}
Then 
\begin{equation*}
    \sigma_{\ess} (\mathcal{L}) = [\kappa, + \infty),
\end{equation*}
and we can conclude that for (\ref{sls-again}) we can take the interval 
$I$ described in Assumptions {\bf (A)}, {\bf (B1)}, {\bf (B2)}, 
and {\bf (B3)} to be $I = (-\infty, \kappa)$. 

Next, in order to describe the Lagrangian subspaces $\ell (x; \lambda)$
and $\tilde{\ell} (x; \lambda)$ specified in our general assumptions 
{\bf (B1)}, {\bf (B2)}, and {\bf (B3)}
we'll need a characterization of 
solutions to (\ref{sls-again}) that lie left in $\mathbb{R}$, along
with a characterization of solutions to (\ref{sls-again}) that 
lie right in $\mathbb{R}$. For this, we begin by fixing some $\lambda < \kappa$
and looking for solutions of (\ref{asls}) of the form 
$\phi (x; \lambda) = e^{\mu (\lambda) x} r (\lambda)$,
where in this case $\mu$ is a scalar function of $\lambda$
and $r$ is a vector-valued function of $\lambda$ with 
$r(\lambda) \in \mathbb{C}^n$. We find that 
\begin{equation*}
    (-\mu^2 P_{\pm} + V_{\pm} - \lambda Q_{\pm}) r = 0,
\end{equation*}
which we can rearrange as 
\begin{equation*}
    P_{\pm}^{-1} (V_{\pm} - \lambda Q_{\pm}) r =  \mu^2 r. 
\end{equation*}
Since the matrices $P_{\pm}$ are positive definite, it's natural to work 
with the inner product 
\begin{equation} \label{sls-ip}
    (r, s)_{\pm} := (P_{\pm} r, s),
\end{equation}
and it's clear that for $\lambda < \kappa$, the matrices
$P_{\pm}^{-1} (V_{\pm} - \lambda Q_{\pm})$ are self-adjoint and 
positive definite with this inner product. We conclude that
the values $\mu^2$ will be positive real values, and 
that the associated eigenvectors can be chosen to be orthonormal
with respect to (\ref{sls-ip}). For each of the $n$ values
of $\mu^2$, we can associate two values $\pm \sqrt{\mu^2}$. 
By a choice of labeling, we can split these
into $n$ negative values $\{\mu_k^{\pm}\}_{k=1}^n$ and 
$n$ positive values $\{\mu_k^{\pm}\}_{k=n+1}^{2n}$, with the
correspondence (again, by labeling convention)
\begin{equation*}
    \mu_{n+k}^{\pm} (\lambda) = - \mu_{k}^{\pm} (\lambda),
    \quad k = 1, 2, \dots, n.
\end{equation*}
For each $k \in \{1, 2, \dots, n\}$, we denote by $r_k^{\pm} (\lambda)$
the eigenvector of $P_{\pm}^{-1} (V_{\pm} - \lambda Q_{\pm})$ with 
associated eigenvalue $\mu_k^{\pm\, 2} = \mu_{n + k}^{\pm\, 2} (\lambda)$.
I.e., 
\begin{equation*}
    P_{\pm}^{-1} (V_{\pm} - \lambda Q_{\pm}) r_k^{\pm} =  \mu_k^{\pm\, 2} r_k^{\pm},
    \quad \forall \, k \in \{1, 2, \dots, n\}. 
\end{equation*}

In order to place (\ref{sls-again}) in our general framework, we set 
$y = {y_1 \choose y_2} = {\phi \choose P(x) \phi'}$, so that we have
\begin{equation} \label{sls-A-equation}
    y' = \mathbb{A} (x; \lambda) y;
    \quad \mathbb{A} (x; \lambda) = 
    \begin{pmatrix}
    0 & P(x)^{-1} \\
    V(x) - \lambda Q(x) & 0
    \end{pmatrix},
\end{equation}
or equivalently
\begin{equation} \label{sls-hammy}
    Jy' = \mathbb{B} (x; \lambda) y;
    \quad \mathbb{B} (x; \lambda) = 
    \begin{pmatrix}
    \lambda Q(x) - V(x) & 0 \\
    0 & P(x)^{-1}
    \end{pmatrix}.
\end{equation}
We see immediately that under Assumption {\bf (SL1)}, 
our general Assumptions {\bf (A)} holds. 

If we set 
\begin{equation*}
    \mathbb{A}_{\pm} (\lambda) 
    := \begin{pmatrix}
     0 & P_{\pm}^{-1} \\
    V_{\pm} - \lambda Q_{\pm} & 0
    \end{pmatrix}, 
\end{equation*}
then under our Assumptions {\bf (SL2)} we have
the relations 
\begin{equation*}
     \quad |\mathbb{A} (x; \lambda) - \mathbb{A}_{\pm} (\lambda)| \le \tilde{C} e^{- \tilde{\eta} |x|}
        \quad  {\rm a.e.}\, x \gtrless \pm \tilde{M}. 
\end{equation*}
for some constants $\tilde{C}, \tilde{M} \ge 0$, $\tilde{\eta} > 0$. 

The values $\{\mu_k^{\pm} (\lambda)\}_{k=1}^{2n}$ described above comprise a labeling
of the eigenvalues of $\mathbb{A}_{\pm} (\lambda)$. If we let $\{\mathbf{r}_k^{\pm}\}_{k=1}^{2n}$
denote the eigenvectors of $\mathbb{A}_{\pm} (\lambda)$ respectively associated with 
these eigenvectors, then we find 
\begin{equation} \label{sls-bold-r}
    \mathbf{r}_k^{\pm} (\lambda) = {r_k^{\pm} (\lambda) \choose \mu_k^{\pm} (\lambda) P_{\pm} r_k^{\pm} (\lambda)}; 
    \quad \mathbf{r}_{n+k}^{\pm} (\lambda) = {r_k^{\pm} (\lambda) \choose - \mu_k^{\pm} (\lambda) P_{\pm} r_k^{\pm} (\lambda)};
    \quad k = 1, 2, \dots, n.
\end{equation}
We'll set 
\begin{equation*}
    R_{\pm} (\lambda) = (r_1^{\pm} (\lambda) \,\,\, r_2^{\pm} (\lambda) \,\,\, \dots \,\,\, r_n^{\pm} (\lambda))
\end{equation*}
and 
\begin{equation*}
    D_{\pm} (\lambda) = \diag (\mu_1^{\pm} (\lambda) \,\,\, \mu_2^{\pm} (\lambda) \,\,\, \dots \,\,\, \mu_n^{\pm} (\lambda)),
\end{equation*}
allowing us to express a frame for the eigenspace of $\mathbb{A}_- (\lambda)$ associated 
with its positive eigenvalues as 
\begin{equation*}
    \mathbf{X}_- (\lambda) 
    = {R_- (\lambda) \choose - P_- R_- (\lambda) D_- (\lambda)}.
\end{equation*}
Likewise, we can express a frame for the eigenspace of $\mathbb{A}_- (\lambda)$ associated 
with its negative eigenvalues as 
\begin{equation*}
    \mathbf{X}_-^g (\lambda) 
    = {R_- (\lambda) \choose P_- R_- (\lambda) D_- (\lambda)},
\end{equation*}
where the superscript $g$ indicates that solutions to 
(\ref{asls}) associated with negative eigenvalues of $\mathbb{A}_-$ will grow as 
$x$ tends to $-\infty$. In the same way, we can 
express a frame for the eigenspace of $\mathbb{A}_+ (\lambda)$ associated 
with its negative eigenvalues as 
\begin{equation*}
    \tilde{\mathbf{X}}_+ (\lambda) 
    = {R_+ (\lambda) \choose P_+ R_+ (\lambda) D_+ (\lambda)},
\end{equation*}
and a frame for the eigenspace of $\mathbb{A}_+ (\lambda)$ associated 
with its positive eigenvalues as 
\begin{equation*}
    \tilde{\mathbf{X}}_+^g (\lambda) 
    = {R_+ (\lambda) \choose - P_+ R_+ (\lambda) D_+ (\lambda)}.
\end{equation*}

The following lemma can be adapted directly from 
Lemma 2.1 in \cite{HS2020}, and we refer the reader to that
reference for the proof.

\begin{lemma} \label{sls-ode-lemma} 
Assume {\bf (SL1)} and {\bf (SL2)} hold, and let $\{\mu_k^{\pm} (\lambda)\}_{k=1}^{2n}$
and $\{\mathbf{r}_k^{\pm} (\lambda)\}_{k=1}^{2n}$ be as described just above.
Then there exists a family of bases
$\{\mathbf{y}_{k}^- (\cdot; \lambda)\}_{k=n+1}^{2n}$,
$\lambda \in (-\infty, \kappa)$,
for the spaces of solutions to (\ref{sls-A-equation}) that lie
left in $\mathbb{R}$, and a family of bases
$\{\mathbf{y}_{k}^+ (\cdot; \lambda)\}_{k=1}^{n}$, 
$\lambda \in (-\infty, \kappa)$,
for the spaces of solutions to (\ref{sls-A-equation}) that lie
right in $\mathbb{R}$. Respectively, we can choose these
so that 
\begin{equation*}
\begin{aligned}
  \mathbf{y}_{n+k}^- (x; \lambda)
  &= e^{- \mu_k^- (\lambda) x} (\mathbf{r}_{n+k}^- (\lambda) + \mathbf{E}_{n+k}^- (x; \lambda)),
  \quad k = 1, 2, \dots, n, \\
   \mathbf{y}_{k}^+ (x; \lambda)
  &= e^{\mu_k^+ (\lambda) x} (\mathbf{r}_{k}^+ (\lambda) + \mathbf{E}_{k}^+ (x; \lambda)),
  \quad k = 1, 2, \dots, n,
\end{aligned}
\end{equation*}
where for any fixed interval $[\lambda_1, \lambda_2]$, with $\lambda_1 < \lambda_2 < \kappa$,
there exist a constant $\delta > 0$ so that for each 
$k \in \{1, 2, \dots, n\}$
\begin{equation*}
    \mathbf{E}_{n+k}^- (x; \lambda) = \mathbf{O} (e^{- \delta |x|}), 
    \quad x \to - \infty; 
    \quad \quad  \mathbf{E}_{k}^+ (x; \lambda) = \mathbf{O} (e^{- \delta |x|}), 
    \quad x \to + \infty,
\end{equation*}
uniformly for $\lambda \in [\lambda_1, \lambda_2]$. 

Moreover, there exists a $\lambda$-dependent family of bases
$\{\mathbf{y}_{k}^- (\cdot; \lambda)\}_{k=1}^{n}$, 
$\lambda \in (-\infty, \kappa)$,
for the spaces of solutions to (\ref{sls-A-equation}) that do not 
lie left in $\mathbb{R}$, and a $\lambda$-dependent family of bases
$\{\mathbf{y}_{k}^+ (\cdot; \lambda)\}_{k=n+1}^{2n}$, 
$\lambda \in (-\infty, \kappa)$,
for the spaces of solutions to (\ref{sls-A-equation}) that 
do not lie right in $\mathbb{R}$. Respectively, we can choose these
so that 
\begin{equation*}
\begin{aligned}
  \mathbf{y}_{k}^- (x; \lambda)
  &= e^{\mu_k^- (\lambda) x} (\mathbf{r}_{k}^- (\lambda) + \mathbf{E}_{k}^- (x; \lambda)),
  \quad k = 1, 2, \dots, n, \\
   \mathbf{y}_{n+k}^+ (x; \lambda)
  &= e^{- \mu_k^+ (\lambda) x} (\mathbf{r}_{n+k}^+ (\lambda) + \mathbf{E}_{n+k}^+ (x; \lambda)),
  \quad k = 1, 2, \dots, n,
\end{aligned}
\end{equation*}
where for any fixed interval $[\lambda_1, \lambda_2]$, with $\lambda_1 < \lambda_2 < \kappa$,
there exist a constant $\delta > 0$ so that for each 
$k \in \{1, 2, \dots, n\}$
\begin{equation*}
    \mathbf{E}_{k}^- (x; \lambda) = \mathbf{O} (e^{- \delta |x|}), 
    \quad x \to - \infty; 
    \quad \quad  \mathbf{E}_{n+k}^+ (x; \lambda) = \mathbf{O} (e^{- \delta |x|}), 
    \quad x \to + \infty,
\end{equation*}
uniformly for $\lambda \in [\lambda_1, \lambda_2]$. 
\end{lemma}

In addition to the structural assertions of Lemma \ref{sls-ode-lemma},
we need to establish the continuity and differentiability in 
$\lambda$ specified in Assumption {\bf (B1)}. For this, we 
take advantage of the observation that we can work with any 
valid frames for $\ell (x; \lambda)$ and $\tilde{\ell} (x; \lambda$.

\begin{lemma} \label{sls-extension-lemma1}
Assume {\bf (SL1)} and {\bf (SL2)} hold, and for 
each $\lambda \in (-\infty, \kappa)$ 
let $\{\mathbf{y}_{k}^- (\cdot; \lambda)\}_{k=n+1}^{2n}$
and $\{\mathbf{y}_{k}^+ (\cdot; \lambda)\}_{k=1}^{n}$
be as described in Lemma \ref{sls-ode-lemma}. If 
$\ell (x; \lambda)$ and $\tilde{\ell} (x; \lambda)$
respectively denote the Lagrangian subspaces with 
frames 
\begin{equation} \label{sls-lframe}
    \mathbf{X} (x; \lambda)
    = (\mathbf{y}_{n+1}^- (x; \lambda) \,\,\, \mathbf{y}_{n+2}^- (x; \lambda) \,\,\, \cdots \,\,\, \mathbf{y}_{2n}^- (x; \lambda)),
\end{equation}
and 
\begin{equation} \label{sls-rframe}
    \tilde{\mathbf{X}} (x; \lambda)
    = (\mathbf{y}_{1}^+ (x; \lambda) \,\,\, \mathbf{y}_{2}^+ (x; \lambda) \,\,\, \cdots \,\,\, \mathbf{y}_{n}^+ (x; \lambda)),
\end{equation}
then $\ell, \tilde{\ell} \in C (\mathbb{R} \times (-\infty, \kappa), \Lambda (n))$. 
\end{lemma}

\begin{lemma} \label{sls-extension-lemma2}
Assume {\bf (SL1)} and {\bf (SL2)} hold, and for some 
fixed $\lambda_0 \in (-\infty, \kappa))$ let
$\{\mathbf{y}_{k}^- (\cdot; \lambda_0)\}_{k=n+1}^{2n}$
and $\{\mathbf{y}_{k}^+ (\cdot; \lambda_0)\}_{k=1}^{n}$
be as described in Lemma \ref{sls-ode-lemma}. Then 
there exists a constant $r_0 > 0$ so that the elements
$\{\mathbf{y}_{k}^- (\cdot; \lambda_0)\}_{k=n+1}^{2n}$
and $\{\mathbf{y}_{k}^+ (\cdot; \lambda_0)\}_{k=1}^{n}$
can be analytically extended in $\lambda$ to the complex 
ball $B (\lambda_0, r_0)$ (centered at $\lambda_0$ with radius
$r_0$). Moreover, the analytic extensions of 
$\{\mathbf{y}_{k}^- (\cdot; \lambda_0)\}_{k=n+1}^{2n}$
comprise a basis for the space of solutions of (\ref{sls-hammy})
that lie left in $\mathbb{R}$, and the analytic extensions 
of $\{\mathbf{y}_{k}^+ (\cdot; \lambda_0)\}_{k=1}^{n}$
comprise a basis for the space of solutions of (\ref{sls-hammy})
that lie right in $\mathbb{R}$. 
\end{lemma}

\begin{remark}
The significance of Lemma \ref{sls-extension-lemma1}
lies in the assertion that in addition to being
continuous in $x$, $\ell$ and $\tilde{\ell}$ are 
continuous in $\lambda$ as well. The significance of Lemma
\ref{sls-extension-lemma2} lies in the assertion that 
we can find frames for $\ell$ and $\tilde{\ell}$,
possibly alternative to (\ref{sls-lframe}) and 
(\ref{sls-rframe}), that are differentiable in 
$\lambda$ for a.e. $\lambda \in (-\infty, \kappa)$.
These lemmas are both proven under more general 
assumptions in Section 2.3 of \cite{HS2020b}. 
\end{remark}

Since $e^{D_- (\lambda) x}$ is an invertible matrix, we can replace the frame
$\mathbf{X} (x; \lambda)$ specified in Lemma \ref{sls-extension-lemma1} 
with $\mathbf{X} (x; \lambda) e^{D_- (\lambda) x}$
(i.e., each of these matrices is a valid frame for $\ell (x; \lambda)$).
From this observation and the estimates of Lemma \ref{sls-ode-lemma}
we see that 
\begin{equation*}
    \lim_{x \to -\infty} \mathbf{X} (x; \lambda) e^{D_- (\lambda) x}
    = \mathbf{X}_- (\lambda).
\end{equation*}
We conclude that the asymptotic Lagrangian subspace $\ell_- (\lambda)$
described in {\bf (B2)} exists, with the choice of frame $\mathbf{X}_- (\lambda)$.
Likewise, since $e^{ - D_+ (\lambda) x}$ is an invertible matrix, we can replace the frame
$\tilde{\mathbf{X}} (x; \lambda)$ specified in Lemma \ref{sls-extension-lemma1}  
with $\tilde{\mathbf{X}} (x; \lambda) e^{- D_+ (\lambda) x}$.
From this observation and the estimates of Lemma \ref{sls-ode-lemma}
we see that 
\begin{equation*}
    \lim_{x \to +\infty} \tilde{\mathbf{X}} (x; \lambda) e^{- D_+ (\lambda) x}
    = \tilde{\mathbf{X}}_+ (\lambda).
\end{equation*}
We conclude that the asymptotic Lagrangian subspace $\tilde{\ell}_+ (\lambda)$
described in {\bf (B2)} exists, with the choice of frame $\tilde{\mathbf{X}}_+ (\lambda)$.

We've now established that Assumptions {\bf (A)} and {\bf (B1)} hold, along with the 
first part of {\bf (B2)}. For the second part of {\bf (B2)}, we need to show that 
\begin{equation*}
    \ell_- (\lambda) \cap \tilde{\ell}_+ (\lambda) = \{0\},
    \quad \forall \, \lambda \in (-\infty, \kappa).
\end{equation*}
According to Lemma 2.2 of \cite{HS2019}, it suffices to show that 
the matrix $\mathbf{X}_- (\lambda)^* J \tilde{\mathbf{X}}_+ (\lambda)$
has a trivial kernel for all $\lambda \in (-\infty, \kappa)$. To 
verify this, we compute
\begin{equation} \label{sls-b2}
    \begin{aligned}
     \mathbf{X}_- (\lambda)^* J \tilde{\mathbf{X}}_+ (\lambda)
     &= (R_- (\lambda)^* \,\,\, - D_- (\lambda) R_- (\lambda)^* P_-^*)
     \begin{pmatrix}
     - P_+ R_+ (\lambda) D_+ (\lambda) \\
     R_+ (\lambda)
     \end{pmatrix} \\
     &= - R_- (\lambda)^* P_+ R_+ (\lambda) D_+ (\lambda)
     - D_- (\lambda)^* R_- (\lambda)^* P_-^* R_+ (\lambda). 
    \end{aligned}
\end{equation}
Here, $D_- (\lambda)$ and $P_{\pm}$ are self-adjoint. Also, by 
orthonormality with respect to the inner products $(\cdot, \cdot)_{\pm}$,
we have the relations
\begin{equation} \label{sls-orthonormality}
    R_{\pm} (\lambda)^* P_{\pm} R_{\pm} (\lambda) = I.
\end{equation}
Using these relations, we can compute 
\begin{equation*}
    (R_- (\lambda)^* P_+ R_+ (\lambda))^{-1}
    = (P_+ R_+ (\lambda))^{-1} (R_- (\lambda)^*)^{-1}
    = R_{+} (\lambda)^* P_- R_- (\lambda).
\end{equation*}
If we multiply (\ref{sls-b2}) on the left by 
$(R_- (\lambda)^* P_+ R_+ (\lambda))^{-1}$ we obtain 
\begin{equation*}
    D_+ (\lambda) 
    + R_{+} (\lambda)^* P_- R_- (\lambda) D_- (\lambda) R_- (\lambda)^* P_- R_+ (\lambda),
\end{equation*}
which is self-adjoint and negative definite (since the eigenvalues of the diagonal
matrices $D_{\pm} (\lambda)$ are all strictly negative). In particular, this 
matrix is non-singular, and we can conclude that 
$\mathbf{X}_- (\lambda)^* J \tilde{\mathbf{X}}_+ (\lambda)$ 
is non-singular as well, which is what we hoped to show. 

This leaves us with {\bf (B3)}, for which we first observe that 
\begin{equation*}
    \mathbb{B}_{\lambda} (x; \lambda)
    = \begin{pmatrix}
    Q(x) & 0 \\
    0 & 0
    \end{pmatrix}.
\end{equation*}
We see that 
\begin{equation*}
    \mathbf{X} (x; \lambda)^* \mathbb{B}_{\lambda} (x; \lambda) \mathbf{X} (x; \lambda)
    = X(x; \lambda)^* Q(x) X (x; \lambda),
\end{equation*}
so that 
\begin{equation*}
     \int_{-\infty}^c \mathbf{X} (x; \lambda)^* \mathbb{B}_{\lambda} (x; \lambda) \mathbf{X} (x; \lambda) dx
    = \int_{-\infty}^c X(x; \lambda)^* Q(x) X (x; \lambda) dx.
\end{equation*}
Since $Q (x)$ is positive definite for a.e. $x \in \mathbb{R}$, 
the right-hand side of this last relation is positive definite 
for any $c \in \mathbb{R}$, which is more than we need for 
{\bf (B3)}. We've now established that under our Assumptions 
{\bf (SL1)} and {\bf (SL2)} on (\ref{sls-again}), our general
Assumptions {\bf (A)}, {\bf (B1)}, {\bf (B2)}, and {\bf (B3)}
all hold. This establishes the first part of Theorem \ref{sls-theorem}.

\subsection{Exchanging the Target Space} \label{sls-exchange-section}

Next, we will use H\"ormander's index to show that the target spaces
$\tilde{\ell}_+ (\lambda_1)$ and $\tilde{\ell}_+ (\lambda_2)$ can 
be replaced by the Dirichlet plane $\ell_D$ with frame 
$\mathbf{X}_D = {0 \choose I}$, and that the resulting flow in this 
case is monotonic (i.e., the direction is the same for each conjugate point). 
For this argument, $\lambda_1$ and $\lambda_2$
are interchangeable, and we'll focus on the latter. 
As discussed in \cite{Ho2020} (see also \cite{Du1976, ZWZ2018}), the difference 
\begin{equation*}
    \mas (\ell (\cdot; \lambda_2), \tilde{\ell}_+ (\lambda_2); (-\infty, +\infty])
    -  \mas (\ell (\cdot; \lambda_2), \ell_D; (-\infty, +\infty])
\end{equation*}
depends only on the fixed Lagrangian subspaces $\ell_D$, $\tilde{\ell}_+ (\lambda_2)$,
$\ell_- (\lambda_2)$ and $\ell_+ (\lambda_2)$, allowing its specification as 
an index often referred to as H\"ormander's index
and denoted
\begin{equation} \label{sls-hormander}
    s (\ell_D, \tilde{\ell}_+ (\lambda_2); \ell_- (\lambda_2), \ell_+ (\lambda_2)). 
\end{equation}
In order to evaluate H\"ormander's index, we'll use the interpolation-space
approach of \cite{Ho2020}, writing 
\begin{equation} \label{sls-hormander}
\begin{aligned}
    s (\ell_D, &\tilde{\ell}_+ (\lambda_2); \ell_- (\lambda_2), \ell_+ (\lambda_2))
    = \mathcal{I} (\ell_+ (\lambda_2); \mathbf{X}_D, {I \choose \tilde{Y}_+ (\lambda_2) \tilde{X}_+ (\lambda_2)^{-1}}) \\
    & - \mathcal{I} (\ell_- (\lambda_2); \mathbf{X}_D, {I \choose \tilde{Y}_+ (\lambda_2) \tilde{X}_+ (\lambda_2)^{-1}}),
    \end{aligned}
\end{equation}
where in order to place ourselves in the context of \cite{Ho2020}, we have 
replaced $\tilde{\mathbf{X}}_+ (\lambda_2)$ with its normalized frame. 
(See \cite{Ho2020} for a full discussion of the notation $\mathcal{I} (\cdot; \cdot, \cdot)$;
for the current purposes equation (\ref{mathcal-i-formula}) just below 
is all we'll need.)
According to Section 3.3 of \cite{Ho2020}, we can write 
\begin{equation} \label{mathcal-i-formula}
\begin{aligned}
    \mathcal{I} &(\ell_+ (\lambda_2); \mathbf{X}_D, {I \choose \tilde{Y}_+ (\lambda_2) \tilde{X}_+ (\lambda_2)^{-1}})
    = n_- (\tilde{Y}_+ (\lambda_2) \tilde{X}_+ (\lambda_2)^{-1} - Y_+ (\lambda_2) X_+ (\lambda_2)^{-1}) \\
    &+ n_0 (\tilde{Y}_+ (\lambda_2) \tilde{X}_+ (\lambda_2)^{-1} - Y_+ (\lambda_2) X_+ (\lambda_2)^{-1}),
\end{aligned}
\end{equation}
where for any self-adjoint $n \times n$ matrix $M$, $n_- (M)$ denotes the number of
negative eigenvalues of $M$ and $n_0 (M)$ denotes the dimension of the kernel 
of $M$. 

In the current setting, 
\begin{equation*}
\tilde{Y}_+ (\lambda_2) \tilde{X}_+ (\lambda_2)^{-1}
= P_+ R_+ (\lambda_2) D_+ (\lambda_2) R_+ (\lambda_2)^* P_+.
\end{equation*}
In the event that $\lambda_2$ is not an eigenvalue of (\ref{sls-again}),
we must have $\mathbf{X}_+ (\lambda_2) = \tilde{\mathbf{X}}_+^g (\lambda_2)$
so that 
\begin{equation*}
{Y}_+ (\lambda_2) {X}_+ (\lambda_2)^{-1}
= - P_+ R_+ (\lambda_2) D_+ (\lambda_2) R_+ (\lambda_2)^* P_+.
\end{equation*}
In this way, we see that 
\begin{equation*}
    \tilde{Y}_+ (\lambda_2) \tilde{X}_+ (\lambda_2)^{-1} - Y_+ (\lambda_2) X_+ (\lambda_2)^{-1}
    = 2 P_+ R_+ (\lambda_2) D_+ (\lambda_2) R_+ (\lambda_2)^* P_+,
\end{equation*}
and this final matrix is self-adjoint and negative definite. We conclude that 
\begin{equation} \label{sls-hormander1}
    \mathcal{I} (\ell_+ (\lambda_2); \mathbf{X}_D, {I \choose \tilde{Y}_+ (\lambda_2) \tilde{X}_+ (\lambda_2)^{-1}})
    = n.
\end{equation}
Likewise, 
\begin{equation*}
\begin{aligned}
    \mathcal{I}& (\ell_- (\lambda_2); \mathbf{X}_D, {I \choose \tilde{Y}_+ (\lambda_2) \tilde{X}_+ (\lambda_2)^{-1}})
    = n_- (\tilde{Y}_+ (\lambda_2) \tilde{X}_+ (\lambda_2)^{-1} - Y_- (\lambda_2) X_- (\lambda_2)^{-1}) \\
    &+ n_0 (\tilde{Y}_+ (\lambda_2) \tilde{X}_+ (\lambda_2)^{-1} - Y_- (\lambda_2) X_- (\lambda_2)^{-1}).
\end{aligned}
\end{equation*}
In this case, 
\begin{equation*}
\begin{aligned}
    \tilde{Y}_+ &(\lambda_2) \tilde{X}_+ (\lambda_2)^{-1} - Y_- (\lambda_2) X_- (\lambda_2)^{-1}
    = P_+ R_+ (\lambda_2) D_+ (\lambda_2) R_+ (\lambda_2)^* P_+ \\
    & \quad + P_- R_- (\lambda_2) D_- (\lambda_2) R_- (\lambda_2)^* P_-. 
\end{aligned}
\end{equation*}
This is a sum of two self-adjoint negative definite operators, and so 
it is negative definite. We conclude that 
\begin{equation*}
     \mathcal{I} (\ell_- (\lambda_2); \mathbf{X}_D, {I \choose \tilde{Y}_+ (\lambda_2) \tilde{X}_+ (\lambda_2)^{-1}})
     = n,
\end{equation*}
and combining this with (\ref{sls-hormander1}), we see that 
\begin{equation*}
   s (\ell_D, \tilde{\ell}_+ (\lambda_2); \ell_- (\lambda_2), \ell_+ (\lambda_2)) = 0, 
\end{equation*}
so that 
\begin{equation*}
    \mas (\ell (\cdot; \lambda_2), \tilde{\ell}_+ (\lambda_2); (-\infty, +\infty])
    =  \mas (\ell (\cdot; \lambda_2), \ell_D; (-\infty, +\infty]).
\end{equation*}

Before turning to the case in which $\lambda_2$ is an eigenvalue 
of (\ref{sls-again}), we observe that conjugate points
arising in the calculation of $\mas (\ell (\cdot; \lambda_2), \ell_D; (-\infty, + \infty])$
all have the same sign (negative). To see this, we employ Lemma 1.1 
of \cite{Ho2020}, which asserts (in the current setting) that in 
order to conclude monotonicity, we need to check two things: (1)
If $P_D$ denotes projection onto the Dirichlet subspace, then
the matrix $P_D \mathbb{B} (x; \lambda_2) P_D$ is non-negative for a.e. 
$x \in \mathbb{R}$; and (2) if $y(x; \lambda_2)$ is any non-trivial 
solution of (\ref{sls-hammy}) with $y(x; \lambda_2) \in \ell_D$
for all $x$ in some interval $[a, b]$, $a < b$, then 
\begin{equation*}
    \int_a^b (\mathbb{B} (x; \lambda_2) y (x; \lambda_2), y(x; \lambda_2)) dx > 0.
\end{equation*}
For (1), we observe that for any $v = {v_1 \choose v_2} \in \mathbb{C}^{2n}$, we have 
$P_D v = {0 \choose v_2}$, so that 
\begin{equation*}
    v^* P_D \mathbb{B} (x; \lambda_2) P_D v
    = (0 \,\,\, v_2^*) 
    \begin{pmatrix}
    \lambda_2 Q(x) - V(x) & 0 \\
    0 & P(x)^{-1}
    \end{pmatrix}
    {0 \choose v_2} 
    = v_2^* P(x)^{-1} v_2 \ge 0,
\end{equation*}
where the final inequality follows from Assumption {\bf (SL2)}. For (2), 
suppose $y(x; \lambda_2)$ is any non-trivial solution of (\ref{sls-hammy})
so that $y(x; \lambda_2) \in \ell_D$ for all $x$ in some interval $[a, b]$, $a < b$.
Then, in particular, $\phi (x; \lambda_2) = 0$ for all such $x$, and 
since $\phi (x; \lambda_2)$ is absolutely continuous on $\mathbb{R}$
we can conclude that $\phi' (x; \lambda_2) = 0$ for a.e. $x \in (a, b)$.
But then $y(x; \lambda_2) = 0$ for a.e. $x \in (a, b)$, contradicting our 
assumption that $y(x; \lambda_2)$ is non-trivial. We conclude that 
Items (1) and (2) both hold, and from Lemma 1.1 of \cite{Ho2020}
we can conclude that conjugate points
arising in the calculation of $\mas (\ell (\cdot; \lambda_2), \ell_D; (-\infty, + \infty])$
all have the same sign (negative). If $\lambda_2$ is not an eigenvalue of 
(\ref{sls-again}), we can now write 
\begin{equation} \label{sls-kernel-sum}
    \mas (\ell (\cdot; \lambda_2), \ell_D; (-\infty, + \infty])
    = - \sum_{x \in \mathbb{R}} \dim (\ell (x; \lambda_2) \cap \ell_D).
\end{equation}
In writing this relation, we've taken advantage of the observations
that $\ell_- (\lambda_2) \cap \ell_D = \{0\}$ and that since 
we are currently assuming that $\lambda_2$ is not an eigenvalue 
of (\ref{sls-again}), $\ell_+ (\lambda_2) \cap \ell_D = \{0\}$
(because $\ell_+ (\lambda_2) = \tilde{\ell}_+^g (\lambda_2)$, and 
$\tilde{\ell}_+^g (\lambda_2) \cap \ell_D = \{0\}$).
(These claims are easily checked by using the frames for $\ell_- (\lambda_2)$
and $\tilde{\ell}_+^g (\lambda_2)$.) We conclude that the left-hand
side of (\ref{sls-kernel-sum}) can be replaced by the Maslov index
$\mas (\ell (\cdot; \lambda_2), \ell_D; [-L, + L])$ for any $L$ 
sufficiently large, and correspondingly the right-hand side can 
be replaced by $- \sum_{x \in (-L, L)} \dim (\ell (x; \lambda_2) \cap \ell_D)$.
In addition, it's clear from 
monotonicity that the values of $x$ for which  
$\dim (\ell (x; \lambda_2) \cap \ell_D) \ne 0$ form a discrete
set, so the 
right-hand side of (\ref{sls-kernel-sum}) is a finite sum. 

In the event that $\lambda_2$ is an eigenvalue for (\ref{sls-again}),
we have to take a different approach, because we no longer have 
an explicit expression for the frame $\mathbf{X}_+ (\lambda_2)$
(which now comprises a combination of solutions that lie right 
in $\mathbb{R}$ and solutions that do not lie right in $\mathbb{R}$).
Nonetheless, according to Section 3.1 of \cite{Ho2020} we always 
have 
\begin{equation*}
     -n \le \mathcal{I} (\ell_- (\lambda_2); \mathbf{X}_D, {I \choose \tilde{Y}_+ (\lambda_2) \tilde{X}_+ (\lambda_2)^{-1}}) \le n.
\end{equation*}
In this way, we see that 
\begin{equation*}
    s (\ell_D, \tilde{\ell}_+ (\lambda_2); \ell_- (\lambda_2), \ell_+ (\lambda_2))    
    = \mathcal{I} (\ell_- (\lambda_2); \mathbf{X}_D, {I \choose \tilde{Y}_+ (\lambda_2) \tilde{X}_+ (\lambda_2)^{-1}}) - n \\
    \le 0,
\end{equation*}
from which it follows that 
\begin{equation*}
    \mas (\ell (\cdot; \lambda_2), \tilde{\ell}_+ (\lambda_2); (-\infty, +\infty])
    \le  \mas (\ell (\cdot; \lambda_2), \ell_D; (-\infty, +\infty]).
\end{equation*}
Using our observation that Theorem \ref{main-theorem} applies in this 
case, we conclude that 
\begin{equation} \label{previous-inequality}
    \mathcal{N} ([\lambda_1, \lambda_2)) 
    \ge -  \mas (\ell (\cdot; \lambda_2), \ell_D; (-\infty, +\infty])
    +  \mas (\ell (\cdot; \lambda_1), \tilde{\ell}_+ (\lambda_1); (-\infty, +\infty]).
\end{equation}
(Here, the Maslov index at $\lambda_1$ has been left unchanged.)

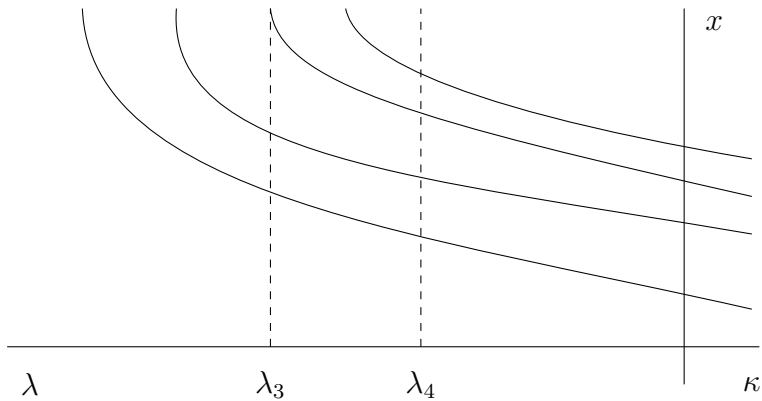
\begin{figure}[ht]
\begin{center}
\begin{tikzpicture}
\draw (-5,0) -- (5,0);	%horizontal line
\draw (4,-.5) -- (4,4.5);	%vertical line
\node at (4.4,4.3) {$x$};
\node at (-4.7,-.5) {$\lambda$};
\node at (-1.5,-.5) {$\lambda_3$};
\node at (.5,-.5) {$\lambda_4$};
\node at (4.9,-.5) {$\kappa$};
\draw[dashed] (-1.5,0) -- (-1.5,4.5);	%vertical line
\draw[dashed] (.5,0) -- (.5,4.5);	%vertical line
%
%SPECTRAL CURVES
\draw (-4,4.5) .. controls (-3.9,2) and (.5, 1.5) .. (4.9,.5);
\draw (-2.75,4.5) .. controls (-2.9,2.5) and (.5, 2.25) .. (4.9,1.5);
\draw (-1.5,4.5) .. controls (-1.4,3.5) and (.5, 3) .. (4.9,2.0);
\draw (-.5,4.5) .. controls (-.4,4) and (.5, 3.25) .. (4.9,2.5);
\end{tikzpicture}
\end{center}
\caption{Monotonic spectral curves when the target space is $\ell_D$. \label{D}}
\end{figure}

In order to obtain the opposite inequality, we observe that 
monotonicity implies that for any pair $\lambda_3 < \lambda_4 < \kappa$ 
the number of conjugate points along the shelf at $\lambda_3$ will be 
less than or equal to the number of conjugate points along the shelf 
at $\lambda_4$. (See Figure \ref{D} and the detailed discussion in 
Section 2.3 of \cite{Ho2020}.) In particular, for 
any $\epsilon > 0$, 
\begin{equation*}
     - \mas (\ell (\cdot; \lambda_2), \ell_D; (-\infty, +\infty])
    \ge - \mas (\ell (\cdot; \lambda_2 - \epsilon), \ell_D; (-\infty, +\infty])
\end{equation*}
(keeping in mind that by monotonicity each side of this inequality is 
non-negative). 
Since the eigenvalues of (\ref{sls-again}) are discrete, we can take $\epsilon$
small enough so that $\lambda_2 - \epsilon$ is not an eigenvalue of (\ref{sls-again}),
and moreover (\ref{sls-again}) has no eigenvalues on the interval 
$[\lambda_2 - \epsilon, \lambda_2)$. This allows us to write
\begin{equation*}
    \begin{aligned}
        \mathcal{N} ([\lambda_1, \lambda_2)) &= \mathcal{N} ([\lambda_1, \lambda_2 - \epsilon)) \\
        &= - \mas (\ell (\cdot; \lambda_2 - \epsilon), \ell_D; (-\infty, +\infty])
        +  \mas (\ell (\cdot; \lambda_1), \tilde{\ell}_+ (\lambda_1); (-\infty, +\infty]) \\
        &\le  - \mas (\ell (\cdot; \lambda_2), \ell_D; (-\infty, +\infty])
        +  \mas (\ell (\cdot; \lambda_1), \tilde{\ell}_+ (\lambda_1); (-\infty, +\infty]).
    \end{aligned}
\end{equation*}
Combining these observations with our previous inequality (\ref{previous-inequality}), we conclude that 
\begin{equation*}
    \mathcal{N} ([\lambda_1, \lambda_2))
    = - \mas (\ell (\cdot; \lambda_1), \ell_D; (-\infty, +\infty])
        +  \mas (\ell (\cdot; \lambda_1), \tilde{\ell}_+ (\lambda_1); (-\infty, +\infty]),
\end{equation*}
and consequently  
\begin{equation*}
    \mas (\ell (\cdot; \lambda_2), \tilde{\ell}_+ (\lambda_2); (-\infty, +\infty])
    =  \mas (\ell (\cdot; \lambda_1), \ell_D; (-\infty, +\infty]).
\end{equation*}

%\begin{figure}[ht] 
%\begin{center}\includegraphics[%
%  width=6.5cm,
%  height=5cm]{MequalsMonRfigD.eps}
%\end{center}
%\caption{Monotonic spectral curves when the target space is %$\ell_D$. \label{D}}
%\end{figure}

As in the case in which $\lambda_2$ is not an eigenvalue of (\ref{sls-again}),
We continue to have $\ell_- (\lambda_2) \cap \ell_D = \{0\}$, 
but if $\lambda_2$ is an eigenvalue of (\ref{sls-again}) we 
may have $\ell_+ (\lambda_2) \cap \ell_D \ne \{0\}$. Nonetheless, 
by monotonicity the direction associated with $+\infty$ as a 
conjugate point must be clockwise, and so an arrival does not
affect the Maslov index. This allows us to conclude that 
(\ref{sls-kernel-sum}) continues to hold even when $\lambda_2$
is an eigenvalue of (\ref{sls-again}). The same considerations
hold for $\lambda_1$, allowing us to write 
\begin{equation*}
    \begin{aligned}
        \mathcal{N} ([\lambda_1, \lambda_2))
        &= \sum_{x \in \mathbb{R}} \dim (\ell (x; \lambda_2) \cap \ell_D)
        - \sum_{x \in \mathbb{R}} \dim (\ell (x; \lambda_1) \cap \ell_D) \\
        &= \sum_{x \in \mathbb{R}} \dim \ker (\mathbf{X} (x; \lambda_2)^* J \mathbf{X}_D)
        - \sum_{x \in \mathbb{R}} \dim \ker (\mathbf{X} (x; \lambda_1)^* J \mathbf{X}_D), 
    \end{aligned}
\end{equation*}
where in obtaining this second equality, we have observed from Lemma 2.2 of 
\cite{HS2019} that if $\mathbf{X}_1$ and $\mathbf{X}_2$ are frames for any 
two Lagrangian subspaces $\ell_1$ and $\ell_2$ then 
\begin{equation*}
    \dim (\ell_1 \cap \ell_2) = \dim \ker (\mathbf{X}_1^* J \mathbf{X}_2).
\end{equation*}
For these latter calculations, we have (recalling 
$\mathbf{X} (x; \lambda_2) = {X (x; \lambda_2) \choose Y(x; \lambda_2)}$)
\begin{equation*}
    \mathbf{X} (x; \lambda_2)^* J \mathbf{X}_D = - X (x; \lambda_2)^*,
\end{equation*}
and since $\dim \ker (- X (x; \lambda_2)^*) = \dim \ker (X(x; \lambda_2))$,
we can write 
\begin{equation} \label{sls-kernel-sums}
     \mathcal{N} ([\lambda_1, \lambda_2))
     = \sum_{x \in \mathbb{R}} \dim \ker X (x; \lambda_2)
        - \sum_{x \in \mathbb{R}} \dim \ker X (x; \lambda_1),
\end{equation}
which is precisely the second assertion in Theorem \ref{sls-theorem}.

\subsection{Eliminating the Left Shelf} \label{sls-elimination-section}

In this section, we will check that we can take $\lambda_1$ sufficiently negative
so that there are no conjugate points along the left shelf. We begin by 
noting that a point $(s, \lambda_1) \in \mathbb{R} \times (-\infty, \kappa)$
will be a conjugate point for the Lagrangian subspaces $\ell (x; \lambda)$
and $\tilde{\ell}_+ (\lambda_1)$ if and only if $\lambda_1$ is an eigenvalue 
for the half-line problem 
\begin{equation*}
    \begin{aligned}
        - (P(x)\phi')' &+ V(x) \phi = \lambda Q(x) \phi;
        \quad x \in (-\infty, s) \\
        &\tilde{\mathbf{X}}_+ (\lambda_1)^* J {\phi(s) \choose P(s) \phi(s)} = 0.
    \end{aligned}
\end{equation*}
We will use an energy argument to show that the set of eigenvalues for this 
problem is bounded below, independently of $s$. To this end, suppose $\lambda$ is an eigenvalue, 
and let $\phi(x; \lambda)$ denote a corresponding eigenfunction. If we take
an $L^2 ((-\infty, s), \mathbb{C}^n)$ inner product of the system with $\phi$, we obtain 
\begin{equation*}
    -\int_{-\infty}^s ((P(x) \phi')',\phi) dx
    + \int_{-\infty}^s (V(x) \phi, \phi) dx
    = \lambda \int_{-\infty}^s (Q(x) \phi, \phi) dx.
\end{equation*}
For the first integral, we can integrate by parts to 
write 
\begin{equation*}
     -\int_{-\infty}^s ((P(x) \phi')',\phi) dx
     = - (P(s) \phi' (s), \phi (s)) + \int_{-\infty}^s (P(x) \phi', \phi') dx.
\end{equation*}
The key point here is the boundary term, and for this, we observe 
that our boundary condition can be expressed as 
\begin{equation*}
    \begin{aligned}
        0 &= (R_+ (\lambda_1)^* \,\,\, D_+ (\lambda_1) R_+ (\lambda_1)^* P_+) 
        {- P(s) \phi' (s) \choose \phi(s)} \\
        & = - R_+ (\lambda_1)^* P(s) \phi' (s) 
        + D_+ (\lambda_1) R_+ (\lambda_1)^* P_+ \phi (s). 
    \end{aligned}
\end{equation*}
Recalling the relation $(R_+ (\lambda_1)^*)^{-1} = P_+ R_+ (\lambda_1)$, 
we can solve for $P(s) \phi' (s)$ in terms of $\phi(s)$ to get
\begin{equation*}
  P(s) \phi' (s) = P_+ R_+ (\lambda_1) D_+ (\lambda_1) R_+ (\lambda_1)^* P_+ \phi (s).  
\end{equation*}
We see that the boundary term can be expressed as 
\begin{equation*}
    - (P(s) \phi' (s), \phi (s)) 
    = - (P_+ R_+ (\lambda_1) D_+ (\lambda_1) R_+ (\lambda_1)^* P_+ \phi (s), \phi (s)).
\end{equation*}
The matrix $P_+ R_+ (\lambda_1) D_+ (\lambda_1) R_+ (\lambda_1)^* P_+$ is negative 
definite, so we can conclude that 
\begin{equation*}
- (P(s) \phi' (s), \phi (s)) \ge 0    
\end{equation*}
for all $\phi (s) \in \mathbb{C}^n$. It follows that 
\begin{equation*}
     -\int_{-\infty}^s ((P(x) \phi')',\phi) dx
     \ge \int_{-\infty}^s (P(x) \phi', \phi') dx
     \ge \theta_P \|\phi'\|_{L^2 ((-\infty, s), \mathbb{C}^n)}^2.
\end{equation*}
We also have 
\begin{equation*}
    \int_{-\infty}^s (Q(x) \phi, \phi) dx
    \ge \theta_Q \|\phi\|_{L^2 ((-\infty, s), \mathbb{C}^n)}^2,
\end{equation*}
and 
\begin{equation*}
    \Big| \int_{-\infty}^s (V(x) \phi, \phi) dx \Big|
    \le C_V \|\phi\|_{L^2 ((-\infty, s), \mathbb{C}^n)}^2.
\end{equation*}
(We emphasize that this is the only place in the analysis in which 
we use the boundedness of $V$ assumed in {\bf (SL1)}.) 
Combining these observations, we see that for $\lambda < 0$ 
we can write 
\begin{equation*}
    \begin{aligned}
        \lambda \theta_Q \|\phi\|_{L^2 ((-\infty, s), \mathbb{C}^n)}^2
        &\ge \lambda \int_{-\infty}^s (Q(x) \phi, \phi) dx
        = -\int_{-\infty}^s ((P(x) \phi')',\phi) dx
    + \int_{-\infty}^s (V(x) \phi, \phi) dx \\
    &\ge \theta_P \|\phi'\|_{L^2 ((-\infty, s), \mathbb{C}^n)}^2
    - C_V \|\phi\|_{L^2 ((-\infty, s), \mathbb{C}^n)}^2,
    \end{aligned}
\end{equation*}
from which we see that 
\begin{equation*}
    \lambda \ge - \frac{C_V}{\theta_Q}.
\end{equation*}
We conclude that for $\lambda < - (C_V/\theta_Q)$ there are no 
conjugate points. In particular, if $\lambda_1 < - (C_V/\theta_Q)$,
then 
\begin{equation*}
    \sum_{x \in \mathbb{R}} \dim \ker (X (x; \lambda_1)) = 0,
\end{equation*}
so by taking $\lambda_1 < - (C_V/\theta_Q)$ in (\ref{sls-kernel-sums})
we see that 
\begin{equation*}
  \mathcal{N} ((-\infty, \lambda_2))
  = \sum_{x \in \mathbb{R}} \dim \ker X (x; \lambda_2), 
\end{equation*}
which is precisely the final claim of Theorem \ref{sls-theorem}.
\hfill $\square$

\subsection{Traveling Waves} \label{traveling-waves-section}

As noted in the Introduction, if we want to analyze the stability
of a traveling-wave solution $\bar{u} (x - st)$ to the Allen-Cahn 
equation (\ref{allen-cahn-equation}), we need to understand
the eigenvalues of 
\begin{equation} \label{sls-s}
    H_s \phi := - \phi'' - s \phi' + V(x) \phi = \lambda \phi,
    \quad V(x) = D^2F (\bar{u} (x)),
\end{equation}
which is not self-adjoint for $s \ne 0$ (even if $V(x)$ is self-adjoint).
If we set $y = {y_1 \choose y_2} = {\phi \choose \phi'}$, then 
we obtain 
\begin{equation*}
    J y' = \mathcal{B} (x; \lambda) y;
    \quad \mathcal{B} (x; \lambda)
    = \begin{pmatrix}
    \lambda I - V(x) & s I \\
    0 & I
    \end{pmatrix},
\end{equation*}
where we're using $\mathcal{B} (x; \lambda)$ in order to reserve 
the notation $\mathbb{B} (x; \lambda)$ for self-adjoint matrices.

In this case, we can readily place the analysis in the setting of 
(\ref{hammy}) by making the change of variables $\zeta = e^{\frac{s}{2} x} y$,
for which we find 
\begin{equation} \label{sls-s-hammy}
    J \zeta' = \mathbb{B} (x; \lambda) \zeta;
    \quad \mathbb{B} (x; \lambda)
    = \begin{pmatrix}
    \lambda I - V(x) & \frac{s}{2} I \\
    \frac{s}{2} & I
    \end{pmatrix}.
\end{equation}
If $V$ satisfies the same assumptions as stated in 
{\bf (SL1)} and {\bf (SL2)}, then our analysis of 
(\ref{sls-hammy}) can be carried out with only 
minor adjustments, and we can conclude precisely 
the claims stated in Theorem \ref{sls-theorem}. 
In fact, as shown in \cite{HLS2018}, the limit conditions 
can be relaxed from exponential rate to the 
following.

\medskip
{\bf (SL2)$^\prime$} There exist self-adjoint matrices
$V_{\pm}$ so that the limits $\lim_{x \to \pm \infty} V(x) = V_{\pm}$ 
exist, and for each $M \in \mathbb{R}$, 
\begin{equation*}
    \int_{M}^{+\infty} (1+|x|) |V(x) - V_+| dx < \infty,
    \quad \int_{-\infty}^{M} (1+|x|) |V(x) - V_-| dx < \infty.
\end{equation*}
\medskip

For convenient reference, we state this assertion as a 
theorem. For a full proof, though by different calculations 
in some places, the reader is referred to \cite{HLS2018}. 

\begin{theorem} \label{sls-theorem-s}
For (\ref{sls-s}), let Assumptions {\bf (SL1)} (on $V$) and {\bf (SL2)$^{\mathbf{\prime}}$} 
hold. Then $\sigma_p (H_s) \subset \mathbb{R}$,
and for $\kappa$ specified as in (\ref{kappa-defined}),
except with $Q_{\pm} = I$,  {\bf (A)},
{\bf (B1)}, {\bf (B2)}, and {\bf (B3)} all hold for 
(\ref{sls-s-hammy}) with $I = (- \infty, \kappa)$. We
conclude that the result of Theorem \ref{main-theorem} holds for 
all intervals $[\lambda_1, \lambda_2]$, $\lambda_1 < \lambda_2 < \kappa$. 
In addition, if $\mathcal{N} ([\lambda_1, \lambda_2))$ 
denotes the number of eigenvalues, counted with multiplicity, 
that (\ref{sls-s}) has on the interval $[\lambda_1, \lambda_2)$,
and we express the frame $\mathbf{X} (x; \lambda)$ from 
{\bf (B1)} as $\mathbf{X} (x; \lambda) = {X(x; \lambda) \choose Y(x; \lambda)}$,
then 
\begin{equation*}
    \mathcal{N} ([\lambda_1, \lambda_2))
    = \sum_{x \in \mathbb{R}} \dim \ker X(x; \lambda_2)
    - \sum_{x \in \mathbb{R}} \dim \ker X(x; \lambda_1),
\end{equation*}
and 
\begin{equation*}
    \mathcal{N} ((-\infty, \lambda_2)) 
    = \sum_{x \in \mathbb{R}} \dim \ker X(x; \lambda_2). 
\end{equation*}
\end{theorem}

Unfortunately, our approach to handling traveling waves 
$\bar{u} (x - st)$ does not readily extend to more
general Allen-Cahn type systems such as 
\begin{equation*}
    u_t + DF(u) = B u_{xx},
    \quad (x, t) \in \mathbb{R} \times \mathbb{R}_+,
    \quad u(x, t) \in \mathbb{C}^n,
\end{equation*}
for which the diffusion matrix $B$ is not the identity
matrix. For an interesting step in this direction, we 
refer the reader to the recent result \cite{BCCJM2020}.

We conclude this section by mentioning a second, more complicated,
case in which the current method can be applied in the analysis
of traveling waves. In particular, we consider equations
\begin{equation} \label{allen-cahn-m}
    u_t + M DF(u) = u_{xx},
    \quad (x, t) \in \mathbb{R} \times \mathbb{R}_+,
    \quad u(x, t) \in \mathbb{C}^n,
\end{equation}
for which $M$ denotes a constant, invertible, self-adjoint
$n \times n$ matrix. In order to analyze the stability of 
a traveling-wave solution $\bar{u} (x - st)$ to (\ref{allen-cahn-m}),
we use moving coordinates as before and linearize, leading
to the eigenvalue problem 
\begin{equation} \label{allen-cahn-m-e}
    - \phi'' - s \phi' + M D^2F(\bar{u} (x)) = \lambda \phi,
    \quad (x, t) \in \mathbb{R} \times \mathbb{R}_+,
    \quad u(x, t) \in \mathbb{C}^n.
\end{equation}
In order to place (\ref{allen-cahn-m-e})
in the current framework, we set $y = {y_1 \choose y_2} = {\phi \choose M^{-1} \phi'}$,
so that 
\begin{equation*}
    Jy' = \mathcal{B} (x; \lambda) y;
    \quad \mathcal{B} (x; \lambda) =
    \begin{pmatrix}
    \lambda M^{-1} - D^2F(\bar{u} (x)) & sI \\
    0 & M
    \end{pmatrix}. 
\end{equation*}
If we now set $\zeta = e^{\frac{s}{2} x} y$ as before, we obtain 
the system 
\begin{equation*}
    J\zeta' = \mathbb{B} (x; \lambda) \zeta;
    \quad \mathbb{B} (x; \lambda) =
    \begin{pmatrix}
    \lambda M^{-1} - D^2 F(\bar{u} (x)) & \frac{s}{2} I \\
    \frac{s}{2}I & M
    \end{pmatrix},
\end{equation*}
which has the form of (\ref{hammy}). 

The verification of our general assumptions {\bf (A)}, 
{\bf (B1)},  {\bf (B2)}, and  {\bf (B3)} requires 
additional assumptions on $M$, and we won't pursue a
full analysis here. We note, however, that the following
important case was analyzed in \cite{ CJ2020, Co2019}:
$M = QS$, where $Q$ is a diagonal matrix with either 
$+1$ or $-1$ in each diagonal entry and $S$ is a positive
diagonal matrix.

\section{Differential-Algebraic Sturm-Liouville Systems} \label{sls-a-section}

Following Section 5.4 in \cite{HS2019}, we consider differential-algebraic 
Sturm-Liouville systems 
\begin{equation} \label{sls-a-again}
    \mathcal{L}_a \phi = - (P(x)\phi')' + V(x) \phi = \lambda \phi,
    \quad x \in \mathbb{R}, \quad \phi (x) \in \mathbb{C}^n,
\end{equation}
with degenerate matrix 
\begin{equation*}
    P(x) = \begin{pmatrix}
    P_{11} (x) & 0 \\
    0 & 0
    \end{pmatrix}.
\end{equation*}
We make the following assumptions on $P$ and $V$:

\medskip
{\bf (DA1)} For some $0 < m < n$, $P_{11} \in \AC_{\loc} (\mathbb{R}, \mathbb{C}^{m \times m})$,
with $P_{11} (x)$ self-adjoint for all $x \in \mathbb{R}$; also, 
$V \in C(\mathbb{R}, \mathbb{C}^{n \times n})$, with $V(x)$ self-adjoint
for all $x \in \mathbb{R}$. In addition, there 
exists a constant $\theta_{P_{11}} > 0$ so that for any $v \in \mathbb{C}^n$,
\begin{equation*}
    (P_{11} (x)v, v) \ge \theta_{P_{11}} |v|^2
\end{equation*}
for all $x \in \mathbb{R}$.

\medskip
{\bf (DA2)} There exist self-adjoint matrices $P_{11}^{\pm}, V^{\pm}$, along
with constants $C$, $M$, and $\eta > 0$, so that 
\begin{equation*}
\begin{aligned}
    |P_{11} (x) - P_{11}^{\pm}| \le C e^{- \eta |x|}, 
    \quad x \gtrless \pm M;
    &\quad \quad |P_{11}' (x)| \le C e^{- \eta |x|}, 
    \quad x \gtrless \pm M; \\
    |V (x) - V^{\pm}| \le C e^{- \eta |x|}, 
    \quad x \gtrless \pm M.
\end{aligned}
\end{equation*}

\begin{remark}
We note that in this case boundedness of $V$ follows from {\bf (DA1)}
and {\bf (DA2)}. 
\end{remark}

For notational convenience, we'll write 
\begin{equation*}
    V(x) = 
    \begin{pmatrix}
    V_{11} (x) & V_{12} (x) \\
    V_{12} (x)^* & V_{22} (x)
    \end{pmatrix},
\end{equation*}
where $V_{11} (x)$ is an $m \times m$ matrix, 
$V_{12} (x)$ is an $m \times (n-m)$ matrix,
and $V_{22} (x)$ is an $(n - m) \times (n - m)$
matrix. We'll write $\phi = {\phi_1 \choose \phi_2}$,
where $\phi_1 (x; \lambda) \in \mathbb{C}^m$ and 
$\phi_2 (x; \lambda) \in \mathbb{C}^{n - m}$, 
allowing us to express (\ref{sls-a-again}) as 
\begin{equation} \label{sls-a-components}
    \begin{aligned}
        - (P_{11} (x) \phi_1')' + V_{11} (x) \phi_1 + V_{12} (x) \phi_2 &= \lambda \phi_1 \\
        V_{12} (x)^* \phi_1 + V_{22} (x) \phi_2 &= \lambda \phi_2.
    \end{aligned}
\end{equation}
We will take as our domain for $\mathcal{L}_a$, the set 
\begin{equation*}
    \begin{aligned}
        \mathcal{D}_a &:= \{\phi = (\phi_1, \phi_2) \in L^2 (\mathbb{R}, \mathbb{C}^m) \times L^2 (\mathbb{R}, \mathbb{C}^{n-m}): \\
        &\phi_1, \phi_1' \in \AC_{\loc} (\mathbb{R}, \mathbb{C}^m), \, 
        \mathcal{L}_a \phi \in L^2 (\mathbb{R}, \mathbb{C}^m) \times L^2 (\mathbb{R}, \mathbb{C}^{n-m}) \}.
    \end{aligned}
\end{equation*}

According to \cite{ALMS1994}, the essential spectrum of $\mathcal{L}_a$ 
will contain the ranges of the eigenvalues of $V_{22} (x)$ as 
$x$ ranges over $\mathbb{R}$. Precisely, we'll let $\{\nu_k (x)\}_{k=1}^{n-m}$
denote the eigenvalues of $V_{22} (x)$, and we'll denote by 
$\mathcal{R}_k$ the closure of the range of $\nu_k: \mathbb{R} \to \mathbb{R}$.
Then 
\begin{equation*}
    \bigcup_{k=1}^{n-m} \mathcal{R}_k \subset \sigma_{\ess} (\mathcal{L}_a).
\end{equation*}
For all $\lambda \notin \cup_{k=1}^{n-m} \mathcal{R}_k$, we can solve the 
second equation in (\ref{sls-a-components}) for $\phi_2$, giving
\begin{equation*}
    \phi_2 (x; \lambda) 
    = (\lambda I - V_{22} (x))^{-1} V_{12} (x)^* \phi_1 (x; \lambda).
\end{equation*}
Upon substitution of this expression for $\phi_2$ into the first 
equation in (\ref{sls-a-components}), we obtain an equation 
for $\phi_1$, 
\begin{equation} \label{phi1-equation}
    - (P_{11} (x) \phi_1')' + \mathbf{V} (x; \lambda) \phi_1 
    = \lambda \phi_1,
\end{equation}
where we've set
\begin{equation*}
    \mathbf{V} (x; \lambda) 
    := V_{11} (x) + V_{12} (x) (\lambda I - V_{22} (x))^{-1} V_{12} (x)^*.
\end{equation*} 

We can now analyze (\ref{phi1-equation}) similarly as we 
analyzed (\ref{sls-again}). First, for $\lambda \notin \cup_{k=1}^{n-m} \mathcal{R}_k$, the 
matrices $(\lambda I - V_{22}^{\pm})$ are non-singular, and 
we can consider the limiting system 
\begin{equation} \label{sls-a-limit}
    - P_{11}^{\pm} \phi_1'' + \mathbf{V}^{\pm} (\lambda) \phi_1
    = \lambda \phi_1,
\end{equation}
where 
\begin{equation*}
    \mathbf{V}^{\pm} (\lambda)
    := V_{11}^{\pm} + V_{12}^{\pm} (\lambda I - V_{22}^{\pm})^{-1} V_{12}^{\pm \, *}.
\end{equation*}
Similarly as with (\ref{sls-again}), we can check that in 
addition to the set $\cup_{k=1}^{n-m} \mathcal{R}_k$, the 
essential spectrum of $\mathcal{L}_a$ includes all values 
$\lambda$ for which $\phi_1 (x) = e^{ikx} r_1$ solves this 
equation for some constant scalar $k \in \mathbb{R}$ and 
constant vector $r_1 \in \mathbb{C}^m$. In this case, 
we have 
\begin{equation*}
    \{k^2 P_{11}^{\pm} +  \mathbf{V}^{\pm} (\lambda) \}r_1 
    = \lambda r_1. 
\end{equation*}
Computing an inner product of this system with $r_1$, 
we see that 
\begin{equation*}
    k^2 (P_{11}^{\pm} r_1, r_1) + (\mathbf{V}^{\pm} (\lambda) r_1, r_1)
    = \lambda |r_1|^2.
\end{equation*}
Since the matrices $P_{11}^{\pm}$ are positive definite, we 
see that in order for $\lambda$ to satisfy this relationship,
we must have 
\begin{equation*}
    \lambda \ge \frac{(\mathbf{V}^{\pm} (\lambda) r_1, r_1)}{|r_1|^2}.
\end{equation*}
If we set 
\begin{equation*}
    \kappa^{\pm} (\lambda) := \inf_{r_1 \in \mathbb{C}^m \backslash \{0\}} 
    \frac{(\mathbf{V}^{\pm} (\lambda) r_1, r_1) }{|r_1|^2}
\end{equation*}
(i.e., the smallest eigenvalue of he matrix $\mathbf{V}^{\pm} (\lambda)$),
then in order to avoid this part of the essential spectrum we must
take $\lambda$ so that 
\begin{equation*}
    \lambda \in \mathcal{R}_0 := \Big{\{}\lambda \in \mathbb{R}: 
    \lambda < \min \{\kappa^- (\lambda), \kappa^+ (\lambda) \} \Big{\}}. 
\end{equation*}
With this notation in place, 
we see that we can consider any interval 
$I \subset \mathbb{R}$ so that 
\begin{equation} \label{sls-a-i}
    I \cap \bigcup_{k=0}^{n-m} \mathcal{R}_k = \emptyset.
\end{equation}

As an important example case, we observe that we can take any interval $I$ that lies 
entirely below the essential spectrum. In order to 
characterize the bottom of the essential spectrum 
more precisely, we begin by setting
\begin{equation*}
    \kappa_1^{\pm} := \inf_{r_1 \in \mathbb{C}^m \backslash \{0\}} 
    \frac{(V_{11}^{\pm} r_1, r_1)}{|r_1|^2}
\end{equation*}
and 
\begin{equation*}
    \kappa_2^{\pm} := \inf_{r_2 \in \mathbb{C}^{n-m} \backslash \{0\}} 
    \frac{(V_{22}^{\pm} r_2, r_2)}{|r_2|^2}.
\end{equation*}
By spectral mapping, the eigenvalues of $(\lambda I - V_{22}^{\pm})^{-1}$
will be $(\lambda - \nu_k^{\pm})^{-1}$, where $\{\nu_k^{\pm}\}_{k=1}^{n-m}$
denote the eigenvalues of $V_{22}^{\pm}$. In this case, we're taking 
$\lambda$ below the set $\cup_{k=1}^{n-m} \mathcal{R}_k$, so in 
particular below the eigenvalues of $V_{22}^{\pm}$. It follows that 
\begin{equation*}
    \inf_{r_2 \in \mathbb{C}^{n-m} \backslash \{0\}} 
    \frac{((\lambda I - V_{22}^{\pm})^{-1}r_2, r_2)}{|r_2|^2}
    = (\lambda - \kappa_2^{\pm})^{-1} < 0,
\end{equation*}
and so 
\begin{equation*}
 ((\lambda I - V_{22}^{\pm})^{-1}r_2, r_2) 
 \ge (\lambda - \kappa_2^{\pm})^{-1} |r_2|^2.    
\end{equation*}
This allows us to compute
\begin{equation*}
    \begin{aligned}
     \inf_{r_1 \in \mathbb{C}^m \backslash \{0\}} &
    \frac{(V_{12}^{\pm} (\lambda I - V_{22}^{\pm})^{-1} V_{12}^{\pm \, *} r_1, r_1)}{|r_1|^2}
    = \inf_{r_1 \in \mathbb{C}^m \backslash \{0\}} 
    \frac{((\lambda I - V_{22}^{\pm})^{-1} V_{12}^{\pm \, *} r_1, V_{12}^{\pm \, *} r_1)}{|r_1|^2} \\
    &\ge \inf_{r_1 \in \mathbb{C}^m \backslash \{0\}} 
    \frac{(\lambda - \kappa_2^{\pm})^{-1} |V_{12}^{\pm \, *} r_1|^2}{|r_1|^2}
    = (\lambda - \kappa_2^{\pm})^{-1} 
    \inf_{r_1 \in \mathbb{C}^m \backslash \{0\}} 
    \frac{(V_{12}^{\pm} V_{12}^{\pm \, *} r_1, r_1) }{|r_1|^2} \\
    & = (\lambda - \kappa_2^{\pm})^{-1} \rho^{\pm}, 
    \end{aligned}
\end{equation*}
where $\rho^{\pm}$ denote the lowest eigenvalues of 
$V_{12}^{\pm} V_{12}^{\pm \, *}$. In summary, we can 
ensure that $\lambda \in \mathcal{R}_0$ by taking 
$\lambda$ to satisfy
\begin{equation*}
    \lambda < \kappa_1^{\pm} + \frac{\rho^{\pm}}{(\lambda - \kappa_2^{\pm})}
\end{equation*}
(i.e., satisfy two inequalities). 
Since $\lambda - \kappa_2^{\pm} < 0$, we can express this relation 
as the quadratic inequality 
\begin{equation*}
    \lambda^2 - (\kappa_1^{\pm} + \kappa_2^{\pm}) \lambda 
    + \kappa_1^{\pm} \kappa_2^{\pm} - \rho^{\pm, \, 2} > 0.
\end{equation*}
(We emphasize that we are taking $\lambda$ below $\cup_{k=1}^{n-m} \mathcal{R}_k$,
so this does not assert that large positive values of $\lambda$ are 
admissible.) Upon solving this quadratic inequality, we find 
that admissible values of $\lambda$ include those values 
below $\cup_{k=1}^{n-m} \mathcal{R}_k$ that also satisfy 
the inequality 
\begin{equation*}
    \lambda < \min_{\pm} \Big{\{}\frac{1}{2} \Big((\kappa_1^{\pm} + \kappa_2^{\pm})
    - \sqrt{(\kappa_1^{\pm} - \kappa_2^{\pm})^2 + 4 \rho^{\pm \, 2}}\Big) \Big{\}}.
\end{equation*}

Returning to the general case, we 
let $I$ denote any interval satisfying (\ref{sls-a-i}), not necessarily 
below $\sigma_{\ess}(\mathcal{L}_a)$. 
For $\lambda \in I$, we are now in a position to 
develop frames $\mathbf{X} (x; \lambda)$ and 
$\tilde{\mathbf{X}} (x; \lambda)$ as described in 
{\bf (B1)}. For this, we begin by looking for solutions
to (\ref{sls-a-limit}) of the form 
$\phi (x; \lambda) = e^{\mu (\lambda) x} r(\lambda)$,
where $\mu: I \to \mathbb{R}$ and $r_1: I \to \mathbb{C}^n$.
We find, 
\begin{equation*}
    \{-\mu^2 P_{11}^{\pm} + \mathbf{V}^{\pm} (\lambda) - \lambda I \}r_1 
    = 0.
\end{equation*}
The allowable values of $\mu^2$ are precisely the eigenvalues 
of $(P_{11}^{\pm})^{-1} (\mathbf{V}^{\pm} (\lambda) - \lambda I)$, which is self-adjoint
with respect to the inner product 
$(r, s)_{P_{11}^{\pm}} := (P_{11}^{\pm} r, s)$. We conclude 
that these eigenvalues will be real-valued, and that we can 
choose the associated eigenvectors to be orthonormal 
with respect to this inner product. In addition, for 
$\lambda \in I$, we have 
\begin{equation*}
    \lambda < \inf_{r_1 \in \mathbb{C}^m \backslash \{0\}} 
    \frac{(\mathbf{V}^{\pm} (\lambda) r_1, r_1)}{|r_1|^2},
\end{equation*}
so that $\mathbf{V}^{\pm} (\lambda) - \lambda I$ is a positive
matrix. We conclude that $\mu^2$ takes only positive real 
values, leading to $n$ negative values for $\mu$ and $n$
positive values. We will denote these values 
$\{\mu_k^{\pm} (\lambda)\}_{k=1}^{2n}$, with the first 
$n$ values negative, the second $n$ values positive, 
and the relation $\mu_{n+k}^{\pm} (\lambda) = - \mu_k^{\pm} (\lambda)$
for all $k \in \{1, 2, \dots, n\}$. We denote the 
corresponding eigenvectors $\{r_k^{\pm} (\lambda)\}_{k=1}^n$
so that 
\begin{equation*}
(P_{11}^{\pm})^{-1} (\mathbf{V}^{\pm} (\lambda) - \lambda I) r_k^{\pm}    
= (\mu_k^{\pm})^2 r_k^{\pm}, \quad \forall \, k \in \{1, 2, \dots, n\}. 
\end{equation*}

In order to place (\ref{phi1-equation}) in our general framework, we
will set $y = {y_1 \choose y_2} = {\phi_1 \choose P_{11} (x) \phi_1'}$
so that we have 
\begin{equation} \label{sls-a-A-ode}
y' = \mathbb{A} (x; \lambda) y;
\quad \mathbb{A} (x; \lambda)
= \begin{pmatrix}
0 & P_{11} (x)^{-1} \\
\mathbf{V} (x; \lambda) - \lambda I & 0 
\end{pmatrix},
\end{equation}
or equivalently 
\begin{equation} \label{sls-a-hammy}
J y' = \mathbb{B} (x; \lambda) y;
\quad \mathbb{B} (x; \lambda)
= \begin{pmatrix}
\lambda I -  \mathbf{V} (x; \lambda) & 0 \\
0 & P_{11} (x)^{-1} 
\end{pmatrix}. 
\end{equation}

For $\lambda \in I$, we see that $\mathbb{B} (\cdot; \lambda) \in L^1_{\loc} (\mathbb{R}, \mathbb{C}^{2m \times 2m})$,
and it's clear that $\mathbb{B} (x; \lambda)$ is self-adjoint for all $x \in \mathbb{R}$. 
We also need to compute $\mathbb{B}_{\lambda} (x; \lambda)$, and for this, we first observe
that 
\begin{equation} \label{V-derivative}
    \begin{aligned}
        \mathbf{V}_{\lambda} (x; \lambda) 
        &= - V_{12} (x) (\lambda I - V_{22} (x))^{-2} V_{12} (x)^* \\
        &= - \Big( (\lambda I - V_{22} (x))^{-1} V_{12} (x)^* \Big)^* (\lambda I - V_{22} (x))^{-1} V_{12} (x)^*.
    \end{aligned}
\end{equation}
Recalling that $V \in C(\mathbb{R}, \mathbb{C}^{n \times n})$, 
and that for $\lambda \in I$, we have $\lambda \notin \sigma (V_{22} (x)) \cup \sigma (V_{22}^{\pm})$
for all $x \in \mathbb{R}$, we see that 
$\mathbf{V}_{\lambda} (\cdot; \lambda) \in L^1_{\loc} (\mathbb{R}, \mathbb{C}^{m \times m})$,
and consequently $\mathbb{B}_{\lambda} (\cdot; \lambda) \in L^1_{\loc} (\mathbb{R}, \mathbb{C}^{2m \times 2m})$.
This establishes Assumptions {\bf (A)}. 

For {\bf (B1)}, we will proceed precisely as we did with 
Sturm-Liouville Systems. Similarly as in (\ref{sls-bold-r}),
the values $\{\mu_k^{\pm} (\lambda)\}_{k=1}^{2n}$ described above comprise a labeling
of the eigenvalues of $\mathbb{A}_{\pm} (\lambda) := \lim_{x \to \pm \infty} \mathbb{A} (x; \lambda)$. 
If we let $\{\mathbf{r}_k^{\pm} (\lambda)\}_{k=1}^{2n}$
denote the eigenvectors of $\mathbb{A}_{\pm} (\lambda)$ respectively associated with 
these eigenvectors, then we find 
\begin{equation} \label{sls-a-bold-r}
    \mathbf{r}_k^{\pm} (\lambda) = {r_k^{\pm} (\lambda) \choose \mu_k^{\pm} (\lambda) P_{11}^{\pm} r_k^{\pm} (\lambda)}; 
    \quad \mathbf{r}_{n+k}^{\pm} (\lambda) = {r_k^{\pm} (\lambda) \choose - \mu_k^{\pm} (\lambda) P_{11}^{\pm} r_k^{\pm} (\lambda)};
    \quad k = 1, 2, \dots, n.
\end{equation}
The following lemma can be established by a proof almost
identical to the proof of Lemma \ref{sls-ode-lemma}.

\begin{lemma} \label{sls-a-ode-lemma} 
Assume {\bf (DA1)} and {\bf (DA2)} hold, and let $I$ be 
as in (\ref{sls-a-i}). Also, let $\{\mu_k^{\pm} (\lambda)\}_{k=1}^{2m}$
and $\{\mathbf{r}_k^{\pm} (\lambda)\}_{k=1}^{2m}$ be as described just above.
Then there exists a family of bases
$\{\mathbf{y}_{k}^- (\cdot; \lambda)\}_{k=m+1}^{2m}$,
$\lambda \in I$,
for the spaces of solutions to (\ref{sls-a-A-ode}) that lie
left in $\mathbb{R}$, and a family of bases
$\{\mathbf{y}_{k}^+ (\cdot; \lambda)\}_{k=1}^{m}$, 
$\lambda \in I$,
for the spaces of solutions to (\ref{sls-a-A-ode}) that lie
right in $\mathbb{R}$. Respectively, we can choose these
so that 
\begin{equation*}
\begin{aligned}
  \mathbf{y}_{m+k}^- (x; \lambda)
  &= e^{- \mu_k^- (\lambda) x} (\mathbf{r}_{m+k}^- (\lambda) + \mathbf{E}_{m+k}^- (x; \lambda)),
  \quad k = 1, 2, \dots, m, \\
   \mathbf{y}_{k}^+ (x; \lambda)
  &= e^{\mu_k^+ (\lambda) x} (\mathbf{r}_{k}^+ (\lambda) + \mathbf{E}_{k}^+ (x; \lambda)),
  \quad k = 1, 2, \dots, m,
\end{aligned}
\end{equation*}
where for any fixed interval $[\lambda_1, \lambda_2] \subset I$,
there exists a constant $\delta > 0$ so that for each 
$k \in \{1, 2, \dots, m\}$
\begin{equation*}
    \mathbf{E}_{m+k}^- (x; \lambda) = \mathbf{O} (e^{- \delta |x|}), 
    \quad x \to - \infty; 
    \quad \quad  \mathbf{E}_{k}^+ (x; \lambda) = \mathbf{O} (e^{- \delta |x|}), 
    \quad x \to + \infty,
\end{equation*}
uniformly for $\lambda \in [\lambda_1, \lambda_2]$. 

Moreover, there exists a $\lambda$-dependent family of bases
$\{\mathbf{y}_{k}^- (\cdot; \lambda)\}_{k=1}^{m}$, 
$\lambda \in I$,
for the spaces of solutions to (\ref{sls-a-A-ode}) that do not 
lie left in $\mathbb{R}$, and a $\lambda$-dependent family of bases
$\{\mathbf{y}_{k}^+ (\cdot; \lambda)\}_{k=m+1}^{2m}$, 
$\lambda \in I$,
for the spaces of solutions to (\ref{sls-a-A-ode}) that 
do not lie right in $\mathbb{R}$. Respectively, we can choose these
so that 
\begin{equation*}
\begin{aligned}
  \mathbf{y}_{k}^- (x; \lambda)
  &= e^{\mu_k^- (\lambda) x} (\mathbf{r}_{k}^- (\lambda) + \mathbf{E}_{k}^- (x; \lambda)),
  \quad k = 1, 2, \dots, m, \\
   \mathbf{y}_{m+k}^+ (x; \lambda)
  &= e^{- \mu_k^+ (\lambda) x} (\mathbf{r}_{m+k}^+ (\lambda) + \mathbf{E}_{m+k}^+ (x; \lambda)),
  \quad k = 1, 2, \dots, m,
\end{aligned}
\end{equation*}
where for any fixed interval $[\lambda_1, \lambda_2] \subset I$,
there exist a constant $\delta > 0$ so that for each 
$k \in \{1, 2, \dots, m\}$
\begin{equation*}
    \mathbf{E}_{k}^- (x; \lambda) = \mathbf{O} (e^{- \delta |x|}), 
    \quad x \to - \infty; 
    \quad \quad  \mathbf{E}_{m+k}^+ (x; \lambda) = \mathbf{O} (e^{- \delta |x|}), 
    \quad x \to + \infty,
\end{equation*}
uniformly for $\lambda \in [\lambda_1, \lambda_2]$. 
\end{lemma}

Precisely as in the case of Sturm-Liouville Systems, we require
the following two auxiliary lemmas, which are again taken 
from \cite{HS2020} (with a straightforward modification in 
this case, extending the result from cases in which 
$\mathbb{B} (x; \lambda)$ is linear in $\lambda$ to 
cases in which it is analytic in $\lambda$). 

\begin{lemma} \label{sls-a-extension-lemma1}
Assume {\bf (DA1)} and {\bf (DA2)} hold, and for 
each $\lambda \in I$ (with $I$ as in (\ref{sls-a-i})) 
let $\{\mathbf{y}_{k}^- (\cdot; \lambda)\}_{k=m+1}^{2m}$
and $\{\mathbf{y}_{k}^+ (\cdot; \lambda)\}_{k=1}^{m}$
be as described in Lemma \ref{sls-a-ode-lemma}. If 
$\ell (x; \lambda)$ and $\tilde{\ell} (x; \lambda)$
respectively denote the Lagrangian subspaces with 
frames 
\begin{equation} \label{sls-a-lframe}
    \mathbf{X} (x; \lambda)
    = (\mathbf{y}_{m+1}^- (x; \lambda) \,\,\, \mathbf{y}_{m+2}^- (x; \lambda) \,\,\, \cdots \,\,\, \mathbf{y}_{2m}^- (x; \lambda)),
\end{equation}
and 
\begin{equation} \label{sls-a-rframe}
    \tilde{\mathbf{X}} (x; \lambda)
    = (\mathbf{y}_{1}^+ (x; \lambda) \,\,\, \mathbf{y}_{2}^+ (x; \lambda) \,\,\, \cdots \,\,\, \mathbf{y}_{m}^+ (x; \lambda)),
\end{equation}
then $\ell, \tilde{\ell} \in C (\mathbb{R} \times I, \Lambda (n))$. 
\end{lemma}

\begin{lemma} \label{sls-a-extension-lemma2}
Assume {\bf (DA1)} and {\bf (DA2)} hold, and for some 
fixed $\lambda_0 \in I$ (with $I$ as in (\ref{sls-a-i})) let
$\{\mathbf{y}_{k}^- (\cdot; \lambda_0)\}_{k=m+1}^{2m}$
and $\{\mathbf{y}_{k}^+ (\cdot; \lambda_0)\}_{k=1}^{m}$
be as described in Lemma \ref{sls-a-ode-lemma}. Then 
there exists a constant $r_0 > 0$ so that the elements
$\{\mathbf{y}_{k}^- (\cdot; \lambda_0)\}_{k=m+1}^{2m}$
and $\{\mathbf{y}_{k}^+ (\cdot; \lambda_0)\}_{k=1}^{m}$
can be analytically extended in $\lambda$ to the complex 
ball $B (\lambda_0, r_0)$ (centered at $\lambda_0$ with radius
$r_0$). Moreover, the analytic extensions of 
$\{\mathbf{y}_{k}^- (\cdot; \lambda_0)\}_{k=m+1}^{2m}$
comprise a basis for the space of solutions of (\ref{sls-a-hammy})
that lie left in $\mathbb{R}$, and the analytic extensions 
of $\{\mathbf{y}_{k}^+ (\cdot; \lambda_0)\}_{k=1}^{m}$
comprise a basis for the space of solutions of (\ref{sls-a-hammy})
that lie right in $\mathbb{R}$. 
\end{lemma}

Proceeding similarly as with Sturm-Liouville Systems, we can 
respectively replace the frames $\mathbf{X} (x; \lambda)$
and $\tilde{\mathbf{X}} (x; \lambda)$ specified in 
(\ref{sls-a-lframe}) and (\ref{sls-a-rframe}) 
with $\mathbf{X} (x; \lambda) e^{D_- (\lambda) x}$
and $\tilde{\mathbf{X}} (x; \lambda) e^{- D_+ (\lambda) x}$, where 
\begin{equation*}
    D_{\pm} (\lambda)
    = \diag (\mu_1^{\pm} (\lambda) \,\,\, \mu_2^{\pm} (\lambda) \,\,\, \dots \,\,\, \mu_m^{\pm} (\lambda)).
\end{equation*}
It follows that the frames for $\ell_- (\lambda)$ and $\tilde{\ell}_+ (\lambda)$ 
can be taken respectively to be  
\begin{equation*}
    \mathbf{X}_- (\lambda)
    = {R_- (\lambda) \choose - P_{11}^- R_- (\lambda) D_- (\lambda)};
    \quad
    \tilde{\mathbf{X}}_+ (\lambda)
    = {R_+ (\lambda) \choose P_{11}^+ R_+ (\lambda) D_+ (\lambda)},
\end{equation*}
where 
\begin{equation*}
    R_{\pm} (\lambda) = (r_1^{\pm} (\lambda) \,\,\, r_2^{\pm} (\lambda) \,\,\, \dots 
    \,\,\, r_m^{\pm} (\lambda)).
\end{equation*}
This establishes {\bf (B1)} and the first part of {\bf (B2)}, and the second 
part of {\bf (B2)} can be established precisely as for Sturm-Liouville systems.  

For {\bf (B3)}, we have 
\begin{equation*}
    \mathbb{B}_{\lambda} (x; \lambda)
    = \begin{pmatrix}
    I - \mathbf{V}_{\lambda} (x; \lambda) & 0 \\
    0 & 0
    \end{pmatrix},
\end{equation*}
from which we see that 
\begin{equation*}
    \mathbf{X} (x; \lambda)^* \mathbb{B}_{\lambda} (x; \lambda) \mathbf{X} (x; \lambda)
    = X (x; \lambda)^* (I - \mathbf{V}_{\lambda} (x; \lambda)) X (x; \lambda).
\end{equation*}
It's clear from (\ref{V-derivative}) that $- \mathbf{V}_{\lambda} (x; \lambda)$ is
non-negative, and so $I - \mathbf{V}_{\lambda} (x; \lambda)$ is positive 
definite. From this observation, {\bf (B3)} follows immediately as in 
Section \ref{sls-section}. We conclude that the assumptions of Theorem \ref{main-theorem}
hold in this case, and this gives the first part of Theorem \ref{sls-a-theorem}. 

For the remainder of Theorem \ref{sls-a-theorem}, it follows from the 
structure of $\mathbf{X} (x; \lambda)$, $\tilde{\mathbf{X}} (x; \lambda)$,
and $\mathbb{B} (x; \lambda)$ that the relevant calculations from
Section \ref{sls-section} can be used to show that the target spaces 
$\tilde{\ell}_+ (\lambda_1)$ and $\tilde{\ell}_+ (\lambda_2)$ can 
be replaced with the Dirichlet space $\ell_D$, and also that the
Maslov index with $\ell_D$ as the target space has monotonic
conjugate points. 

Last, suppose $[\lambda_1, \lambda_2] \subset I$ lies entirely below
the essential spectrum of $\mathcal{L}_a$. Then $\lambda_1$ can be 
chosen sufficiently negative so that there are no conjugate points 
along the vertical shelf at $\lambda_1$. To see this, we again proceed 
as in Section \ref{sls-section}, observing that a point
$(s, \lambda_1) \in \mathbb{R} \times I$ will be conjugate if 
and only if $\lambda_1$ is an eigenvalue of the half-line problem 
\begin{equation*}
    \begin{aligned}
      - (P_{11} (x) \phi_1')' + \mathbf{V} (x; \lambda) \phi_1 &= \lambda \phi_1 \\
      \tilde{\mathbf{X}}_+ (\lambda_1)^* J {\phi_1 (s) \choose P_{11} (s) \phi_1'(s)}
      &= 0.
    \end{aligned}
\end{equation*}
Proceeding as in Section \ref{sls-section}, the only new aspect is the term 
\begin{equation*}
    \int_{-\infty}^s (\mathbf{V} (x; \lambda) \phi_1 (x; \lambda), \phi_1 (x; \lambda)) dx,
\end{equation*}
which we bound (in absolute value) by 
$\|\mathbf{V} (\cdot; \lambda)\|_{L^{\infty} (-\infty, s)} \|\phi_1 (\cdot; \lambda)\|_{L^{2} (-\infty, s)}^2$.
In the current setting, 
\begin{equation*}
|\mathbf{V} (x; \lambda)| \le |V_{11} (x)| + |V_{12} (x)| |(\lambda I - V_{22} (x))^{-1}| |V_{12} (x)^*|,    
\end{equation*}
where $|\cdot|$ denotes any matrix norm. Using the facts that $\lambda \in I$, $V \in C (\mathbb{R}, \mathbb{C}^{n \times n})$,
and the limit conditions {\bf (DA2)}, we conclude that $|\mathbf{V} (x; \lambda)|$ is bounded
independently of $x$ and $\lambda$ (for $\lambda < \lambda_2$). We conclude that if we 
take $\lambda_1$ sufficiently negative, 
there will be no conjugate points along the vertical shelf at $\lambda_1$. 
This completes the proof of Theorem \ref{sls-a-theorem}. 
\hfill $\square$

\section{Fourth Order Potential Systems} \label{fourth-section}

In this section, we apply Theorem \ref{main-theorem} to fourth-order
potential systems
\begin{equation} \label{fourth}
    \mathcal{L} \phi := - \phi'''' + V(x) \phi = \lambda \phi;
    \quad x \in \mathbb{R}, 
    \quad \phi (x; \lambda) \in \mathbb{C}^n.
\end{equation}
In order to ensure that our general assumptions {\bf (A)},
{\bf (B1)}, {\bf (B2)}, and {\bf (B3)} hold, we make the following
assumptions on the coefficient matrix $V$. 

\medskip
{\bf (FP1)} We take $V \in C (\mathbb{R}, \mathbb{C}^{n \times n})$, 
with $V(x)$ self-adjoint for all $x \in \mathbb{R}$.

\medskip
{\bf (FP2)} We assume the limits $\lim_{x \to \pm \infty} V(x) = V_{a}$
exist and agree, and 
\begin{equation*}
    \int_{-\infty}^{+\infty} (1+|x|) (V(x) - V_a) dx < \infty.
\end{equation*}

\begin{remark}
We emphasize that in this case we take the endstates $\lim_{x \to \pm \infty} V(x) = V_{\pm}$
to agree. This corresponds with cases in which the PDE 
\begin{equation*}
    u_t + DF(u) = - u_{xxxx}; 
    \quad (x, t) \in \mathbb{R} \times \mathbb{R}_+,
    \quad u(x, t) \in \mathbb{C}^n, 
\end{equation*}
is linearized about a stationary solution $\bar{u} (x)$ for which 
the endstates $u_{\pm}$ agree. If $V_- \ne V_+$, the analysis 
becomes substantially more technical, and we leave such cases to
future studies. 
\end{remark}

We take as our domain for $\mathcal{L}$ the set $H^4 (\mathbb{R}, \mathbb{C}^n)$,
noting from \cite{We1987} that with this choice $\mathcal{L}$ is self-adjoint. 
As with Sturm-Liouville systems, the essential spectrum of $\mathcal{L}$
is determined by the asymptotic problem
\begin{equation} \label{fourth-asymptotic}
\phi'''' + V_{a} \phi = \lambda \phi.
\end{equation}
Precisely, if we look for solutions of the form 
$\phi (x) = e^{ikx} r$, then the essential spectrum 
of (\ref{fourth-asymptotic}) is precisely the collection
of $\lambda \in \mathbb{R}$ for which $\phi (x) = e^{i k x}$ solves
(\ref{fourth-asymptotic}) for some $k\in \mathbb{R}$
and $r \in \mathbb{C}^n$. Upon substitution of 
$\phi (x) = e^{ikx} r$ in to (\ref{fourth-asymptotic}),
we obtain 
\begin{equation*}
    (k^4 I + V_{a})r = \lambda r
    \quad \implies \quad
    k^4 |r|^2 + (V_{a} r, r) = \lambda |r|^2.
\end{equation*}
We see from this that if we set 
\begin{equation} \label{fourth-kappa-defined}
    \kappa := \inf_{r \ne 0} \frac{(V_{a} r, r)}{|r|^2},
\end{equation}
then $\sigma_{\ess} (\mathcal{L}) = [\kappa, \infty)$. This will allow us 
to take the interval $I$ in Assumptions {\bf (A)}, {\bf (B1)},
{\bf (B2)}, and {\bf (B3)} to be 
$I = (-\infty, \kappa)$. 

In order to characterize the Lagrangian subspaces 
$\ell (x; \lambda)$ and $\tilde{\ell} (x; \lambda)$ 
described in Assumption {\bf (B1)}, we will need
a lemma analogous to Lemma \ref{sls-ode-lemma}. 
In order to develop such a lemma, we begin by looking
for solutions of (\ref{fourth-asymptotic}) of the 
form $\phi (x; \lambda) = e^{\mu (\lambda) x} r$,
where in this case $\mu$ is a real-valued function 
of $\lambda$ and $r$ is a constant vector 
$r \in \mathbb{C}^n$. We see that 
\begin{equation*}
    (\mu^4 I + V_{a} - \lambda I) r = 0,
\end{equation*}
so in particular, the allowable values of 
$\lambda - \mu^4$ are eigenvalues of the 
matrix $V_{a}$. I.e., if we denote the 
eigenvalues of $V_{a}$ by $\{\nu_k\}_{k=1}^n$,
then each allowable value of $\mu^4$ must
satisfy
\begin{equation*}
\lambda - \mu^4 = \nu_k   
\end{equation*}
for some $\nu_k \in \sigma (V_{a})$. 
Each such $\nu_k$ will correspond with 
four values of $\mu$, and we will denote 
the full collection of such values $\{\mu_k\}_{k=1}^{4n}$,
indexed so that for each $k \in \{1, 2 \dots, n\}$,
\begin{equation*}
    \begin{aligned}
      \mu_k (\lambda) 
      &= (-\frac{1}{\sqrt{2}} - i \frac{1}{\sqrt{2}}) \sqrt[4]{\nu_k - \lambda};
      \quad  \mu_{n+k} (\lambda) 
      = (-\frac{1}{\sqrt{2}} + i \frac{1}{\sqrt{2}}) \sqrt[4]{\nu_k - \lambda}; \\
      \mu_{2n + k} (\lambda) 
      &= (\frac{1}{\sqrt{2}} + i \frac{1}{\sqrt{2}}) \sqrt[4]{\nu_k - \lambda};
      \quad  \mu_{3n+k} (\lambda) 
      = (\frac{1}{\sqrt{2}} - i \frac{1}{\sqrt{2}}) \sqrt[4]{\nu_k - \lambda}.
    \end{aligned}
\end{equation*}
We note that with this choice of indexing, we have the relations
\begin{equation*}
    \mu_{2n + k} (\lambda) = - \mu_k (\lambda); 
    \quad \mu_{3n+k} (\lambda) = - \mu_{n+k} (\lambda); 
    \quad \forall \, k \in \{1, 2, \dots, n\}.
\end{equation*}

For each $k \in \{1, 2 \dots, n\}$, the values $\mu_k (\lambda)$,
$\mu_{n+k} (\lambda)$, $\mu_{2n + k} (\lambda)$, and 
$\mu_{3n+k} (\lambda)$ all correspond with the same eigenvector
of $V_{a}$, which we denote $r_k$. For the set 
$\{\mu_k (\lambda)\}_{k=1}^n$, we can express this as 
\begin{equation*}
    (\mu_k (\lambda)^4 I + V_{a} - \lambda I) r_k = 0.
\end{equation*}
Since the matrix $V_{a}$ is self-adjoint, we can choose the 
collection $\{r_k\}_{k=1}^n$ to be orthonormal.
We will set 
\begin{equation*}
    R = (r_1 \,\,\, r_2 \,\,\, \cdots \,\,\, r_n),
\end{equation*}
for which orthonormality can be expressed as $R^* R = I$. 

In order to place (\ref{fourth}) in our general framework, we will 
express it as a first-order system. For this, it will be convenient 
to make the choices $y_1 = \phi$, $y_2 = \phi''$, $y_3 = - \phi'''$,
and $y_4 = - \phi'$, for which we find 
\begin{equation} \label{fourth-ode-A}
    y' = \mathbb{A} (x; \lambda) y; \quad 
    \mathbb{A} (x; \lambda) 
    = \begin{pmatrix}
    0 & 0 & 0 & -I \\
    0 & 0 & -I & 0 \\
    V(x) - \lambda I & 0 & 0 & 0 \\
    0 & -I & 0 & 0 
    \end{pmatrix}, 
\end{equation}
or equivalently 
\begin{equation} \label{fourth-hammy}
    J y' = \mathbb{B} (x; \lambda) y; \quad 
    \mathbb{B} (x; \lambda) 
    = \begin{pmatrix}
    \lambda I - V(x) & 0 & 0 & 0 \\
    0 & I & 0 & 0 \\
    0 & 0 & 0 & -I \\
    0 & 0 & -I & 0 
    \end{pmatrix}.
\end{equation}
(We refer the reader to \cite{HJK2018} for a full discussion of 
the motivation behind these choices for the vector $y$.)
The values $\{\mu_k\}_{k=1}^{4n}$ are precisely the 
eigenvalues of the matrix 
\begin{equation*}
    \mathbb{A}_{a} (\lambda) := \lim_{x \to \pm \infty} \mathbb{A} (x; \lambda),
\end{equation*}
and it's straightforward to check that the associated eigenvectors
are respectively
\begin{equation*}
    \quad \mathbf{r}_{pn+k} (\lambda)
    = \begin{pmatrix}
    r_k \\ (\mu_{pn+k})^2 r_k \\ - (\mu_{pn+k})^3 r_j \\ - \mu_{pn+k} r_k
    \end{pmatrix};
    \quad p = 0, 1, 2, 3.
\end{equation*}

The following lemma can be proven in almost precisely the same way 
as Lemma 2.2 in \cite{HLS2018}. 

\begin{lemma} \label{fourth-ode-lemma} 
Assume {\bf (FP1)} and {\bf (FP2)} hold, and let $\{\mu_k (\lambda)\}_{k=1}^{4n}$
and $\{\mathbf{r}_k (\lambda)\}_{k=1}^{4n}$ be as described just above.
Then there exists a family of bases
$\{\mathbf{y}_{k}^- (\cdot; \lambda)\}_{k=2n+1}^{2n}$,
$\lambda \in (-\infty, \kappa)$,
for the spaces of solutions to (\ref{fourth-ode-A}) that lie
left in $\mathbb{R}$, and a family of bases
$\{\mathbf{y}_{k}^+ (\cdot; \lambda)\}_{k=1}^{2n}$, 
$\lambda \in (-\infty, \kappa)$,
for the spaces of solutions to (\ref{fourth-ode-A}) that lie
right in $\mathbb{R}$. 
Respectively, we can choose these so that 
\begin{equation*}
\begin{aligned}
  \mathbf{y}_{2n+k}^- (x; \lambda)
  &= e^{- \mu_k (\lambda) x} (\mathbf{r}_{2n+k} (\lambda) + \mathbf{E}_{2n+k}^- (x; \lambda)),
  \quad k = 1, 2, \dots, 2n, \\
   \mathbf{y}_{k}^+ (x; \lambda)
  &= e^{\mu_k (\lambda) x} (\mathbf{r}_{k} (\lambda) + \mathbf{E}_{k}^+ (x; \lambda)),
  \quad k = 1, 2, \dots, 2n,
\end{aligned}
\end{equation*}
where for any fixed interval $[\lambda_1, \lambda_2] \subset (-\infty, \kappa)$,
$\lambda_1 < \lambda_2$, and for any 
$k \in \{1, 2, \dots, 2n\}$
\begin{equation*}
    \mathbf{E}_{2n+k}^- (x; \lambda) = \mathbf{O} ((1+|x|)^{-1}), 
    \quad x \to - \infty; 
    \quad \quad  \mathbf{E}_{k}^+ (x; \lambda) = \mathbf{O} ((1+|x|)^{-1}), 
    \quad x \to + \infty,
\end{equation*}
uniformly for $\lambda \in [\lambda_1, \lambda_2]$. 

Moreover, there exists a $\lambda$-dependent family of bases
$\{\mathbf{y}_{k}^- (\cdot; \lambda)\}_{k=1}^{2n}$, 
$\lambda \in (-\infty, \kappa)$,
for the spaces of solutions to (\ref{fourth-ode-A}) that do not 
lie left in $\mathbb{R}$, and a $\lambda$-dependent family of bases
$\{\mathbf{y}_{k}^+ (\cdot; \lambda)\}_{k=2n+1}^{4n}$, 
$\lambda \in (-\infty, k)$,
for the spaces of solutions to (\ref{fourth-ode-A}) that 
do not lie right in $\mathbb{R}$. Respectively, we can choose these
so that 
\begin{equation*}
\begin{aligned}
  \mathbf{y}_{k}^- (x; \lambda)
  &= e^{\mu_k (\lambda) x} (\mathbf{r}_{k} (\lambda) + \mathbf{E}_{k}^- (x; \lambda)),
  \quad k = 1, 2, \dots, 2n, \\
   \mathbf{y}_{2n+k}^+ (x; \lambda)
  &= e^{- \mu_k (\lambda) x} (\mathbf{r}_{2n+k} (\lambda) + \mathbf{E}_{2n+k}^+ (x; \lambda)),
  \quad k = 1, 2, \dots, 2n,
\end{aligned}
\end{equation*}
where for any fixed interval $[\lambda_1, \lambda_2] \subset (-\infty, \kappa)$, 
$\lambda_1 < \lambda_2$, and for any 
$k \in \{1, 2, \dots, 2n\}$
\begin{equation*}
    \mathbf{E}_{k}^- (x; \lambda) = \mathbf{O} ((1+|x|)^{-1}), 
    \quad x \to - \infty; 
    \quad \quad  \mathbf{E}_{2n+k}^+ (x; \lambda) = \mathbf{O} ((1+|x|)^{-1}), 
    \quad x \to + \infty,
\end{equation*}
uniformly for $\lambda \in [\lambda_1, \lambda_2]$. 
\end{lemma}

Precisely as in the previous cases, we require
the following two auxiliary lemmas, which are again taken 
from \cite{HS2020}. 

\begin{lemma} \label{fourth-extension-lemma1}
Assume {\bf (FP1)} and {\bf (FP2)} hold, and for 
each $\lambda \in (-\infty, \kappa)$ 
let $\{\mathbf{y}_{k}^- (\cdot; \lambda)\}_{k=2n+1}^{4n}$
and $\{\mathbf{y}_{k}^+ (\cdot; \lambda)\}_{k=1}^{2n}$
be as described in Lemma \ref{fourth-ode-lemma}. If 
$\ell (x; \lambda)$ and $\tilde{\ell} (x; \lambda)$
respectively denote the Lagrangian subspaces with 
frames 
\begin{equation} \label{fourth-lframe}
    \mathbf{X} (x; \lambda)
    = (\mathbf{y}_{2n+1}^- (x; \lambda) \,\,\, \mathbf{y}_{2n+2}^- (x; \lambda) \,\,\, \cdots \,\,\, \mathbf{y}_{4n}^- (x; \lambda)),
\end{equation}
and 
\begin{equation} \label{fourth-rframe}
    \tilde{\mathbf{X}} (x; \lambda)
    = (\mathbf{y}_{1}^+ (x; \lambda) \,\,\, \mathbf{y}_{2}^+ (x; \lambda) \,\,\, \cdots \,\,\, \mathbf{y}_{2n}^+ (x; \lambda)),
\end{equation}
then $\ell, \tilde{\ell} \in C (\mathbb{R} \times (-\infty, \kappa), \Lambda (n))$. 
\end{lemma}

\begin{lemma} \label{fourth-extension-lemma2}
Assume {\bf (FP1)} and {\bf (FP2)} hold, and for some 
fixed $\lambda_0 \in (-\infty, \kappa)$ let
$\{\mathbf{y}_{k}^- (\cdot; \lambda_0)\}_{k=2n+1}^{4n}$
and $\{\mathbf{y}_{k}^+ (\cdot; \lambda_0)\}_{k=1}^{2n}$
be as described in Lemma \ref{fourth-ode-lemma}. Then 
there exists a constant $r_0 > 0$ so that the elements
$\{\mathbf{y}_{k}^- (\cdot; \lambda_0)\}_{k=2n+1}^{4n}$
and $\{\mathbf{y}_{k}^+ (\cdot; \lambda_0)\}_{k=1}^{2n}$
can be analytically extended in $\lambda$ to the complex 
ball $B (\lambda_0, r_0)$ (centered at $\lambda_0$ with radius
$r_0$). Moreover, the analytic extensions of 
$\{\mathbf{y}_{k}^- (\cdot; \lambda_0)\}_{k=2n+1}^{4n}$
comprise a basis for the space of solutions of (\ref{fourth-hammy})
that lie left in $\mathbb{R}$, and the analytic extensions 
of $\{\mathbf{y}_{k}^+ (\cdot; \lambda_0)\}_{k=1}^{2n}$
comprise a basis for the space of solutions of (\ref{fourth-hammy})
that lie right in $\mathbb{R}$. 
\end{lemma}

We will set 
\begin{equation*}
    \mathcal{D} (\lambda) 
    = \diag (\mu_1 (\lambda) \,\,\, \mu_2 (\lambda) \,\,\, \dots \,\,\, \mu_{2n} (\lambda)),
\end{equation*}
and we note that our labeling conventions have been chosen so that 
\begin{equation*}
    - \mathcal{D} (\lambda) 
    = \diag (\mu_{2n+1} (\lambda) \,\,\, \mu_{2n+2} (\lambda) \,\,\, \dots \,\,\, \mu_{4n} (\lambda)).
\end{equation*}
If we replace $\mathbf{X} (x; \lambda)$ with $\mathbf{X} (x; \lambda) e^{\mathcal{D} (\lambda) x}$
and $\tilde{\mathbf{X}} (x; \lambda)$ with $\tilde{\mathbf{X}} (x; \lambda) e^{- \mathcal{D} (\lambda) x}$,
we readily see that the asymptotic Lagrangian subspaces $\ell_- (\lambda)$ and $\tilde{\ell}_+ (\lambda)$
are well-defined with respective frames 
\begin{equation} \label{fourth-asym-frames}
    \mathbf{X}_- (\lambda) 
    = \begin{pmatrix}
    R & R \\ R D (\lambda)^2 & R (D (\lambda)^*)^2\\
    R D (\lambda)^3 & R (D (\lambda)^*)^3 \\ R D (\lambda) & R D (\lambda)^*
    \end{pmatrix};
    \quad \tilde{\mathbf{X}}_+ (\lambda) 
    = \begin{pmatrix}
    R & R \\ R D (\lambda)^2 & R (D (\lambda)^*)^2\\
    - R D (\lambda)^3 & - R (D (\lambda)^*)^3 \\ - R D (\lambda) & - R D (\lambda)^*
    \end{pmatrix}.
\end{equation}
Likewise, we obtain asymptotic frames associated with solutions that do not lie left (respectively right)
in $\mathbb{R}$, and we see from Lemma \ref{fourth-ode-lemma} that these will be 
$\mathbf{X}_-^g (\lambda) = \tilde{\mathbf{X}}_+ (\lambda)$ and
$\tilde{\mathbf{X}}_+^g (\lambda) =  \mathbf{X}_- (\lambda)$. 

We need to check directly that $\mathbf{X}_- (\lambda)$, $\tilde{\mathbf{X}}_+ (\lambda)$,
$\mathbf{X}_-^g (\lambda)$, and $\tilde{\mathbf{X}}_+^g (\lambda)$ are frames for 
Lagrangian subspaces. The calculation is the same for each case, so we provide details
only for the first. If we compute $\mathbf{X}_- (\lambda)^* J \mathbf{X}_- (\lambda)$,
and use the orthogonality relation $R^* R = I$, we obtain a diagonal $2n \times 2n$
matrix with upper left $n \times n$ submatrix 
\begin{equation*}
    - D (\lambda)^3 + (D (\lambda)^*)^3 - (D (\lambda)^*)^2 D (\lambda) + D (\lambda)^* D (\lambda)^2
\end{equation*}
and lower right $n \times n$ submatrix
\begin{equation*}
    - (D (\lambda)^*)^3 + D (\lambda)^3 - D (\lambda)^2 D (\lambda)^* + D (\lambda) (D (\lambda)^*)^2.
\end{equation*}
The entries of $D (\lambda)^2$ are purely imaginary, so that 
$D (\lambda)^2 = - (D (\lambda)^*)^2$. If follows immediately that the two 
matrix expressions above are both 0, and we can conclude that 
$\mathbf{X}_- (\lambda)^* J \mathbf{X}_- (\lambda) = 0$. 

For the second part of Assumption {\bf (B2)}, we will verify that the matrix
$\mathbf{X}_- (\lambda)^* J \tilde{\mathbf{X}}_+ (\lambda)$ is 
non-singular for all $\lambda < \kappa$. Computing directly 
as with the calculation of $\mathbf{X}_- (\lambda)^* J \mathbf{X}_- (\lambda)$
just above, we find that 
\begin{equation*}
    \mathbf{X}_- (\lambda)^* J \tilde{\mathbf{X}}_+ (\lambda)
    = \begin{pmatrix}
    0 & 4 (D(\lambda)^*)^3 \\
    4 D(\lambda)^3 & 0
    \end{pmatrix},
\end{equation*}
and since the matrix $D(\lambda)$ is diagonal with non-zero 
entries (for $\lambda < \kappa$), we can conclude that 
$\mathbf{X}_- (\lambda)^* J \tilde{\mathbf{X}}_+ (\lambda)$
is non-singular.

In order to verify Assumption {\bf (B3)} in this case, we begin 
by observing that 
\begin{equation*}
    \mathbb{B}_{\lambda} (x; \lambda)
    = \begin{pmatrix}
    I & 0 & 0 & 0 \\
    0 & 0 & 0 & 0 \\
    0 & 0 & 0 & 0 \\
    0 & 0 & 0 & 0 
    \end{pmatrix},
\end{equation*}
so that 
\begin{equation*}
    \mathbf{X} (x; \lambda)^*  \mathbb{B}_{\lambda} (x; \lambda) \mathbf{X} (x; \lambda)
    = X_1 (x; \lambda)^* X_1 (x; \lambda),
\end{equation*}
where the $n \times 2n$ matrix $X_1 (x; \lambda)$ comprises the first $n$ rows
of each column of $\mathbf{X} (x; \lambda)$. We compute 
\begin{equation*}
    \int_{-\infty}^c \mathbf{X} (x; \lambda)^* \mathbb{B}_{\lambda} (x; \lambda) \mathbf{X} (x; \lambda) dx
    = \int_{-\infty}^c X_1 (x; \lambda)^* X_1 (x; \lambda) dx.
\end{equation*}
The columns of $X_1 (x; \lambda)$ are $2n$ linearly independent solutions 
of (\ref{fourth}), and so this matrix is positive definite by precisely 
the same considerations as discussed in Section \ref{sls-section}. 

We have now verified Assumptions {\bf (A)}, {\bf (B1)}, {\bf (B2)},
and {\bf (B3)} for this case, and so we can apply Theorem \ref{main-theorem}
to obtain the first claim in Theorem \ref{fourth-theorem}. For the 
second claim in Theorem \ref{fourth-theorem}, we will proceed as 
in the previous sections, using H\"ormander's index to replace 
$\tilde{\ell}_+ (\lambda_1)$ and $\tilde{\ell}_+ (\lambda_2)$
with a target frame $\mathbf{X}_T$ with respect to which the
calculations of the Maslov indices are monotonic. As discussed 
in \cite{Ho2020}, a natural frame to work with is 
\begin{equation*}
\mathbf{X}_T
= \begin{pmatrix}
0 & 0 \\
0 & I \\
I & 0 \\
0 & 0
\end{pmatrix}.
\end{equation*}
It is straightforward to check that $\mathbf{X}_T$ is the frame
for a Lagrangian subspace of $\mathbb{C}^{4n}$, and we denote this 
subspace $\ell_T$. 

Focusing on the case $\lambda = \lambda_2$, we recall from 
Section \ref{sls-exchange-section} that the difference 
\begin{equation*}
    \mas (\ell (\cdot; \lambda_2), \tilde{\ell}_+ (\lambda_2); (-\infty, +\infty])
    -  \mas (\ell (\cdot; \lambda_2), \ell_T; (-\infty, +\infty])
\end{equation*}
depends only on the fixed Lagrangian subspaces
$\ell_T$, $\tilde{\ell}_+ (\lambda_2)$, $\ell_- (\lambda_2)$ and $\ell_+ (\lambda_2)$,
and corresponds with H\"ormander's index
\begin{equation} \label{sls-hormander}
    s (\ell_T, \tilde{\ell}_+ (\lambda_2); \ell_- (\lambda_2), \ell_+ (\lambda_2)). 
\end{equation}
In order to evaluate H\"ormander's index, we'll again use the interpolation-space
approach of \cite{Ho2020}, and for this we need to work with a frame for 
$\tilde{\ell}_+ (\lambda_2)$ for which that analysis holds. For this, 
we introduce the inverse of the matrix 
\begin{equation*}
    \begin{pmatrix}
    R & R \\
    - R D (\lambda_2) & - R D(\lambda_2)^*
    \end{pmatrix},
\end{equation*}
which we see by inspection is 
\begin{equation} \label{fourth-tilde-M}
    \tilde{M} (\lambda_2) =
    \begin{pmatrix}
    (D (\lambda_2) - D(\lambda_2)^*)^{-1} & 0 \\
    0 & (D (\lambda_2) - D(\lambda_2)^*)^{-1}
    \end{pmatrix}
    \begin{pmatrix}
    -D(\lambda_2)^* R^* & - R^* \\
    D (\lambda_2) R^* & R^*
    \end{pmatrix}.
\end{equation}
We will replace the frame $\tilde{\mathbf{X}}_+ (\lambda_2)$ with the 
frame $\tilde{\mathbf{X}}_+ (\lambda_2) \tilde{M} (\lambda_2)$. For 
notational purposes, we can express this new frame as 
\begin{equation*}
   \tilde{\mathbf{X}}_+ (\lambda_2) \tilde{M} (\lambda_2)
   = \begin{pmatrix}
   I & 0 \\
   \tilde{X}_{21} & \tilde{X}_{22} \\
   \tilde{X}_{31} & \tilde{X}_{31} \\
   0 & I
   \end{pmatrix};
   \quad 
   \begin{pmatrix}
   \tilde{X}_{21} & \tilde{X}_{22} \\
   \tilde{X}_{31} & \tilde{X}_{31} 
   \end{pmatrix}
   = \begin{pmatrix}
   R D(\lambda_2)^2 & R (D(\lambda_2)^*)^2 \\
   - R D(\lambda_2)^3 & - R (D(\lambda_2)^*)^3 
   \end{pmatrix} \tilde{M} (\lambda_2).
\end{equation*}
In order to apply the development of \cite{Ho2020}, 
we need to check two conditions on the frames
$\mathbf{X}_T$ and $\tilde{\mathbf{X}}_+ (\lambda_2) \tilde{M} (\lambda_2)$.
First, we need to verify that 
\begin{equation} \label{first-condition}
    \mathbf{X}_T^* J (\tilde{\mathbf{X}}_+ (\lambda_2) \tilde{M} (\lambda_2))
    + (\tilde{\mathbf{X}}_+ (\lambda_2) \tilde{M} (\lambda_2))^* J \mathbf{X}_T = 0.
\end{equation}
To see this, we compute directly to find 
\begin{equation*}
    \mathbf{X}_T^* J (\tilde{\mathbf{X}}_+ (\lambda_2) \tilde{M} (\lambda_2))
    = \begin{pmatrix}
    I & 0 \\
    0 & -I
    \end{pmatrix}; 
    \quad (\tilde{\mathbf{X}}_+ (\lambda_2) \tilde{M} (\lambda_2))^* J \mathbf{X}_T
    = \begin{pmatrix}
    -I & 0 \\
    0 & I
    \end{pmatrix}, 
\end{equation*}
from which (\ref{first-condition}) is immediate. The second condition we need 
to check is that the matrix $\mathbf{X}_T^* J (\tilde{\mathbf{X}}_+ (\lambda_2) \tilde{M} (\lambda_2))$
is non-singular, and this is immediately clear from the previous calculations. 
We can conclude from \cite{Ho2020} that 
\begin{equation} \label{fourth-hormander}
\begin{aligned}
    s (\ell_T, \tilde{\ell}_+ (\lambda_2); \ell_- (\lambda_2), \ell_+ (\lambda_2))
    & = \mathcal{I} (\ell_+ (\lambda_2); \mathbf{X}_T, \tilde{\mathbf{X}}_+ (\lambda_2) \tilde{M} (\lambda_2)) \\
     &- \mathcal{I} (\ell_- (\lambda_2); \mathbf{X}_T, \tilde{\mathbf{X}}_+ (\lambda_2) \tilde{M} (\lambda_2)),
\end{aligned}
\end{equation}
where as in Section \ref{sls-section} the notation $\mathcal{I} (\cdot; \cdot, \cdot)$
has been taken directly from \cite{Ho2020}.
If $\lambda_2$ is not an eigenvalue for (\ref{fourth}), then $\ell_+ (\lambda_2)$ is the 
Lagrangian subspace with frame $\tilde{\mathbf{X}}_+^g (\lambda_2)$, which, as noted above,
is equal to $\mathbf{X}_- (\lambda_2)$. I.e., $\ell_+ (\lambda_2) = \ell_- (\lambda_2)$,
and so clearly H\"ormander's index is 0. We can conclude that if 
$\lambda_2$ is not an eigenvalue for (\ref{fourth}), then 
\begin{equation*}
    \mas (\ell (\cdot; \lambda_2), \tilde{\ell}_+ (\lambda_2); (-\infty, +\infty])
    =  \mas (\ell (\cdot; \lambda_2), \ell_T; (-\infty, +\infty]). 
\end{equation*}

Before turning to the case in which $\lambda_2$ is an eigenvalue 
of (\ref{fourth}), we observe that conjugate points
arising in the calculation of $\mas (\ell (\cdot; \lambda_2), \ell_T; (-\infty, + \infty])$
all have the same sign (negative). To see this, we employ Lemma 1.1 
of \cite{Ho2020}, which asserts (in the current setting) that in 
order to conclude monotonicity, we need to check two things: (1)
If $P_T$ denotes projection onto the Lagrangian subspace $\ell_T$, then
the matrix $P_T \mathbb{B} (x; \lambda_2) P_T$ is non-negative for a.e.
$x \in \mathbb{R}$; and (2) if $y(x; \lambda_2)$ is any non-trivial 
solution of (\ref{fourth-hammy}) with $y(x; \lambda_2) \in \ell_T$
for all $x$ in some interval $[a, b]$, $a < b$, then 
\begin{equation*}
    \int_a^b (\mathbb{B} (x; \lambda_2) y (x; \lambda_2), y(x; \lambda_2)) dx > 0.
\end{equation*}
For (1), we observe that 
\begin{equation*}
    v = \begin{pmatrix}
    v_1 \\
    v_2 \\
    v_3 \\
    v_4
    \end{pmatrix} \implies
    P_T v = \begin{pmatrix}
    0 \\
    v_2 \\
    v_3 \\
    0
    \end{pmatrix},
\end{equation*}
and consequently 
\begin{equation*}
    v^* P_T \mathbb{B} (x; \lambda_1) P_T v
    = (0 \,\,\, v_2^* \,\,\, v_3^* \,\,\, 0) 
    \begin{pmatrix}
    \lambda I - V(x) & 0 & 0 & 0 \\
    0 & I & 0 & 0 \\
    0 & 0 & 0 & -I \\
    0 & 0 & -I & 0 
    \end{pmatrix}
    \begin{pmatrix}
    0 \\
    v_2 \\
    v_3 \\
    0
    \end{pmatrix}
    = |v_2|^2 \ge 0.
   \end{equation*}
For (2), suppose $y(x; \lambda_2)$ is any non-trivial solution of (\ref{fourth-hammy})
so that $y(x; \lambda_2) \in \ell_T$ for all $x$ in some interval $[a, b]$, $a < b$.
Then, in particular, $\phi (x; \lambda_2) = 0$ for all such $x$, and 
since $\phi (x; \lambda_2)$ and its first two derivatives are absolutely 
continuous on $\mathbb{R}$ we can conclude that $\phi^{(k)} (x; \lambda_2) = 0$, 
$k = 1, 2, 3$, for a.e. $x \in (a, b)$.
But then $y(x; \lambda_2) = 0$ for a.e. $x \in (a, b)$, contradicting our 
assumption that $y(x; \lambda_2)$ is non-trivial. We conclude that 
Items (1) and (2) both hold, and from Lemma 1.1 of \cite{Ho2020}
we can conclude that conjugate points
arising in the calculation of $\mas (\ell (\cdot; \lambda_2), \ell_T; (-\infty, + \infty])$
all have the same sign (negative). If $\lambda_2$ is not an eigenvalue of 
(\ref{fourth}), we can now write 
\begin{equation} \label{fourth-kernel-sum}
    \mas (\ell (\cdot; \lambda_2), \ell_T; (-\infty, + \infty])
    = - \sum_{x \in \mathbb{R}} \dim (\ell (x; \lambda_2) \cap \ell_T).
\end{equation}
In writing this relation, we've taken advantage of the observations
that $\ell_- (\lambda_2) \cap \ell_T = \{0\}$ and that since 
we are currently assuming that $\lambda_2$ is not an eigenvalue 
of (\ref{fourth}), $\ell_+ (\lambda_2) \cap \ell_T = \{0\}$
(because $\ell_+ (\lambda_2) = \tilde{\ell}_+^g (\lambda_2)$ and 
$\tilde{\ell}_+^g (\lambda_2) \cap \ell_T = \{0\}$;
these claims are easily checked by using the frames for $\ell_- (\lambda_2)$
and $\tilde{\ell}_+^g (\lambda_2)$.) We conclude that the left-hand
side of (\ref{fourth-kernel-sum}) can be replaced by the Maslov index
$\mas (\ell (\cdot; \lambda_2), \ell_T; [-L, + L])$ for any $L$ 
sufficiently large, and correspondingly the right-hand side can 
be replaced by $- \sum_{x \in (-L, L)} \dim (\ell (x; \lambda_2) \cap \ell_T)$.
In addition, it's clear from 
monotonicity that the values of $x$ for which  
$\dim (\ell (x; \lambda_2) \cap \ell_T) \ne 0$ form a discrete set, 
so the right-hand side of (\ref{fourth-kernel-sum}) is a finite sum. 

We can use monotonicity precisely as in Section 
\ref{sls-exchange-section} to cover the case in which $\lambda_2$
is an eigenvalue of (\ref{fourth}). The same considerations
hold for $\lambda_1$, allowing us to write 
\begin{equation*}
    \begin{aligned}
        \mathcal{N} ([\lambda_1, \lambda_2))
        &= \sum_{x \in \mathbb{R}} \dim (\ell (x; \lambda_2) \cap \ell_T)
        - \sum_{x \in \mathbb{R}} \dim (\ell (x; \lambda_1) \cap \ell_T) \\
        &= \sum_{x \in \mathbb{R}} \dim \ker (\mathbf{X} (x; \lambda_2)^* J \mathbf{X}_T)
        - \sum_{x \in \mathbb{R}} \dim \ker (\mathbf{X} (x; \lambda_1)^* J \mathbf{X}_T), 
    \end{aligned}
\end{equation*}
where in obtaining this second equality, we have observed from Lemma 2.2 of 
\cite{HS2019} that if $\mathbf{X}_1$ and $\mathbf{X}_2$ are frames for any 
two Lagrangian subspaces $\ell_1$ and $\ell_2$ then 
\begin{equation*}
    \dim (\ell_1 \cap \ell_2) = \dim \ker (\mathbf{X}_1^* J \mathbf{X}_2).
\end{equation*}
For these latter calculations, if we write 
\begin{equation*}
    \mathbf{X} (x; \lambda)
    = \begin{pmatrix}
    X_{11} (x; \lambda) & X_{12} (x; \lambda) \\
     X_{21} (x; \lambda) & X_{22} (x; \lambda) \\
      X_{31} (x; \lambda) & X_{32} (x; \lambda) \\
       X_{41} (x; \lambda) & X_{42} (x; \lambda),
    \end{pmatrix}
\end{equation*}
then we have 
\begin{equation*}
    \mathbf{X} (x; \lambda_1)^* J \mathbf{X}_T 
    = \begin{pmatrix}
    X_{11} (x; \lambda) & X_{12} (x; \lambda) \\
    - X_{41} (x; \lambda) & - X_{42} (x; \lambda)
    \end{pmatrix}.
\end{equation*}
If we recall our specifications for the components of $y$ 
in terms of $\phi$, $\phi'$, $\phi''$, and $\phi'''$, we 
see that we can write 
\begin{equation} \label{fourth-kernel-sums}
     \mathcal{N} ([\lambda_1, \lambda_2))
     = \sum_{x \in \mathbb{R}} \dim \ker \Phi (x; \lambda_2)
        - \sum_{x \in \mathbb{R}} \dim \ker \Phi (x; \lambda_1),
\end{equation}
where (for $i = 1, 2$)
\begin{equation*}
    \Phi (x; \lambda_i)
    = \begin{pmatrix}
    \phi_1 (x; \lambda_i) & \phi_2 (x; \lambda_i) & \dots & \phi_{2n} (x; \lambda_i) \\
    \phi_1' (x; \lambda_i) & \phi_2' (x; \lambda_i) & \dots & \phi_{2n}' (x; \lambda_i)
    \end{pmatrix},
\end{equation*}
with $\{\phi_j (x; \lambda_1)\}_{j=1}^{2n}$ comprising a collection of $2n$ 
linearly independent solutions of (\ref{fourth}) that lie left in $\mathbb{R}$.  

Last, we check that we can take $\lambda_1$ sufficiently negative so that 
there are no conjugate points along the left shelf. For this, we begin by 
observing that $(s, \lambda_1) \in \mathbb{R} \times (-\infty, \kappa)$
will be a conjugate point for $\ell (\cdot; \lambda_1)$ and $\tilde{\ell}_+ (\lambda_1)$
if and only if $\lambda_1$ is an eigenvalue for 
\begin{equation*}
    \begin{aligned}
      \phi'''' + V (x) \phi &= \lambda \phi; \\
      \tilde{\mathbf{X}}_+ (\lambda_1)^* J 
      \begin{pmatrix}
      \phi (s) \\ \phi''(s) \\ - \phi''' (s) \\ - \phi'(s)
      \end{pmatrix} &= 0.
    \end{aligned}
\end{equation*}
Suppose $\lambda$ is an eigenvalue for this system, and let $\phi$
denote an associated eigenfunction. If we take an $L^2 ((-\infty, s),\mathbb{C}^n)$ 
inner product of the system with $\phi$, we obtain the integral relation 
\begin{equation} \label{fourth-integral-relation}
    \int_{-\infty}^s (\phi'''',\phi) dx
    + \int_{-\infty}^s (V \phi, \phi) dx 
    = \lambda \|\phi\|^2_{L^2 ((-\infty, s), \mathbb{C}^n)}.
\end{equation}
For the first integral, we integrate by parts twice to obtain 
the relation 
\begin{equation*}
\int_{-\infty}^s (\phi'''',\phi) dx
= \|\phi''\|^2_{L^2 ((-\infty, s), \mathbb{C}^n)}
+ \Big( {\phi (s) \choose - \phi' (s)}, {\phi''' (s) \choose \phi'' (s)} \Big).
\end{equation*}
Recalling (\ref{fourth-asym-frames}), we can express the boundary condition as 
\begin{equation*}
    \begin{pmatrix}
    R^* & (D(\lambda)^*)^2 R^* & - (D(\lambda)^*)^3 R^* & - D(\lambda)^* R^* \\
    R^* & D(\lambda)^2 R^* & - D(\lambda)^3 R^* & - D(\lambda) R^*
    \end{pmatrix}
    \begin{pmatrix}
      \phi''' (s) \\ \phi'(s) \\ \phi (s) \\ \phi''(s)
      \end{pmatrix}
      = 0,
\end{equation*}
or equivalently 
\begin{equation*}
   \begin{pmatrix}
    R^* & - D(\lambda)^* R^* \\
    R^* & - D(\lambda) R^* 
    \end{pmatrix} 
    {\phi'''(s) \choose \phi''(s)}
    = 
    \begin{pmatrix}
    (D(\lambda)^*)^3 R^* & (D(\lambda)^*)^2 R^* \\
    D(\lambda)^3 R^* & D(\lambda)^2 R^*
    \end{pmatrix}
    {\phi (s) \choose - \phi'(s)}.
\end{equation*}
Recalling (\ref{fourth-tilde-M}), we see that 
\begin{equation*}
 {\phi'''(s) \choose \phi''(s)}
 = \tilde{M} (\lambda)^* 
 \begin{pmatrix}
    (D(\lambda)^*)^3 R^* & (D(\lambda)^*)^2 R^* \\
    D(\lambda)^3 R^* & D(\lambda)^2 R^*
    \end{pmatrix}
    {\phi (s) \choose - \phi'(s)},
\end{equation*}
where we can write 
\begin{equation} \label{fourth-middle-three}
    \begin{aligned}
     \tilde{M} (\lambda)^* 
    &\begin{pmatrix}
    (D(\lambda)^*)^3 R^* & (D(\lambda)^*)^2 R^* \\
    D(\lambda)^3 R^* & D(\lambda)^2 R^*
    \end{pmatrix}
    = 
    \begin{pmatrix}
    R & 0 \\
    0 & R
    \end{pmatrix}
    \begin{pmatrix}
    - D(\lambda) & D(\lambda)^* \\
    -I & I
    \end{pmatrix} \\
    & \times \begin{pmatrix}
    (D(\lambda)^* - D(\lambda))^{-1} & 0 \\
    0 & (D(\lambda)^* - D(\lambda))^{-1}
    \end{pmatrix}
    \begin{pmatrix}
    (D(\lambda)^*)^3 & (D(\lambda)^*)^2 \\
    D(\lambda)^3 & D(\lambda)^2
    \end{pmatrix}
    \begin{pmatrix}
    R^* & 0 \\
    0 & R^*
    \end{pmatrix}
    \end{aligned}
\end{equation}
This matrix is clearly similar to the product of the middle
three matrices, and so has the same eigenvalues as that matrix. 
In order to compute these eigenvalues, we set 
\begin{equation*}
    \Lambda (\lambda) 
    := \diag (\sqrt[4]{\nu_1 - \lambda} \,\,\, \sqrt[4]{\nu_2 - \lambda} \,\,\, 
    \dots \,\,\, \sqrt[4]{\nu_n - \lambda}), 
\end{equation*}
so that 
\begin{equation*}
\begin{aligned}
    D(\lambda) &= (-\frac{1}{\sqrt{2}} - i \frac{1}{\sqrt{2}}) \Lambda (\lambda);
    \quad D(\lambda)^* - D(\lambda) = (i \sqrt{2}) \Lambda (\lambda); \\
    \quad D(\lambda)^2 &= i \Lambda(\lambda)^2;
    \quad D(\lambda)^3 = (\frac{1}{\sqrt{2}} - i \frac{1}{\sqrt{2}}) \Lambda (\lambda),
\end{aligned}
\end{equation*}
with corresponding adjoints. These relations allow us to express the product
of the middle three matrices in (\ref{fourth-middle-three}) as 
\begin{equation*}
    \begin{aligned}
        &\begin{pmatrix}
        (\frac{1}{\sqrt{2}} + i \frac{1}{\sqrt{2}}) \Lambda & (- \frac{1}{\sqrt{2}} + i \frac{1}{\sqrt{2}}) \Lambda \\
        -I & I
        \end{pmatrix}
        \begin{pmatrix}
        \frac{1}{i \sqrt{2}} \Lambda^{-1} & 0 \\
        0 & \frac{1}{i \sqrt{2}} \Lambda^{-1}
         \end{pmatrix}
         \begin{pmatrix}
         (\frac{1}{\sqrt{2}} + i \frac{1}{\sqrt{2}}) \Lambda^3 & -i \Lambda^2 \\
         (\frac{1}{\sqrt{2}} - i \frac{1}{\sqrt{2}}) \Lambda^3 & i \Lambda^2
         \end{pmatrix} \\
         & \quad = \begin{pmatrix}
         \sqrt{2} \Lambda^3 & - \Lambda^2 \\
         - \Lambda^2 & \sqrt{2} \Lambda
         \end{pmatrix}.
    \end{aligned}
\end{equation*}
We observe that this matrix is self-adjoint, and it follows
that the full matrix in (\ref{fourth-middle-three}) is self-adjoint.
In addition, we can compute
the eigenvalues of this matrix by computing the roots of 
the characteristic equation
\begin{equation*}
    \begin{aligned}
        \det &\begin{pmatrix}
        \sqrt{2} \Lambda^3 - \sigma I & - \Lambda^2 \\
        - \Lambda^2 & \sqrt{2} \Lambda - \sigma I
        \end{pmatrix}
        = \det ((\sqrt{2} \Lambda^3 - \sigma I) (\sqrt{2} \Lambda - \sigma I) - \Lambda^4) \\
        &= \det (\sigma^2 I - \sqrt{2} (\Lambda + \Lambda^3) \sigma + \Lambda^4).
    \end{aligned}
\end{equation*}
Here, since $\Lambda$ is a diagonal matrix, this determinant is a product
\begin{equation*}
    \prod_{j=1}^n (\sigma^2 - \sqrt{2} (\sqrt[4]{\nu_j - \lambda} + (\sqrt[4]{\nu_j - \lambda})^3) \sigma 
    + \nu_j - \lambda),
\end{equation*}
which clearly can have no roots for $\sigma \le 0$. We conclude that the matrix 
\begin{equation*}
    \tilde{M} (\lambda)^* 
 \begin{pmatrix}
    (D(\lambda)^*)^3 R^* & (D(\lambda)^*)^2 R^* \\
    D(\lambda)^3 R^* & D(\lambda)^2 R^*
    \end{pmatrix}
\end{equation*}
is positive definite, and so 
\begin{equation*}
    \Big( {\phi (s) \choose - \phi' (s)}, {\phi''' (s) \choose \phi'' (s)} \Big)
    = \Big( {\phi (s) \choose - \phi' (s)},  \tilde{M} (\lambda)^* 
 \begin{pmatrix}
    (D(\lambda)^*)^3 R^* & (D(\lambda)^*)^2 R^* \\
    D(\lambda)^3 R^* & D(\lambda)^2 R^*
    \end{pmatrix} {\phi (s) \choose - \phi' (s)} \Big) \ge 0
\end{equation*}
for all $s \in \mathbb{R}$. Returning to 
(\ref{fourth-integral-relation}), we see that 
\begin{equation*}
     \lambda \|\phi\|^2_{L^2 ((-\infty, s), \mathbb{C}^n)}
    \ge \int_{-\infty}^s (V \phi, \phi) dx  
    \ge - \|V\|_{L^{\infty} (\mathbb{R}, \mathbb{C}^{n \times n})} \|\phi\|^2_{L^2 ((-\infty, s), \mathbb{C}^n)},
\end{equation*}
and consequently 
\begin{equation*}
    \lambda \ge - \|V\|_{L^{\infty} (\mathbb{R}, \mathbb{C}^{n \times n})}.
\end{equation*}
In this way, we see that if we take $\lambda_1 < - \|V\|_{L^{\infty} (\mathbb{R}, \mathbb{C}^{n \times n})}$
then there will be no conjugate points along the vertical shelf at $\lambda = \lambda_1$. This gives
the final claim in Theorem \ref{fourth-theorem}. 
\hfill $\square$


\begin{thebibliography}{99}

\bibitem{ALMS1994} F. V. Atkinson, H. Langer, R. Mennicken, and 
A.A. Shkalikov, {\it The essential spectrum of some matrix operators},
Math. Nachr. {\bf 167} (1994) 5 -- 20.

\bibitem{Arnold1967} V. I. Arnol'd, {\it Characteristic class entering in quantization 
conditions}, Func. Anal. Appl. {\bf 1} (1967) 1 -- 14.

\bibitem{BCCJM2020} T. J. Baird, P. Cornwell, G. Cox, C. K. R. T. Jones, 
and R. Marangell, {\em A Maslov index for non-Hamiltonian systems}.
Preprint 2020, arXiv: 2006.14517v1. 

\bibitem{BCJLMS2018} M. Beck, G. Cox, C. K. R. T. Jones, 
Y. Latushkin, K. McQuighan, and A. Sukhtayev, {\it Instability of pulses 
in gradient reaction-diffusion systems: a symplectic approach},
Philos. Trans. Roy. Soc. A {\bf 376} (2018), no. 2117, 20170187, 
20 pp. 

\bibitem{BF1998} B.\ Booss-Bavnbek and K.\ Furutani, {\em The Maslov index: 
a functional analytical definition and the spectral flow formula,}
Tokyo J.\ Math.\ {\bf 21} (1998), 1--34.

\bibitem{BJ1995} A.\ Bose and C.\ K.\ R.\ T.\ Jones,
{\em Stability of the in-phase traveling wave solution in a pair 
of coupled nerve fibers,}
Indiana U. Math. J. {\bf 44} (1995) 189 -- 220.

\bibitem{BM2013} M.\ Beck and S.\ Malham, 
{\em Computing the Maslov index for large systems}, 
Proc. Amer. Math. Soc. {\bf 143} (2015), no. 5, 2159–2173.

\bibitem{Bott1956} R. Bott, {\em On the iteration of closed geodesics and the Sturm 
intersection theory}, Comm. Pure Appl. Math. {\bf 9} (1956) 171--206.

\bibitem{CB2015} F.\ Chardard and T.\ J.\ Bridges,
{\em Transversality of homoclinic orbits, the Maslov index, and 
the symplectic Evans function}, Nonlinearity {\bf 28} (2015)
77--102. 

\bibitem{CDB2009} F.\ Chardard, F.\ Dias and T.\ J.\ Bridges, 
{\em Computing the Maslov index of solitary waves, Part 1: 
Hamiltonian systems on a four-dimensional phase space}, 
Phys. D {\bf 238} (2009) 1841 -- 1867.

\bibitem{CDB2011} F.\ Chardard, F.\ Dias and T.\ J.\ Bridges, 
{\em Computing the Maslov index of solitary waves, Part 2: 
Phase space with dimension greater than four}, 
Phys. D {\bf 240} (2011) 1334 -- 1344.

\bibitem{CH2007} C--N. Chen and Xijun Hu, {\em Maslov index for homoclinic 
orbits of Hamiltonian systems}, Ann. Inst. H. Poincar\'e Anal. Nonlin\'eaire
{\bf 24} (2007) 589--603.

\bibitem{CH2014} C--N. Chen and Xijun Hu, {\em Stability analysis for standing
pulse solutions to FitzHugh--Nagumo equations}, Calculus of Variations
and Partial Differential Equations {\bf 49} (2014) 827--845.

\bibitem{Chardard2009} F.\ Chardard, 
{\em Stability of Solitary Waves}, Doctoral thesis, Centre de Mathematiques et de 
Leurs Applications, 2009. Advisor: T.\ J.\ Bridges. 

\bibitem{CJ2018} P. Cornwell and C. K. R. T. Jones,
{\em On the existence and stability of fast traveling waves in a doubly diffusive FitzHugh--Nagumo system},
SIAM Journal on Applied Dynamical Systems {\bf 17} (2018) 754-787

\bibitem{CJ2020} P. Cornwell and C. K. R. T. Jones,
{\em A stability index for traveling waves in activator-inhibitor systems}
Proceedings of the Royal Society of Edinburgh: Section A Mathematics 
{\bf 150} (2020) 517--548.

\bibitem{CJLS2016} G.\ Cox, C.\ K.\ R.\ T.\ Jones, Y.\ Latushkiun, and A.\ Sukhtayev,
{\em The Morse and Maslov indices for multidimensional Schr\"odinger operators
with matrix-valued potentials}, Trans. Amer. Math. Soc. {\bf 368} (2016) 8145-8207. 

\bibitem{CLM1994}  S.\ Cappell, R.\ Lee and E.\ Miller, {\em On the Maslov index},
Comm.\ Pure Appl.\ Math.\ {\bf 47}  (1994), 121--186.

\bibitem{Co2019} P. Cornwell, 
{\em Opening the Maslov box for traveling waves in skew-gradient systems: counting eigenvalues and proving (in) stability},
Indiana U. Math. J. {\bf 68} (2019) 1801-1832.

\bibitem{DJ2011} J.\ Deng and C.\ Jones, {\em Multi-dimensional Morse Index Theorems 
and a symplectic view of elliptic boundary value problems,}
Trans. Amer. Math. Soc. {\bf 363} (2011) 1487 -- 1508.

\bibitem{Du1976} J. J. Duistermaat, {\em On the Morse index in variational 
calculus}, Adv. Math. {\bf 21} (1976) 173--195.

\bibitem{F2004} K.\ Furutani, {\em Fredholm-Lagrangian-Grassmannian and the Maslov index}, 
Journal of Geometry and Physics {\bf 51} (2004) 269 -- 331.

\bibitem{He1981} D. Henry, 
{\em Geometric theory of semilinear parabolic equations}, 
Lect. Notes Math. \textbf{840}, Springer-Verlag, Berlin-New York, 1981.

\bibitem{HJK2018} P. Howard, S. Jung, and B. Kwon,
{\it The Maslov index and spectral counts for Hamiltonian systems 
on [0,1]}, J. Dynamics and Differential Equations {\bf 30} (2018)
1703--1729.

\bibitem{HLS2017} P. Howard, Y. Latushkin, and A. Sukhtayev,
{\it The Maslov index for Lagrangian pairs on $\mathbb{R}^{2n}$}, 
J. Mathematical Analysis and Applications {\bf 451} (2017) 794-821.

\bibitem{HLS2018} P. Howard, Y. Latushkin, and A. Sukhtayev,
{\it  The Maslov and Morse indices for Schr\"odinger operators on $\mathbb{R}$}, 
Indiana U. Mathematics Journal {\bf 67} (2018) 1765-1815.

\bibitem{Ho2020} P. Howard, {\em H\"ormander's index and Oscillation Theory},
Preprint 2020. 

\bibitem{HS2016} P. Howard and A. Sukhtayev, 
{\em The Maslov and Morse indices for Schr\"odinger operators on $[0,1]$}, 
J. Differential Equations {\bf 260} (2016), no. 5, 4499-4549.

\bibitem{HS2019} P. Howard and A. Sukhtayev, {\it Renormalized oscillation theory
for linear Hamiltonian systems on $[0, 1]$ via the Maslov index},
Preprint 2019, arXiv 1808.08264.

\bibitem{HS2020} P. Howard and A. Sukhtayev, 
{\em The Maslov and Morse Indices 
for Sturm-Liouville Systems on the Half-Line}, 
Discrete and Continuous Dynamical Systems A {\bf 40(2)} (2020)
983--1012. 

\bibitem{HS2020b} P. Howard and A. Sukhtayev, 
{\em Renormalized oscillation theory for singular linear Hamiltonian 
systems}, Preprint 2020, arXiv 2009.10681. 

\bibitem{J1988a}  C.\ K.\ R.\ T.\ Jones,
{\em Instability of standing waves for nonlinear Schr\"odinger-type equations},
Ergodic Theory Dynam. Systems {\bf 8} (1988) 119 -- 138.

\bibitem{J1988b}  C.\ K.\ R.\ T.\ Jones,
{\em An instability mechanism for radially symmetric standing
waves of a nonlinear Schr\"odinger equation},
J. Differential Equations {\bf 71} (1988) 34 -- 62.

\bibitem{JM2012} C.\ K.\ R.\ T.\ Jones and R.\ Marangell,
{\em The spectrum of travelling wave solutions to the Sine-Gordon
equation}, Discrete and Cont. Dyn. Sys. {\bf 5} (2012) 925 -- 937.

\bibitem{LS2018} Y. Latushkin and S. Sukhtaiev,
{\em The Maslov Index and the spectra of second order elliptic operators}, 
Advances in Mathematics {\bf 329} (2018) 422–-486.

\bibitem{KP2013} T. Kapitula and K. Promislow, 
{\em Spectral and dynamical stability of nonlinear waves}, 
Springer, New York, 2013.

\bibitem{Maslov1965a} V.\ P.\ Maslov, {\it Theory of perturbations and 
asymptotic methods}, Izdat. Moskov. Gos. Univ. Moscow, 1965. French 
tranlation Dunod, Paris, 1972. 

\bibitem{Morse1934} H. C. M. Morse, {\it The calculus of variations in the large},
AMS Coll. Publ. {\bf 18} (1934).

\bibitem{Ph1996} J.\ Phillips, {\em Selfadjoint Fredholm operators and spectral flow}, 
Canad.\ Math.\ Bull. {\bf 39} (1996), 460--467.

\bibitem{RS1993} J.\ Robbin and D.\ Salamon,
{\em The Maslov index for paths}, Topology {\bf 32} (1993) 827 -- 844.

%\bibitem{Sm1965} S. Smale, {\em On the Morse index theorem}, 
%J. Math. Mech. 14 (1965), 1049–1055.

\bibitem{Sturm1836} C. Sturm, {\it M\'emoire sur les \'equations diff\'erentielles lin\'eaires du second ordre},
J. math. pures appl. {\bf 1} (1836) 106-186.

\bibitem{We1987} J. Weidmann, {\it Spectral theory of ordinary differential 
operators}, Springer-Verlag 1987.

\bibitem{ZWZ2018} Y. Zhou, L. Wu, and C. Zhu, 
{\it H\"ormander index in the finite-dimensional case},
Math. China {\bf 13} (2018) 725--761. 

\end{thebibliography}
\end{document}